\documentclass[a4paper]{amsbook}
\usepackage{amssymb}
\usepackage{amsmath}
\usepackage{amsfonts}
\usepackage{amsthm}
\usepackage{amscd}
\usepackage{times}
\usepackage{txfonts}

\theoremstyle{plain}
\newtheorem{thm}{Theorem}[subsection]
\newtheorem{lem}[thm]{Lemma}
\newtheorem{cor}[thm]{Corollary}
\newtheorem{prop}[thm]{Proposition}
\theoremstyle{definition}
\newtheorem{defn}[thm]{Definition}

\newtheorem{rem}[thm]{Remark}

\newcommand{\bZ}{{\mathbb Z}}
\newcommand{\bQ}{{\mathbb Q}}
\newcommand{\bR}{{\mathbb R}}
\newcommand{\bI}{{\mathbb I}}
\newcommand{\bJ}{{\mathbb J}}

\newcommand{\bH}{{\mathbb H}}
\newcommand{\fm}{{\mathfrak m}}
\newcommand{\fD}{{\mathfrak D}}
\newcommand{\fR}{{\mathfrak R}}
\newcommand{\fB}{{\mathfrak B}}

\newcommand{\on}[1]{\mathop{\mathrm{#1}}\nolimits}

\renewcommand{\mid}{\hskip.02in;\hskip.02in}
\renewcommand{\bibname}{{\sc References}}

\renewcommand{\thesection}{\S\thechapter.%
\arabic{section}}
\renewcommand{\thesubsection}{\thechapter.%
\arabic{section}.\arabic{subsection}}
\renewcommand{\thesubsubsection}{\thechapter.%
\arabic{section}.\arabic{subsection}.\arabic{subsubsection}}

\makeatletter
\newcounter{step}
\newenvironment{step}[1]{
\refstepcounter{step}\par
{\bf Step}\ $\mathbf\thestep$.\quad{\it #1}
\par\noindent}{}
\newcounter{case}
\newenvironment{case}[1]{
\refstepcounter{case}\par
\underline{Case\ $\thecase$}.\quad{\it #1}
\par\noindent}{}
\renewenvironment{proof}[1][\proofname]{\par
\normalfont \topsep6\p@\@plus6\p@\relax
\trivlist\itemindent\normalparindent
\item[\hskip\labelsep\scshape 
#1\@addpunct{.}]\ignorespaces
\setcounter{step}0
\setcounter{case}0
}\@endtheorem
\makeatother

\newcounter{komoku}[subsection] 
\renewcommand{\thekomoku}{\thechapter.\arabic{section}.%
\arabic{subsection}.\arabic{komoku}} 
\newcommand{\komoku}[1]{\refstepcounter{komoku}%
\hskip0in\par\noindent\thekomoku\ {\bf#1}\ } 
\newcounter{subkomoku}[komoku]

\addtocontents{toc}{{\noindent\bf Part II.\quad
Basic invariants associated to the idealistic filtration
}}
\makeatletter
\renewcommand{\tocsubsection}[3]{%
\phantom{\S}{\@ifnotempty{#2}{\ignorespaces#1 #2.\quad}}#3}
\makeatother

\makeatletter
\newenvironment{bibpart}[1]{%
\normalfont\footnotesize\labelsep .5em\relax
\renewcommand\theenumiv{\arabic{enumiv}}\let\p@enumiv\@empty
\list{\@biblabel{\theenumiv}}{\settowidth\labelwidth{\@biblabel{#1}}%
\leftmargin\labelwidth \advance\leftmargin\labelsep
\usecounter{enumiv}}%
\sloppy \clubpenalty\@M \widowpenalty\clubpenalty
\sfcode`\.=\@m
}{%
\def\@noitemerr{\@latex@warning{Empty `thebibliography' environment}}%
\endlist
}
\makeatother

\begin{document}
\title{
{\huge\sc Toward resolution of singularities
\linebreak
over a field of \linebreak
positive characteristic \linebreak
{\ }\linebreak
(The Idealistic Filtration Program)}
\vskip.6in 
{\it Dedicated to Professor Heisuke Hironaka}
\vskip.8in 
{Part\ II.
\vskip0in
Basic invariants associated to the idealistic filtration and their
properties}
\vskip2.8in 
}
\author{\sc Hiraku\ Kawanoue \& \sc Kenji Matsuki}
\date{today}
\maketitle

\tableofcontents
\begin{chapter}{Introduction to Part II}
\begin{section}{Overview of the series.}\label{0.1}
This is the second of the series of papers under the title

``Toward resolution of singularities  over a field of positive
characteristic

(the Idealistic Filtration Program)''
\begin{center} 
\begin{tabular}{ll}
Part I. &Foundation; the language of the idealistic filtration
\\ 
Part II. &Basic invariants associated to the idealistic filtration
\\ 
&and their properties
\\ 
Part III. &Transformations and modifications of the idealistic filtration
\\ 
Part IV.  &Algorithm in the framework of the idealistic filtration.
\end{tabular} 
\end{center} 
For a brief summary of the entire series, including its goal and the
overview of the Idealistic Filtration Program, we refer
the reader to the introduction in Part I.

Here we will concentrate ourselves on the outline of Part II, which is
presented in the next section.
\end{section} 

\begin{section}{Outline of Part II.}\label{0.2}
As described in the overview of the Idealistic Filtration Program 
(cf. \S 0.2 in Part I), we construct a strand of invariants, whose
maximum locus determines each center of blowup of our algorithm for
resolution of singularities.  The strand of invariants consists of
the units (cf. 0.2.3.2.2 in Part I), each  of which is a triplet of numbers
$(\sigma, \widetilde{\mu},s)$ associated to a certain
idealistic filtration (cf. Chapter 2 in Part I) and a simple normal crossing
divisor $E$ called a boundary.  (To be precise, the invariant 
$\sigma$ is a sequence of numbers indexed by ${\mathbb Z}_{\geq 0}$ 
as described in Definition 3.2.1.1 in Part I.)  The purpose of Part 
II is to establish the fundamental properties of the invariants 
$\sigma$ and
$\widetilde{\mu}$.  They are the main constituents of the unit, while the
remaining factor $s$ can easily be computed as the number of
(certain specified) components in the boundary passing through a given
point, and needs no further mathematical discussion.  Our goal is to
study the intrinsic nature of these invariants associated to a given
idealistic filtration.  The discussion in Part II does not involve the
analysis regarding the
exceptional divisors created by blowups, and hence could only be directly
applied to the situation {\it at year 0} of our algorithm.  The systematic
discussion
on how some subtle adjustments should be made in the presence of the
exceptional divisors {\it after year 0} and on how the strand of
invariants
functions in the algorithm, built upon the analysis in Part
III of the modifications and transformations of an idealistic filtration,
will have to wait for Part IV.

In the appendix, we report a new development, unexpected at the time of
writing Part I, which suggests a possibility of constructing an algorithm
using only the $\fD$-saturation (or $\fD_E$-saturation) and without 
using the $\fR$-saturation,
still within the framework of the Idealistic Filtration Program.  This would
avoid the
problem of termination, which we specified in the introduction to Part I as
the only missing piece toward completing our algorithm in positive
characteristic.  (See \S 0.3. $\mathbf{Current\ status\ of\ the\ 
Idealistic\ Filtration\ Program}$ at the end of the introduction here 
in Part II for further developments and ``evolution'' of IFP up to 
date.)

\vskip.1in

The following is a rough description of the content of each chapter and the
appendix in Part II.

Throughout the description, let $R$ be the
coordinate ring of an affine open subset of a nonsingular variety $W$ of
dimension  $d = \dim W$ over an algebraically closed field $k$ of
positive characteristic
$\on{char}(k) = p > 0$ or of characteristic zero $\on{char}(k) = 0$, where
in the latter case we set $p = \infty$ formally (cf. 0.2.3.2.1
and Definition 3.1.1.1 (2) in Part I).

\begin{subsection}
{Invariant $\sigma$.}\label{0.2.1} Chapter 1 is devoted to the discussion of
the invariant $\sigma$, which is defined for a ${\fD}$-saturated idealistic
filtration
${\bI}$ over $R$ (cf. 2.1.2
in Part I).  The subtle adjustment of the invariant $\sigma$, in the
presence of the exceptional divisor $E$, which is defined for a
${\fD}_E$-saturated idealistic
filtration (cf. 1.2.2. Logarithmic differential operators in Part I), will
be postponed until Part III and Part IV.

\komoku{Leading algebra and its structure.}\label{0.2.1.1} We fix a closed
point $P \in \on{Spec}\hskip.02in R \subset W$, with ${\fm}_P$ denoting the
maximal
ideal for the local ring $R_P$.  The leading algebra $L({\bI}_P)$ of the
localization ${\bI}_P$ of the idealistic filtration ${\bI}$ at $P$ is
defined to be the
graded $k$-subalgebra of $G_P = \bigoplus_{n \in {\bZ}_{\geq 0}} (G_P)_n =
\bigoplus_{n \in {\bZ}_{\geq 0}}{\fm}_P^n/{\fm}_P^{n+1}$ (cf. 3.1.1 in Part
I)
$$L({\bI}_P) = \bigoplus_{n \in {\bZ}_{\geq 0}}L({\bI}_P)_n \subset G_P,$$
where
$$L({\bI}_P)_n = \left\{\overline{f} = (f \bmod{{\fm}_P^{n+1}}
) \mid (f,n) \in {\bI}_P, f \in {\fm}_P^n\right\}.$$
For $e \in {\bZ}_{\geq 0}$ with $p^e \in {\bZ}_{> 0}$, we define the pure
part $L({\bI}_P)_{p^e}^{\on{pure}}$ of $L({\bI}_P)_{p^e}$ by the formula
$$L({\bI}_P)_{p^e}^{\on{pure}} = L({\bI}_P)_{p^e} \cap F^e((G_P)_1) \subset
L({\bI}_P)_{p^e}$$
where $F^e$ is the $e$-th power of the Frobenius map of $G_P$.

The most remarkable structure of the leading algebra $L({\bI}_P)$ is that it
is generated by its pure part (cf. Lemma 3.1.2.1 in Part I), i.e.,
$$L({\bI}_P) = k[L({\bI}_P)^{\on{pure}}] \hskip.1in \on{where} \hskip.1in
L({\bI}_P)^{\on{pure}} = \bigsqcup_{e \in {\bZ}_{\geq
0}}L({\bI}_P)_{p^e}^{\on{pure}}.$$
This follows from the fact that ${\bI}_P$ is ${\fD}$-saturated, since so is
${\bI}$ (cf. compatibility of ${\fD}$-saturation with localization,
discussed in \S
2.4 in Part I).

\komoku{Definition of the invariant $\sigma$ and its
computation.}\label{0.2.1.2} We define the invariant $\sigma(P)$ by the
formula
$$\sigma(P) = \left(d - l_{p^e}^{\on{pure}}(P)\right)_{e 
\in {\bZ}_{\geq 0}} \in
\prod_{e \in {\bZ}_{\geq 0}}{\bZ}_{\geq 0} \hskip.1in \on{where} \hskip.1in
l_{p^e}^{\on{pure}}(P) = \dim_kL({\bI}_P)_{p^e}^{\on{pure}},$$
which reflects the behavior of the pure part of the leading algebra
$L({\bI}_P)$.  Varying $P$ among all the closed points
${\fm}\text{-}\on{Spec}\hskip.02in
R$ (i.e., all the maximal ideals of $R$), we obtain the invariant
$$\sigma:{\fm}\text{-}\on{Spec}\hskip.02in R \rightarrow \prod_{e \in
{\bZ}_{\geq 0}}{\bZ}_{\geq 0}.$$
Recall that Lemma 3.1.2.1 in Part I gives the description of a specific set
of generators for the leading algebra $L({\bI}_P)$ taken from its pure part.
Using
this lemma, we can compute the dimension of the pure part
$l_{p^e}^{\on{pure}}(P) =
\dim_kL({\bI}_P)_{p^e}^{\on{pure}}$ in terms of the dimension of the entire
degree $p^e$ component $l_{p^e}(P) = \dim_kL({\bI}_P)_{p^e}$ and in terms of
the
dimensions of the pure parts $l_{p^{\alpha}}^{\on{pure}}(P)$ for
${\alpha} = 0,
\ldots, e-1$.  That is to say, $l_{p^e}^{\on{pure}}(P)$ can be computed
inductively from $l_{p^e}(P)$ and the dimensions of the pure parts of lower
degree.

\komoku{Upper semi-continuity of the invariant $\sigma$.}\label{0.2.1.3} We
observe that
$l_{p^e}(P)$ can be computed as the rank of a certain ``Jacobian-like''
matrix, and hence is easily seen to be lower semi-continuous as a function
of $P$.  The
upper semi-continuity of the invariant $\sigma = \left(d -
l_{p^e}^{\on{pure}}\right)_{e
\in {\bZ}_{\geq 0}}$ then follows immediately from the inductive computation
of the pure part $l_{p^e}^{\on{pure}}$ described in \ref{0.2.1.2}.  The
upper
semi-continuity of the invariant
$\sigma$ as a function over
${\fm}\text{-}{\on{Spec}}\ R$ also allows us to extend its domain to
${\on{Spec}}\hskip.02in R$.  That is to say, we have the invariant $\sigma$
defined
over the extended domain
$$\sigma:\on{Spec}\hskip.02in R \rightarrow \prod_{e
\in {\bZ}_{\geq 0}}{\bZ}_{\geq 0},$$
which is automatically upper semi-continuous as a function over
${\on{Spec}}\hskip.02in R$.

\komoku{Clarification of the meaning of the upper
semi-continuity.}\label{0.2.1.4} We say by definition that a function $f:X
\rightarrow T$, from a topological
space $X$ to a totally ordered set $T$, is upper semi-continuous if the set
$X_{\geq t} = \{x \in X \mid f(x) \geq t\}$ is closed for any $t \in T$.
When the
target space $T$ is not well-ordered, however, we have to be extra-careful
as we try to see the equivalence of this definition to the other
``well-known''
conditions for the upper semi-continuity.  The target space of the invariant
$\sigma:{\fm}\text{-}\on{Spec}\hskip.02in R
\rightarrow
\prod_{e
\in {\bZ}_{\geq 0}}{\bZ}_{\geq 0}$ is a priori not well-ordered.
Nevertheless, using the fact that $l_{p^e}^{\on{pure}}(P)$ is non-decreasing
as a function of $e
\in {\bZ}_{\geq 0}$ for a fixed $P \in {\fm}\text{-}\on{Spec}\hskip.02in R$,
we observe that the target space for the invariant $\sigma$ can be replaced
by some
well-ordered subset.  It can be seen easily then that the upper
semi-continuity of the invariant $\sigma$ in
the above sense is actually equivalent to the condition that, given a point
$P
\in
{\fm}\text{-}\on{Spec}\hskip.02in R$, there exists a neighborhood
$U_P$ of
$P$ such that
$\sigma(P)
\geq
\sigma(Q)$ for any point $Q \in U_P$.  From this upper semi-continuity,
interpreted in the equivalent condition, it follows that the domain of the
invariant
$\sigma$ can be extended from ${\fm}\text{-}\on{Spec}\hskip.02in R$ to
$\on{Spec}\hskip.02in R$, as mentioned at the end of \ref{0.2.1.3}.  We
summarize the basic
facts surrounding the definition of the upper semi-continuity in Chapter 1
for the sake of clarification.

\komoku{Local behavior of a leading generator system\footnote{We use the 
abbreviation ``LGS'' for the word ``Leading Generator System''. 
Prof. Cossart kindly suggested to us that ``LGS'' could be read 
``Leading {\it Giraud}  System'' in honor of J. Giraud, whose 
contribution (cf. \cite{MR0460712},\cite{MR0384799}) is profound 
in search of the right 
notion of ``a hypersurface of maximal contact'' in positive 
characteristic.}, and its modification

into one which is uniformly pure.}\label{0.2.1.5} Recall that a subset
${\bH} =
\left\{(h_l,p^{e_l})\right\}_{l = 1}^N \subset {\bI}_P$ 
with associated nonnegative
integers $0 \leq e_1 \leq \cdots \leq e_N$ is said to be a leading generator
system of
the idealistic filtration ${\bI}_P$, if the leading terms of its elements
provide a specific set of generators for the leading
algebra
$L({\bI}_P)$.  More precisely, it satisfies the following conditions (cf.
3.1.3 in Part I):
\begin{center} 
\begin{tabular}{ll}
(i) & $h_l \in {\fm}_P^{p^{e_l}}$ and $\overline{h_l} = (h_l 
\bmod {\fm}_P^{p^{e_l}+1}) \in
L({\bI}_P)_{p^{e_l}}^{\on{pure}}$ for
$l = 1, \ldots, N$,
\\ 
(ii) & $\left\{\overline{h_l}^{p^{e - e_l}} \mid e_l \leq e\right\}$ 
consists of $\#\{l
\mid e_l \leq e\}$-distinct elements, and forms a $k$-basis of \\
&$L({\bI}_P)_{p^e}^{\on{pure}}$ for any $e \in {\bZ}_{\geq 0}$.
\\ 
\end{tabular} 
\end{center} 

Since the leading algebra $L({\bI}_P)$ is generated by its pure part
$$L({\bI}_P)^{\on{pure}} = \bigsqcup_{e \in {\bZ}_{\geq
0}}L({\bI}_P)_{p^e}^{\on{pure}}\qquad (\text{cf.\ } \ref{0.2.1.1}), $$
we conclude from
condition (ii) that the leading terms $\left\{\overline{h_l} = (h_l 
\bmod{{\fm}_P^{p^{e_l} + 1})}\right\}_{l = 1}^N$ of ${\bH}$ provide a set of
generators for $L({\bI}_P)$, 
i.e.,$L({\bI}_P) = k[\left\{\overline{h_l}\right\}_{l = 1}^N]$.

A basic question then about the local behavior of a leading generator system
is: 

Does ${\bH}$ remain being a leading generator system of ${\bI}_Q$ for any
closed point $Q$ in a
neighborhood
$U_P$ of
$P$ (if we take
$U_P$ small enough) ?

Even though the answer is {\it no} in general, we show that we can modify
a given leading generator system ${\bH}$ into a new one ${\bH}'$ such that
${\bH}'$
is a leading generator system of ${\bI}_Q$ for any closed point $Q$ in a
neighborhood $U_P$, as long as $Q$ is on the local maximum locus of the
invariant
$\sigma$ (and $Q$ is on the support of ${\bI}$).  The last extra condition
is equivalent to saying $\sigma(Q) = \sigma(P)$ by the upper semi-continuity
of the
invariant
$\sigma$.  We say then ${\bH}'$ is ``{\it uniformly pure}''.  
We will use this modification as the
main tool to derive the upper semi-continuity of the invariant
$\widetilde{\mu}$ in
Chapter 3.

\end{subsection}

\begin{subsection}{Power series expansion.}\label{0.2.2} Chapter 2 is
devoted to the discussion of the power series expansion with respect to a
leading generator
system and its (weakly-)associated regular system of parameters.

\komoku{Similarities between a regular system of parameters and a leading
generator system.}\label{0.2.2.1} If we have a leading generator system
${\bH} =
\left\{(h_l,p^{e_l})\right\}_{l = 1}^N$ {\it in characteristic zero} (for a
${\fD}$-saturated idealistic filtration ${\bI}_P$ over $R_P$ at a closed
point $P \in
\on{Spec}\hskip.02in R
\subset W$), then the elements in the leading generator system are all
concentrated at level 1, i.e., $e_l = 0$ and $p^{e_l} = 1$ for $l
= 1,
\ldots, N$ (cf. Chapter 3 in Part I).  This implies by definition of a
leading generator system that the set of the elements
$H = (h_1, \ldots, h_l)$ forms (a part of) a regular system of parameters
$(x_1, \ldots, x_d)$.  (Say, $h_l = x_l$ for $l = 1, \ldots, N$.)  
{\it In positive characteristic},
this is no longer the case.  However, we can still regard the notion of a
leading generator system as a generalization of the notion of a regular
system of
parameters, and we may expect some similarities between the two notions.
One of such expected similarities is the power series expansion, which we
discuss next. 

\komoku{Power series expansion with respect to a leading generator
system.}\label{0.2.2.2} In
characteristic zero, any
element $f
\in R_P$ has a power series expansion (with respect to the regular system of
parameters $X = (x_1, \ldots, x_d)$, where $h_l = x_l$ for $l = 1, \ldots,
N$, with $H
= (h_1, \ldots, h_N)$ consisting of the elements of a leading generator
system as described in \ref{0.2.2.1})
$$f = \sum_{I \in \left({\bZ}_{\geq 0}\right)^d}c_IX^I = \sum_{B \in
\left({\bZ}_{\geq 0}\right)^N} a_BH^B$$
where $c_I \in k$ and where $a_B$ is a power series in terms of the
remainder $(x_{N+1}, \ldots, x_d)$ of the regular system of
parameters.  

In positive characteristic, we expect to have a power series expansion with
respect to a leading generator system.  More specifically and more
generally, the setting for Chapter 2 is given as follows.  We have a subset
${\mathcal H} = \{h_1, \ldots, h_N\} \subset R_P$ consisting of $N$
elements, and
nonnegative integers $0
\leq e_1
\leq
\cdots \leq e_N$ attached to these elements, satisfying the following
conditions (cf. 4.1.1 in Part I):
\begin{center} 
\begin{tabular}{ll}
(i) & $h_l \in {\fm}_P^{p^{e_l}}$ and $\overline{h_l} = (h_l 
\bmod 
{\fm}_P^{p^{e_l}+1}) = v_l^{p^{e_l}}$ with
$v_l
\in {\fm}_P/{\fm}_P^2$ for $l = 1, \ldots, N$,
\\ 
(ii) & $\{v_l \mid l = 1, \ldots, N\} \subset {\fm}_P/{\fm}_P^2$ consists of
$N$-distinct and $k$-linearly independent \\
&elements in the
$k$-vector space ${\fm}_P/{\fm}_P^2$.
\\ 
\end{tabular} 
\end{center} 
We also take a regular system of parameters $(x_1, \ldots, x_d)$ such that
$$(\on{asc}) \hskip.1in v_l = \overline{x_l} = (x_l 
\bmod{{\fm}_P^2}) \text{\ for\ } l = 1, \ldots, N.$$
(We say that a regular system of parameters $(x_1, \ldots, x_d)$ is
associated to 
$H = (h_1, \ldots, h_N)$ if the above condition $(\on{asc})$ is
satisfied.  For the description of the condition of 
$(x_1, \ldots, x_d)$ being
weakly-associated to $H$, we refer the reader to
Chapter 2.)

Now we claim that any element $f \in R_P$ has a power series expansion of
the form 
$$(\star) \hskip.1in f = \sum_{B \in \left({\bZ}_{\geq 0}\right)^N}a_BH^B
\text{\ where\ }a_B = \sum_{K \in \left({\bZ}_{\geq
0}\right)^d}b_{B,K}X^K,$$
with $b_{B,K}$ being a power series in terms of the remainder
$(x_{N+1},
\ldots, x_d)$ of the regular system of parameters, and with $K = (k_1,
\ldots, k_d)$ varying in the range satisfying the condition
$$0 \leq k_l \leq p^{e_l} - 1 \text{\ for\ }l = 1, \ldots, N \hskip.05in
\on{and} \hskip.05in k_l = 0 \text{\ for\ }l = N+1, \ldots, d.$$
The existence of the power series expansion of the form $(\star)$ and its
uniqueness (with respect to a fixed subset ${\mathcal H}$ and its
chosen (weakly-)associated regular system of parameters $(x_1,
\ldots, x_d)$) follow immediately, and are the results stated independent of
the notion of an idealistic filtration.

\komoku{Formal coefficient lemma.}\label{0.2.2.3} In the general setting as
described in \ref{0.2.2.2}, the discussion on the power series expansion
does not
involve the notion of an idealistic filtration.  
The most interesting and important 
result regarding the power series expansion, however, is obtained when we
introduce and
require the following condition for $H$ to satisfy, involving a
${\fD}$-saturated idealistic filtration ${\bI}_P$ over $R_P$:
\begin{center} 
\begin{tabular}{ll}
(iii) & $(h_l,p^{e_l}) \in {\bI}_P$ for $l = 1, \ldots, N$.
\\
\end{tabular} 
\end{center} 

Now the formal coefficient lemma claims
$$(f,a) \in {\bI}_P, f = \sum_{B \in \left({\bZ}_{\geq 0}\right)^N}a_BH^B
\Longrightarrow (a_B,a - |[B]|) \in \widehat{{\bI}_P} \text{\ for\ any\ }B
\in
\left({\bZ}_{\geq 0}\right)^N.$$ (We recall that, for $B = (b_1, \ldots,
b_N) \in \left({\bZ}_{\geq 0}\right)^N$, we denote $(p^{e_1}b_1, \ldots,
p^{e_N}b_N)$ by
$[B]$ and
$\sum_{l = 1}^Np^{e_l}b_l$ by
$|[B]|$.  For the definition of the completion $\widehat{{\bI}_P}$ of the
idealistic filtration ${\bI}_P$, we refer the reader to \S 2.4
in Part I.)  The statement of the formal coefficient lemma turns out to be
quite useful and powerful.  In fact, Lemma 4.1.4.1 (Coefficient Lemma) in
Part I can be
obtained as a corollary to this formal version in Part II.  We will see some
applications of the formal coefficient lemma not only in Chapter 3 when we
study the
invariant
$\widetilde{\mu}$, but also in Part III when we analyze the modifications
and transformations of an idealistic filtration and in Part IV when we 
give the description of our algorithm.

\end{subsection}

\begin{subsection}{Invariant $\widetilde{\mu}$.}\label{0.2.3} Chapter 3 is
devoted to the discussion of the invariant $\widetilde{\mu}$, which is a
counterpart
in the new setting of the Idealistic Filtration Program to the notion of the
``weak order'' in the classical setting, whose definition involves the
exceptional
divisors.  Naturally, when we carry out our algorithm, the definition of the
invariant $\widetilde{\mu}$ in the middle of its process involves the
exceptional
divisors created by blowups.  It also involves the subtle adjustments we
have to make to the notion of a leading generator system for a
$\fD_E$-saturated idealistic filtration in the presence of the exceptional
divisor $E$ (cf. \ref{0.2.1}).  However, we restrict the discussion of the
invariant
$\widetilde{\mu}$ in Part II to the one with no exceptional divisors taken
into consideration, and hence to the discussion which could only be directly
applied to
the situation {\it at year 0} of the algorithm.  The discussion with the
exceptional divisors taken into consideration, i.e.,
the discussion which can then be applied to the situation {\it after year 0}
of the algorithm, will be postponed until it finds an appropriate place
in Part III
or Part IV, where we will show how we should adjust the arguments in Part II
in the presence of the exceptional divisors.

\komoku{Definition of $\widetilde{\mu}$.}\label{0.2.3.1} Let ${\bI}$ be a
${\fD}$-saturated idealistic filtration over $R$ as before.
Let $P \in \on{Spec}\hskip.02in R \subset W$ be a closed point.  Take a
leading generator system ${\bH}$ for ${\bI}_P$, and let ${\mathcal H}$ be
the set
consisting of its elements.  Recall that in 3.2.2 in Part I we set
$$\mu_{\mathcal H}({\bI}_P) = \inf\left\{\mu_{\mathcal H}(f,a) :=
\frac{\on{ord}_{\mathcal H}(f)}{a} \mid (f,a) \in {\bI}_P, a >
0\right\}$$  where
$$\on{ord}_{\mathcal H}(f) = \sup\{n \in {\bZ}_{\geq 0} \mid f \in {\fm}_P^n
+ ({\mathcal H})\},$$
and that we define the invariant $\widetilde{\mu}(P)$ by the formula
$$\widetilde{\mu}(P) = \mu_{\mathcal H}({\bI}_P).$$
There are two main issues concerning the invariant $\widetilde{\mu}(P)$.
\begin{center}
\begin{tabular}{ll}
Issue 1: &Is $\widetilde{\mu}(P)$ independent of the choice of ${\bH}$ and
hence of ${\mathcal H}$ ?\\
Issue 2: &Is $\widetilde{\mu}$ upper semi-continuous as a function of 
\\
&the
(closed) point $P \in \on{Spec}\hskip.02in R \subset W$ ? 
\end{tabular}
\end{center}

\komoku{$\widetilde{\mu}(P)$ is independent of the choice of ${\mathcal
H}$.}\label{0.2.3.2} We settled Issue 1 affirmatively via Coefficient Lemma
in Part I.  We would like to emphasize, on one hand, that we carried out 
the entire argument in Part I at the algebraic level of a local ring. 
This argument, showing that the invariant
$\widetilde{\mu}(P)$ is determined independent of the choice of a leading
generator system, seems to be in contrast to the argument  by W{\l}odarczyk,
where he
uses some (analytic) automorphism of the completion of the local ring,
showing that certain invariants are determined independent of the choice of
a hypersurface
of maximal contact via the notion of homogenization.  Note that the notion
of a leading generator system is a collective substitute for the notion of
a hypersurface of maximal contact. (cf. 0.2.3.2.1 in Part I).

We remark, on the other hand, that we can give an analytic interpretation of
the invariant $\widetilde{\mu}(P)$ using the power series expansion
discussed in
Chapter 2.  In fact, we see that $\on{ord}_{\mathcal H}(f)
=
\on{ord}(a_{\mathbb O})$ where $a_{\mathbb O}$ with ${\mathbb O} 
= (0, \ldots, 0) \in
\left({\bZ}_{\geq 0}\right)^N$ is the ``constant term'' of the power series
expansion of the form $(\star)$.  This explicit interpretation leads to an
alternative way to settle Issue 1, though quite similar in spirit to the
proof
at the algebraic level, via the formal coefficient lemma.  Note that
$\widetilde{\mu}(P)$ is rational, i.e.,
$\widetilde{\mu}(P) \in {\bQ}$, if we assume that
${\bI}$ is of r.f.g. type (and hence that so is
${\bI}_P$). 

\komoku{Upper semi-continuity of $(\sigma,\widetilde{\mu})$.}\label{0.2.3.3}
Regarding Issue 2, we have to emphasize first that the question asking the
upper
semi-continuity of the invariant $\widetilde{\mu}$ by itself is ill-posed,
and its answer is no when literally taken.  The precise and correct question
to ask is
the upper semi-continuity of the pair $(\sigma,\widetilde{\mu})$ with
respect to the lexicographical order.  Since the invariant $\sigma$ is upper
semi-continuous, this is equivalent to asking if the invariant
$\widetilde{\mu}$ is upper semi-continuous along the local maximum locus of
the invariant
$\sigma$.  We settle Issue 2 affirmatively in this precise form. 

The difficulty in studying the behavior of the invariant $\widetilde{\mu}(P)
=
\mu_{\mathcal H}({\bI}_P)$, as we let $P$ vary along the local maximum locus
of the invariant $\sigma$, lies in the fact that we also have to
change the leading generator system
${\bH}$ and hence ${\mathcal H}$ simultaneously.  This is caused by the fact
that our definition of a leading generator system is a priori ``pointwise''
in nature
and hence that we do not know, even if ${\bH}$ is a leading generator system
for ${\bI}_P$ at a point $P$, ${\bH}$ stays being a leading generator system
for
${\bI}_Q$ at a point $Q$ in a neighborhood of $P$.  In general, it does not.
There arises the need to modify a given leading generator system into one
which is
{\it uniformly pure} as discussed in \ref{0.2.1.5}.  With the modified and
uniformly pure leading generator system, the upper semi-continuity at issue
is reduced
to that of the multiplicity of a function in the usual setting.  The upper
semi-continuity can also be verified if we look at the power series
expansion with
respect to a uniformly pure leading generator system, and study the behavior
of its coefficients.
\end{subsection}

\begin{subsection}{Appendix.}\label{0.2.4} In the appendix, we report a new
development, which establishes the nonsingularity principle using only
the $\fD$-saturation 
and without using the $\fR$-saturation.  
Recall that in Part I we established the
nonsingularity principle using both the $\fD$-saturation and $\fR$-saturation
(cf.
0.2.3.2.4 and Chapter 4 in Part I).  This opens
up a possibility of constructing an algorithm, still in the frame work of
the Idealistic Filtration Program, using only
the $\fD$-saturation and without
using the $\fR$-saturation.  Note that the $\fR$-saturation invites 
the problem of
termination, which we specified in the introduction to Part
I as the only missing piece toward completing our algorithm
in positive characteristic.  Therefore, we believe that this new development
is a substantial step forward in our quest for establishing an algorithm for
resolution of singularities in positive characteristic.
\end{subsection}

\vskip.1in

This finishes the description of the outline of Part II.

\end{section} 

\begin{section}{Current status of the 
Idealistic Filtration Program.}\label{0.3} 
It has been more than a year since we posted the original version 
of Part II on the electronic archive in August, 2007.  We would 
like to report on the current status of the IFP, and make a note to Part I.

\begin{subsection}{Current status}\label{0.3.1} Since the advent of 
the new nonsingularity principle as described in 0.2.4, we have been 
pursuing the scheme of constructing an algorithm using only the 
$\fD$-saturation (or $\fD_E$-saturation in the presence of an 
exceptional divisor $E$).  In fact, in characteristic zero, 
the scheme works almost perfectly providing an algorithm 
for local uniformization, with the triplet $(\sigma,\widetilde{\mu},s)$ 
being the unit to constitute the strand of invariants.  (In order to 
obtain the global resolution of singularities, one has to work a little 
bit more to fill in the gap between the maximum locus of the strand and 
the support of the modification of an idealistic filtration.  The gap is 
an anomaly observed when $\widetilde{\mu} = 1$.)  In positive 
characteristic, as we do not use the $\fR$-saturation any more, 
the denominators of the invariant $\widetilde{\mu}$ are well-controlled, 
being no obstruction to showing the termination of the algorithm.  
Recently, however, some ``bad'' examples surfaced; if we try to 
naively follow the analogy to the case in characteristic zero, 
the blowup of a ``$(\sigma,\widetilde{\mu},s)$-permissible'' center would 
lead to the strict increase of the invariant $\widetilde{\mu}$, 
violating the principle that the strand of invariants we construct 
should never increase after blowup.  A few of these examples also 
indicate that the so-called monomial case needs a more careful 
treatment in positive characteristic than in characteristic zero.  In 
order to overcome these pathologies observed in the ``bad'' examples, we 
introduce and insert a new invariant $\widetilde{\nu}$, making the 
quadruple $(\sigma,\widetilde{\mu},\widetilde{\nu},s)$ the new unit to 
constitute the strand of invariants.  The invariant $\widetilde{\nu}$ 
is closely related to the invariant ``$\nu$'' used in 
\cite{MR2427629} and \cite{CP0703302}.  
We are now testing if our algorithm, taking the 
``$(\sigma,\widetilde{\mu},\widetilde{\nu},s)$-permissible'' center in a 
quite explicit way, will provide a solution to the problem of local 
uniformization (and global resolution) in positive characteristic.  We 
want to emphasize that we consider these new developments as the events 
in the process of ``evolution'' of the IFP, rather than mutation, since 
the basic strategy of the IFP, as envisioned in Part I, remains intact 
throughout our project.  We reported the current status of the 
evolution of the IFP at the workshop held at RIMS in December of 2008, 
and we refer the reader to \cite{RIMS2008} for the precise content of the 
report.  More details will be published in our subsequent papers in the 
near future.

\end{subsection}

\begin{subsection}{Roles of $\sigma$ and 
$\widetilde{\mu}$.}\label{0.3.2} Despite all the changes in the 
evolution process of the IFP discussed above, the fundamental roles of 
the invariants $\sigma$ and $\widetilde{\mu}$, as the first two factors 
of the unit constituting the strand of invariants, remain unchanged. 
Therefore, the main portion of Part II, discussing these fundamental 
roles, remain unchanged. 

\end{subsection}

\begin{subsection}{Note to Part I}\label{0.3.3} After Part I was 
published from Publications of RIMS, we learned that the result stated 
as Proposition 2.3.2.4 in Part I has already appeared in \cite{LT74}.  The 
arguments both in Part I and \cite{LT74} are closely related to the 
classical results of Nagata \cite{MR0089836}.  
Due to our negligence, this fact was never mentioned in Part I, 
even though \cite{LT74} was included in the 
references for Part I.

\end{subsection} 

\end{section}

\end{chapter} 
\begin{chapter}{Invariant $\sigma$}
The purpose of this chapter is to investigate the basic properties of the
invariant $\sigma$.

In this chapter, $R$
represents the coordinate ring of an affine open subset
$\on{Spec}\hskip.02in R$ of a nonsingular variety $W$ of $\dim W = d$ over
an algebraically
closed  field $k$ of positive characteristic ${\on{char}}(k) = p$ or of
characteristic zero ${\on{char}}(k) = 0$, where in the latter case we
formally set 
$p = \infty$ (cf. 0.2.3.2.1 and Definition 3.1.1.1 (2) in Part I).

Let ${\bI}$ be a ${\fD}$-saturated idealistic filtration over $R$, and
${\bI}_P$ its localization at a closed point $P \in \on{Spec}\hskip.02in R
\subset W$.

\begin{section}{Definition and its computation.}\label{1.1}

\begin{subsection}{Definition of $\sigma$.}\label{1.1.1} First we recall the
definition, given in \S 3.2 in Part I,
of the invariant $\sigma$ at a closed point $P \in {\on{Spec}}\ R \subset
W$. 

\begin{defn}\label{1.1.1.1}
The invariant $\sigma$ at $P$, which we denote by $\sigma(P)$, is defined to
be the following infinite sequence indexed
by $e \in \bZ_{\geq 0}$
$$\sigma(P) = \left(d - l_{p^0}^{\on{pure}}(P), d - l_{p^1}^{\on{pure}}(P),
\cdots , d 
- l_{p^e}^{\on{pure}}(P), \cdots \right) 
= \left(d - l_{p^e}^{\on{pure}}(P)\right)_{e \in
{\bZ}_{\geq 0}}$$ 
where 
$$ d = \dim W, \quad
l_{p^e}^{\on{pure}}(P) = \dim_kL({\bI}_P)_{p^e}^{\on{pure}}.
$$ 

(We refer the reader to Chapter 3 in Part
I or \ref{0.2.1.1} in the introduction to Part II for the definitions of the
leading algebra $L({\bI}_P)$ of the idealistic filtration ${\bI}_P$, its
degree $p^e$ component
$L({\bI}_P)_{p^e}$, and its pure part $L({\bI}_P)_{p^e}^{\on{pure}}$.)
\end{defn} 

\begin{rem}\label{1.1.1.2}
\item[(1)] 
The reason why we take the infinite sequence $\left(d -
l_{p^e}^{\on{pure}}(P)\right)_{e \in {\bZ}_{\geq 0}}$ instead of the infinite
sequence $\left(l_{p^e}^{\on{pure}}(P)\right)_{e \in {\bZ}_{\geq 0}}$ 
is two-fold:
\begin{enumerate} 
\item[(i)] 
If we consider the infinite sequence 
$\left(l_{p^e}^{\on{pure}}(P)\right)_{e \in
{\bZ}_{\geq 0}}$, it is lower
semi-continuous as a function of $P$.  Taking the negative of each factor
($+ d$) of the sequence, we have our
invariant upper semi-continuous, as we will see below.  (We consider that
the bigger $\left(l_{p^e}^{\on{pure}}(P)\right)_{e \in {\bZ}_{\geq 0}}$ 
is, the better
the singularity is.  Therefore, as the measure of how bad the
singularity is, it is also natural to define our
invariant as its negative.)
\item[(ii)] 
We reduce the problem of resolution of singularities of an abstract variety
$X$ 
to that of embedded resolution.  Therefore, it would be desirable
or even necessary to come up with an algorithm which would induce
the ``same'' process of resolution of singularities, no matter what
ambient variety $W$ we choose for an embedding
$X\hookrightarrow W$ (locally).

While the infinite sequence 
$\left(l_{p^e}^{\on{pure}}(P)\right)_{e \in {\bZ}_{\geq
0}}$ (or its negative $\left(- l_{p^e}^{\on{pure}}(P)\right)_{e 
\in {\bZ}_{\geq 0}}$)
is dependent of the choice of $W$, the infinite
sequence $\left(\dim W -  l_{p^e}^{\on{pure}}(P)\right)_{e \in 
{\bZ}_{\geq 0}}$ is not.
Therefore, the latter
is more appropriate as an invariant toward constructing such an algorithm.
\end{enumerate} 
\item[(2)] 
The dimension of the pure part is non-decreasing as a function of $e
\in \bZ_{\geq 0}$, and is  uniformly bounded from above by
$d =\dim W$, i.e., 
$$0 \leq l_{p^0}^{\on{pure}}(P) \leq l_{p^1}^{\on{pure}}(P) \leq\dotsb \leq
l_{p^{e - 1}}^{\on{pure}}(P) \leq
l_{p^e}^{\on{pure}}(P) \leq \dotsb\leq d = \dim W$$
and hence stabilizes after some point, i.e., there exists
$e_M \in \bZ_{\geq 0}$ such that
$$l_{p^e}^{\on{pure}}(P) = l_{p^{e_M}}^{\on{pure}}(P) \text{ for }e \geq e_M.$$
Therefore, although $\sigma(P)$ is an infinite sequence by definition,
essentially we are only looking at some finite
part of it. 
\item[(3)] 
In characteristic zero, the invariant $\sigma(P)$ consists of
only one term $d - l_{p^0}^{\on{pure}}$, 
while the remaining terms $d - l_{p^e}^{\on{pure}}$ are not defined for 
$e > 0$, as we set $p = \infty$ in characteristic zero.  (However, we 
may still say $\sigma(P)$ is an infinite sequence and write $\sigma(P) 
\in \prod_{e \in {\bZ}_{\geq 0}}{\bZ}_{\geq 0}$, for the sake of 
simplicity of presentation, intentionally ignoring the particular 
situation in characteristic zero.)  Note that $l_{p^0}^{\on{pure}} = 
l_{p^0} = \dim_kL({\bI}_P)_1$ can be regarded as the number indicating 
``how many linearly independent 
hypersurfaces of maximal contact we can take'' for ${\bI}_P$ 
(cf. Chapter 3 in Part I).
\end{rem} 
\end{subsection}

\begin{subsection}{Computation of $\sigma$.}\label{1.1.2}
The next lemma computes $l_{p^e}^{\on{pure}}(P)$ in terms of $l_{p^e}(P)$
and in terms of $l_{p^{\alpha}}^{\on{pure}}(P)$ for ${\alpha} = 0, \ldots,
e-1$, which we can
assume inductively have already been computed.  We also see that
$l_{p^e}(P)$ can be computed as the rank of a certain ``Jacobian-like''
matrix, and
hence that it is lower semi-continuous as a function of $P$.  This
immediately leads to
the lower semi-continuity of the sequence 
$\left(l_{p^e}^{\on{pure}}(P)\right)_{e \in
{\bZ}_{\geq 0}}$ and hence to the upper semi-continuity of
$\sigma(P) = \left(d - l_{p^e}^{\on{pure}}(P)\right)_{e \in 
{\bZ}_{\geq 0}}$ as a
function of $P$.  We will discuss the upper semi-continuity of $\sigma$ in
detail in the next section.

\begin{lem}\label{1.1.2.1} 

$\mathbf{Case: P \not\in \on{Supp}(\bI).}$

In this case, since we assume $\bI$ is $\fD$-saturated 
and since so is $\bI_P$, we observe that 
$$\bI_P = R_P \times {\bR}.$$  
Accordingly, the invariant $\sigma(P)$ takes the absolute 
minimum ${\mathbb O}$ in the value set of the 
\linebreak[4]
invariant $\sigma$, i.e.,
$$\sigma(P) = \left(\sigma(P)(e)\right)_{e \in {\bZ}_{\geq 0}} = 
{\mathbb O} \text{ with }\sigma(P)(e) = 0 \hskip.1in \forall e \in 
{\bZ}_{\geq 0}.$$

$\mathbf{Case: P \in \on{Supp}(\bI).}$

In this case, fixing $e \in {\bZ}_{\geq 0}$, we compute 
$l_{p^e}^{\on{pure}}$ in the following manner:

Suppose we have already computed $l_{p^{\alpha}}^{\on{pure}}(P)$ for
${\alpha} = 0, \ldots, e-1$.

Let $0 \leq e_1 < \cdots < e_K \leq e-1$ be the integers indicating the
places
where $l_{p^{\alpha}}^{\on{pure}}(P)$ jumps, i.e.,
\begin{eqnarray*}
0 &=& l_{p^0}^{\on{pure}}(P) = \cdots = l_{p^{e_1-1}}^{\on{pure}}(P) \\
&<& l_{p^{e_1}}^{\on{pure}}(P) = \cdots = l_{p^{e_2-1}}^{\on{pure}}(P) \\
&& \cdots \\
&<& l_{p^{e_K}}^{\on{pure}}(P) = \cdots = l_{p^{e-1}}^{\on{pure}}(P). \\
\end{eqnarray*}
Introduce variables $\{v_{ij}\}_{i = 1}^K$ where the
second subscript $j$ ranges from $1$ to
$l_{p^{e_i}}^{\on{pure}}(P) - l_{p^{e_{(i-1)}}}^{\on{pure}}(P)$, i.e., $j =
1, \ldots,
l_{p^{e_i}}^{\on{pure}}(P) - l_{p^{e_{(i-1)}}}^{\on{pure}}(P)$.

Then we compute $l_{p^e}^{\on{pure}}(P)$ as follows:
$$l_{p^e}^{\on{pure}}(P) = l_{p^e}(P) - l_{p^e}^{\on{mixed}}(P)$$
where the number $l_{p^e}^{\on{mixed}}(P)$ is by definition given by the
formula below
\begin{displaymath}
\begin{array}{ll}
&l_{p^e}^{\on{mixed}}(P) \\
&= \#\{\on{monomials\ of\ the\ form\ }\prod_{i = 1}^K
\left(v_{ij}^{p^{e_i}}\right)^{b_{ij}} \mid
\sum_{i,j}p^{e_i}b_{ij} = p^e, \hskip.05in p^{e_i}b_{ij} \neq p^e
\hskip.05in \forall ij\}.\\
\end{array}
\end{displaymath}

Moreover, take a set of generators $\{s_1, \ldots, s_r\}$ of the ideal
${\bI}_{p^e}$ of 
the idealistic filtration ${\bI}$ at level $p^e$, i.e., $(s_1, \ldots, s_r)
= {\bI}_{p^e} \subset R$.  Let $(x_1, \ldots, x_d)$ be 
a regular system of parameters at $P$.  Then
$$
l_{p^e}(P) = \on{rank} \left[\partial_{X^I}(s_{\beta})\right]_{|I| =
p^e}^{{\beta} = 1, \ldots, r}.$$
\end{lem}
\begin{proof} $\mathbf{Case: P \not\in \on{Supp}(\bI).}$  
In this case, by definition, there exists an element $(f,a) \in \bI_P$ 
with $a > 0$ such that $\on{ord}_P(f) < a$.  There also exists an 
appropriate differential operator $d$ of degree $\on{ord}_P(f)$ such 
that $d(f) = u$ is a unit of $R_P$.  Then we have
$$(d(f),a - \on{ord}_P(f)) = (u,a - \on{ord}_P(f)) \in \bI_P$$
and hence by condition (differential) in Definition 2.1.2.1 in Part I
$$(\bI_P)_{a - \on{ord}_P(f)} = R_P.$$
This implies by condition (ii) in Definition 2.1.1.1 in Part I that
$$(\bI_P)_{n(a - \on{ord}_P(f))} = R_P \hskip.1in \forall n \in \bZ_{> 0}.$$
We conclude then by condition (iii) in Definition 2.1.1.1 in Part I that
$$\bI_P = R_P \times {\bR}.$$
From this the assertion on $\sigma(P)$ easily follows, since we have 
$L(\bI_P) = G_P$.

\vskip.1in

$\mathbf{Case: P \in \on{Supp}(\bI).}$  Let
$$L({\bI}_P) = \bigoplus_{n \in {\bZ}_{\geq 0}}L({\bI}_P)_n \subset G_P =
\bigoplus_{n \in {\bZ}_{\geq 0}}{\fm}_P^n/{\fm}_P^{n+1}$$
be the leading algebra of ${\bI}_P$.

By Lemma 3.1.2.1 in Part I, we can choose $\{e_1 < \cdots < e_M\} \subset
{\bZ}_{\geq 0}$ and $V_1 \sqcup \cdots \sqcup V_M \subset G_1$
with $V_i =
\{v_{ij}\}_j$ satisfying the following conditions

\begin{center} 
\begin{tabular}{ll}
(i) & $F^{e_i}(V_i) \subset L({\bI}_P)_{p^{e_i}}^{\on{pure}}$ 
for $1 \leq i \leq M$,
\\ 
(ii) & $\bigsqcup_{e_i \leq e}F^e(V_i)$ is a $k$-basis of
$L({\bI}_P)_{p^e}^{\on{pure}}$ for any $e \in {\bZ}_{\geq 0}$.
\\ 
\end{tabular} 
\end{center} 

Since $L({\bI}_P)^{\on{pure}}$ generates $L({\bI}_P)$, 
we have $L({\bI}_P) = 
k[\bigsqcup_{i = 1}^MF^{e_i}(V_i)]$.

From this it follows that
\begin{displaymath}
\begin{array}{ll}
l_{p^e}(P) &= \#\{\on{monomials\ of\ the\ form\ }\prod_{e_i \leq e}
\left(v_{ij}^{p^{e_i}}\right)^{b_{ij}} \mid
\sum_{i,j}p^{e_i}b_{ij} = p^e\}\\
l_{p^e}^{\on{pure}} &= \#\{\on{monomials\ of\ the\ form\
}\left(v_{ij}^{p^{e_i}}\right)^{b_{ij}} \mid e_i \leq e, \hskip.1in
p^{e_i}b_{ij} = p^e\}\\
\end{array}
\end{displaymath}
and hence that
$$l_{p^e}^{\on{pure}} = l_{p^e}(P) - l_{p^e}^{\on{mixed}}$$
where
\begin{displaymath}
\begin{array}{ll}
&l_{p^e}^{\on{mixed}}(P) \\
&= \#\{\on{monomials\ of\ the\ form\ }\prod_{e_i \leq e}
\left(v_{ij}^{p^{e_i}}\right)^{b_{ij}} \mid
\sum_{i,j}p^{e_i}b_{ij} = p^e, \hskip.05in p^{e_i}b_{ij} \neq p^e
\hskip.05in \forall ij\}.\\
\end{array}
\end{displaymath}

\vskip.1in

The assertion in ``Moreover'' part follows from the fact 
that $L({\bI}_P)_{p^e}$ is generated
as a $k$-vector space by the degree $p^e$ terms of the power series
expansions of
$\{s_{\beta}\}_{\beta = 1, \ldots, r}$ with respect to a regular system of
parameters $(x_1, \ldots, x_d)$, i.e.,
$$L({\bI}_P)_{p^e} = \left\langle s_{\beta} \bmod
{\fm}_P^{p^e + 1} \mid \beta = 1, \ldots, r\right\rangle 
= \left\langle s_{\beta}
\bmod (x_1, \ldots,x_d)^{p^e + 1}
\mid \beta = 1, \ldots, r\right\rangle$$ 
and that their coefficients appear as the
entries of the matrix given in the statement, i.e.,
$$s_{\beta} = \sum_{|I| = p^e}\partial_{X^I}(s_{\beta})X^I 
\bmod (x_1, \ldots, x_d)^{p^e + 1}.$$

\vskip.1in

This completes the proof of Lemma \ref{1.1.2.1}.
\end{proof}

\begin{rem}\label{1.1.2.2} Let us consider 
$\tau(P) = \left(l_{p^e}(P)\right)_{e \in
{\bZ}_{\geq 0}}$.  Then noting $l_{p^0}(P) = l_{p^0}^{\on{pure}}(P)$, we
conclude by Lemma \ref{1.1.2.1} (1) that
$\sigma(P)$ determines $\tau(P)$ and vice versa.


In particular, for $P, Q \in {\fm}\text{-}\on{Spec}\hskip.02in R$, we have
\begin{center}
\begin{tabular}{ll}
&$\sigma(P) = \sigma(Q) \Longleftrightarrow \tau(P) = \tau(Q)$\\
&$\sigma(P) \geq \sigma(Q) \Longleftrightarrow \tau(P) \leq \tau(Q)$.\\
\end{tabular}
\end{center}
Therefore, the upper semi-continuity of the invariant $\sigma$, which we
will show in the next section, is equivalent to the lower semi-continuity of
the invariant $\tau$.
\end{rem}
\end{subsection}

\end{section} 
\begin{section}{Upper semi-continuity.}\label{1.2}

\begin{subsection}{Basic facts surrounding the definition of the upper
semi-continuity.}\label{1.2.1} In this subsection, we clarify some basic
facts surrounding the definition of the upper
semi-continuity.  We denote by $f:X \rightarrow T$ a function from a
topological space $X$ to a totally-ordered set $T$.

\begin{defn}\label{1.2.1.1} We say $f$ is upper
semi-continuous if the set 
$$X_{\geq t} := \left\{x \in X \mid f(x)\geq t\right\}$$
is closed for any $t \in T$.
\end{defn} 

\begin{lem}\label{1.2.1.2} Consider the conditions below:
\begin{center} 
\begin{tabular}{ll}
(i) & For any $x \in X$, there exists an open neighborhood $U_x$ such that
\\
&$f(x) \geq f(y)$ for any $y \in U_x$. \\
(ii) & The set $X_{> t} = \left\{x \in X \mid f(x) > t\right\}$ 
is closed for any $t \in T$. \\ 
(iii) & $f$ is upper semi-continuous.\\
\end{tabular} 
\end{center} 
Then we have the following implications:
$$(i) \Longleftrightarrow (ii) \Longrightarrow (iii).$$
Moreover, if $T$ is well-ordered (in the sense that every non-empty subset
has a least element), then conditions (ii) and (iii) are equivalent.
\end{lem}
\begin{proof}  The proof is elementary, and left to the reader as an
exercise.
\end{proof}

\begin{cor}\label{1.2.1.3} For the invariant
$\sigma:{\fm}\text{-}\on{Spec}\hskip.02in R \rightarrow \prod_{e \in
{\bZ}_{\geq 0}}
{\bZ}_{\geq 0}$, where the target space $\prod_{e \in {\bZ}_{\geq 0}}
{\bZ}_{\geq 0}$ is totally ordered with respect to the lexicographical
order, conditions (i), (ii), (iii) in Lemma \ref{1.2.1.2} are all
equivalent.
\end{cor}
\begin{proof}
As mentioned in Remark \ref{1.1.1.2} (2), the dimension of the pure part
$l_{p^e}^{\on{pure}}(P)$ is non-decreasing as a function of $e \in
{\bZ}_{\geq 0}$.  Accordingly, $\sigma(P)(e) = d - l_{p^e}^{\on{pure}}(P)$
is non-increasing as a function of $e \in
{\bZ}_{\geq 0}$.  Therefore, instead of $\prod_{e \in {\bZ}_{\geq
0}}{\bZ}_{\geq 0}$, we may take the subset $T = \left\{(t_e)_{e \in {\bZ}_{\geq
0}} \in \prod_{e \in {\bZ}_{\geq 0}}{\bZ}_{\geq 0} \mid
t_{e_1} \geq t_{e_2}
\on{if} e_1 > e_2\right\} \subset \prod_{e \in {\bZ}_{\geq 0}}
{\bZ}_{\geq 0}$ as the target space for $\sigma$.  Observe that 
$T$ is well-ordered (with respect to the total order induced by 
the one on $\prod_{e \in {\bZ}_{\geq 0}}{\bZ}_{\geq 0}$).  In fact, 
for a non-empty subset $S \subset T$, we can construct its least 
element $s_{\min} =
\left(s_{\min,e}\right)_{e \in {\bZ}_{\geq 0}}$ inductively by
the following formula:
$$s_{\min,e} = \min\left\{s_e \in {\bZ}_{\geq 0} \mid s = 
\left(s_i\right)_{i \in
{\bZ}_{\geq 0}} \in S \on{\ s.t.\ }s_i = s_{\min,i} \on{\ for\ }i < 
e\right\}.$$
Now the statement of the corollary follows from Lemma \ref{1.2.1.2}.
\end{proof}

The following basic description of the stratification into the level sets
can be seen easily, and its proof is left to the reader.

\begin{cor}\label{1.2.1.4} Let $f:X \rightarrow T$ be an upper
semi-continuous function.  Suppose that $X$ is noetherian, and that $T$ is
well-ordered.  Then
$f$ takes only finitely many values over $X$, i.e.,
$$\left\{f(x) \mid x \in X\right\} = \left\{t_1 < \cdots < 
t_n\right\} \subset T.$$
Accordingly, we have a strictly decreasing finite sequence of closed subsets
$$X = X_{\geq t_1} \supsetneqq \cdots \supsetneqq X_{\geq t_n} \supsetneqq
\emptyset,$$ 
which provides the stratification of $X$ into the level sets
$$\left\{x \in X \mid f(x) = t_i\right\} = X_{\geq t_i} 
\setminus X_{\geq t_{i+1}}
\hskip.1in \on{for} \hskip.1in i = 1, \ldots, n.$$
\end{cor}

\end{subsection}

\begin{subsection}{Upper semi-continuity of the invariant
$\sigma$.}\label{1.2.2}

\begin{prop}\label{1.2.2.1} The invariant\ 
$$\sigma:{\fm}\text{-}\on{Spec}\hskip.02in R \rightarrow \prod_{e \in
{\bZ}_{\geq 0}}{\bZ}_{\geq 0}$$ is upper semi-continuous.
\end{prop}
\begin{proof}
Set $X = {\fm}\text{-}\on{Spec}\hskip.023in R$ for notational simplicity.

Given $t = \left(t_e\right)_{e \in {\bZ}_{\geq 0}}
\in \prod_{e \in {\bZ}_{\geq 0}}{\bZ}_{\geq 0}$ and $n \in {\bZ}_{\geq 0}$,
we denote by $t_{\leq n}$ the truncation of $t$ up to the $n$-th term,
i.e., $t_{\leq n} = \left(t_e\right)_{e = 0}^n \in \prod_{e = 0}^n
{\bZ}_{\geq 0}$.

We define $\sigma_{\leq n}: X \rightarrow \prod_{e = 0}^n{\bZ}_{\geq 0}$ by
$\sigma_{\leq n}(x) =
\left(\sigma(x)\right)_{\leq n}$ for $x \in X$.

We also set $X(t_{\leq n}) = \left\{x \in X \mid \sigma_{\leq n}(x) 
\geq t_{\leq n}\right\}$.  Then
$$X_{\geq t} = \bigcap_{n = 0}^{\infty}X(t_{\leq n}).$$

In order to show the upper semi-continuity of $\sigma$, we have to show
$X_{\geq t}$ is closed for any $t = \left(t_e\right)_{e \in {\bZ}_{\geq 0}}
\in \prod_{e \in {\bZ}_{\geq 0}}{\bZ}_{\geq 0}$.  The above equality implies
that it suffices to show $X(t_{\leq n})$ is closed for any $n \in
{\bZ}_{\geq 0}$.  This follows if
the function $\sigma_{\leq n}$ is upper semi-continuous, which we will show
by induction on $n$.

The function $\sigma_{\leq 0} = d - l_{p^0}^{\on{pure}}$ is upper
semi-continuous, since $l_{p^0}^{\on{pure}} = l_{p^0}$ is lower
semi-continuous (cf. Lemma \ref{1.1.2.1}).

Assume we have shown $\sigma_{\leq n-1}$ is upper semi-continuous.  We show
then that the function $\sigma_{\leq n}$ satisfies condition (i) in Lemma
\ref{1.2.1.2} and hence that it is upper semi-continuous.

Suppose we are given $x \in X$.  Since $\sigma_{\leq n-1}$ is upper
semi-continuous (and since the target space $\prod_{e = 0}^{n-1}{\bZ}_{\geq
0}$ is well ordered), there
exists an open neighborhood
$U_x$ such that
$\sigma_{\leq n-1}(x) \geq \sigma_{\leq n-1}(y)$ for any $y \in U_x$  (cf.
Lemma \ref{1.2.1.2}).  Since the function $l_{p^n}$ is lower 
semi-continuous (cf. Lemma \ref{1.1.2.1}), by shrinking
$U_x$ if necessary, we may assume that $l_{p^n}(x) \leq l_{p^n}(y)$ for any
$y \in U_x$.

Take $y \in U_x$.  

If $\sigma_{\leq n-1}(x) > \sigma_{\leq n-1}(y)$, then we obviously have
$\sigma_{\leq n}(x) > \sigma_{\leq n}(y)$.

If $\sigma_{\leq n-1}(x) = \sigma_{\leq n-1}(y)$, then from the definition
it follows that $l_{p^{\alpha}}^{\on{pure}}(x) =
l_{p^{\alpha}}^{\on{pure}}(y)$ for ${\alpha} = 0,
\ldots, n-1$.  This implies by Lemma \ref{1.1.2.1} (1) that
$l_{p^n}^{\on{mixed}}(x) = l_{p^n}^{\on{mixed}}(y)$.  
Therefore, we conclude that
$$l_{p^n}^{\on{pure}}(x) = l_{p^n}(x) - l_{p^n}^{\on{mixed}}(x) \leq
l_{p^n}(y) - l_{p^n}^{\on{mixed}}(y) = l_{p^n}^{\on{pure}}(y),$$
and hence that
$$d - l_{p^n}^{\on{pure}}(x) \geq d - l_{p^n}^{\on{pure}}(y).$$
Thus we have $\sigma_{\leq n}(x) \geq \sigma_{\leq n}(y)$.

This shows that $\sigma_{\leq n}$ satisfies condition (i) in Lemma
\ref{1.2.1.2} and hence that it is upper semi-continuous, and completes the
induction.

Therefore, we conclude $\sigma:{\fm}\text{-}\on{Spec}\hskip.02in R
\rightarrow \prod_{e \in {\bZ}_{\geq 0}}
{\bZ}_{\geq 0}$ is upper semi-continuous.

This completes the proof of Proposition \ref{1.2.2.1}.
\end{proof}

\begin{cor}\label{1.2.2.2} We can extend the domain from
${\fm}\text{-}\on{Spec}\hskip.02in R$ to $\on{Spec}\hskip.02in R$ to have
the invariant $\sigma:\on{Spec}\hskip.02in R \rightarrow
\prod_{e \in {\bZ}_{\geq 0}} {\bZ}_{\geq 0}$, by defining
$$\sigma(Q) = \min\left\{\sigma(P) \mid P \in {\fm}\text{-}\on{Spec}
\hskip.02in R, \hskip.02in P \in \overline{Q}\right\} \hskip.05in 
\on{for} \hskip.05in Q \in \on{Spec}\hskip.02in R.$$
The formula is equivalent to saying that $\sigma(Q)$ is equal to $\sigma(P)$
with $P$ being a general closed point on $\overline{Q}$.  The invariant
$\sigma$ with the
extended domain is also upper semi-continuous.

Moreover, since $\on{Spec}\hskip.02in R$ is noetherian and since $\prod_{e
\in {\bZ}_{\geq 0}}{\bZ}_{\geq 0}$ can be replaced with the well-ordered set
$T$ as described in the proof of Corollary \ref{1.2.1.3}, conditions (i) and
(ii) in Lemma \ref{1.2.1.2}, as well as the assertions of Corollary
\ref{1.2.1.4}, hold for the upper semi-continuous
function
$\sigma:\on{Spec}\hskip.02in R \rightarrow T$.
\end{cor}

\begin{proof} Observe that, given $Q \in \on{Spec}\hskip.023in R$, the
formula for $\sigma(Q)$ is well-defined, since the existence of the minimum
(i.e., the least element) on the right hand
side is guaranteed by the fact that the value set of the invariant $\sigma$
is well-ordered (cf. the proof of Corollary \ref{1.2.1.3}).  Note that there
exists a non-empty dense open subset
$U$ of
$\overline{Q} \cap {\fm}\text{-}\on{Spec}\hskip.02in R$ such that $\sigma(Q)
= \sigma(P)$ for $P \in U$, a fact implied by condition (i) of the upper
semi-continuity of the invariant $\sigma$.
The upper semi-continuity of the invariant $\sigma$ with the extended domain
$\on{Spec}\hskip.023in R$ is immediate from the upper semi-continuity of the
invariant $\sigma$ with the original domain
${\fm}\text{-}\on{Spec}\hskip.023in R$.

The ``Moreover'' part follows immediately from the statements of Lemma
\ref{1.2.1.2}, Corollary \ref{1.2.1.3} and Corollary \ref{1.2.1.4}.

This completes the proof of Corollary \ref{1.2.2.2}.
\end{proof}

\end{subsection}
\end{section} 

\begin{section}{Local behavior of a leading generator system.}\label{1.3}

\begin{subsection}{Definition of a leading generator system and a remark
about the subscripts}\label{1.3.1} We say that a subset ${\bH} =
\left\{(h_l,p^{e_l})\right\}_{l = 1}^N \subset {\bI}_P$ with 
nonnegative integers $0 \leq e_1 \leq \cdots \leq e_N$ 
attached is a leading generator system (of the localization ${\bI}_P$ of the
${\fD}$-saturated idealistic filtration ${\bI}$ over $R$ at a
closed point
$P
\in {\fm}\text{-}\on{Spec}\hskip.02in R$), if the leading terms of its
elements provide a specific set of generators for the leading algebra
$L({\bI}_P)$.  More precisely, it satisfies the following conditions (cf.
Definition 3.1.3.1 in Part I):
\begin{center} 
\begin{tabular}{ll}
(i) & $h_l \in {\fm}_P^{p^{e_l}}$ and $\overline{h_l} = (h_l 
\bmod {\fm}_P^{p^{e_l}+1}) \in
L({\bI}_P)_{p^{e_l}}^{\on{pure}}$ for
$l = 1, \ldots, N$,
\\ 
(ii) & $\left\{\overline{h_l}^{p^{e - e_l}} \mid e_l \leq e\right\}$ 
consists of $\#\left\{l
\mid e_l \leq e\right\}$-distinct elements, and forms a $k$-basis of \\
&$L({\bI}_P)_{p^e}^{\on{pure}}$ for any $e \in {\bZ}_{\geq 0}$.
\\ 
\end{tabular} 
\end{center} 

Since the leading algebra $L({\bI}_P)$ is generated by its pure part
$L({\bI}_P)^{\on{pure}} = \bigsqcup_{e \in {\bZ}_{\geq
0}}L({\bI}_P)_{p^e}^{\on{pure}}$ (cf. \ref{0.2.1.1}), we conclude from
condition (ii) that the leading terms of ${\bH}$
$$\left\{\overline{h_l} = (h_l \bmod {\fm}_P^{p^{e_l} + 1})
\right\}_{l = 1}^N$$
provide a set of generators for $L({\bI}_P)$, 
i.e.,$L({\bI}_P) = k[\left\{\overline{h_l}\right\}_{l =
1}^N]$.

We remark that, for the subscripts of the leading generator system ${\bH}$,
we sometimes use the letter ``$l$'' as above, writing ${\bH} =
\left\{(h_l,p^{e_l})\right\}_{l = 1}^N$ with nonnegative integers 
$0 \leq e_1 \leq \cdots \leq e_N$ attached, and that some other times 
we use the letters $i$ and $j$, writing 
${\bH} = \left\{(h_{ij},p^{e_i})\right\}_{i = 1, \ldots,
M}^j$ with nonnegative integers
$0 \leq e_1 < \cdots < e_M$ attached.  In the latter use of the subscripts,
conditions (i) and (ii) are written as in 3.1.3 of Part I:

\begin{center} 
\begin{tabular}{ll}
(i) & $h_{ij} \in {\fm}_P^{p^{e_i}}$ and $\overline{h_{ij}} = (h_{ij}
\bmod {\fm}_P^{p^{e_i}+1}) \in
L({\bI}_P)_{p^{e_i}}^{\on{pure}}$ for
any $ij$,
\\ 
(ii) & $\left\{\overline{h_{ij}}^{p^{e - e_i}} \mid e_i \leq e\right\}$ 
consists of
$\#\left\{ij \mid e_i \leq e\right\}$-distinct elements, and forms \\
&a $k$-basis of $L({\bI}_P)_{p^e}^{\on{pure}}$ for any $e \in {\bZ}_{\geq
0}$. 
\\ 
\end{tabular} 
\end{center} 

In the future, we use the subscripts in both ways, while choosing one at a
time, depending upon the situation and its convenience.
\end{subsection}

\begin{subsection}{A basic question.}\label{1.3.2} Let
${\bH}$ be a leading generator system of ${\bI}_P$.  If we take a
neighborhood
$U_P$ of
$P$ small enough, then
${\bH}$ is a subset of ${\bI}_Q$ for any closed point $Q \in U_P
\cap {\fm}\text{-}\on{Spec}\hskip.02in R$.  We may then ask the following
question regarding the local behavior of the
leading generator system:

Is
${\bH}$ a leading generator system of
${\bI}_Q$ ?

A moment of thought reveals that the answer to this question in general is
no.  In fact, due to the upper semi-continuity of the invariant $\sigma$, by
shrinking $U_P$ if necessary, we may assume $\sigma(P) \geq \sigma(Q)$ for
any closed point $Q \in U_P \cap {\fm}\text{-}\on{Spec}\hskip.02in R$.  If
$\sigma(P) > \sigma(Q)$, then there is no way
that ${\bH}$ could be a leading generator system of ${\bI}_Q$.  (Note 
that the invariant $\sigma$ is completely determined by 
the leading generator system.)

We refine our question to avoid the obvious calamity as above::

Is ${\bH}$ a leading generator system of ${\bI}_Q$ for any closed point $Q
\in C \cap {\fm}\text{-}\on{Spec}\hskip.02in R \subset U_P \cap
{\fm}\text{-}\on{Spec}\hskip.02in R$ where $C = \{Q \in
U_P
\mid
\sigma(P) = \sigma(Q)\}$ ?

The answer to this question, for an arbitrary leading generator system
${\bH}$ of ${\bI}_P$, is still no.  One of the conditions for ${\bH}$ to be
a
leading generator system of ${\bI}_P$ requires any element $(h_{ij},
p^{e_i}) \in {\bH}$ to be pure at $P$, i.e., $(h_{ij} 
\bmod {\fm}_P^{p^{e_i}+1}) \in L({\bI}_P)_{p^{e_i}}^{\on{pure}}$.
However, even when a closed point $Q \in U_P \cap
{\fm}\text{-}\on{Spec}\hskip.02in R$ satisfies the condition
$Q \in C \cap {\fm}\text{-}\on{Spec}\hskip.02in R$, some element
$(h_{ij},p^{e_i})$ may fail to be pure at $Q$, i.e., $(h_{ij} 
\bmod {\fm}_Q^{p^{e_i}+1}) \not\in L({\bI}_Q)_{p^{e_i}}^{\on{pure}}$,
and hence ${\bH}$ fails to be a leading generator system at $Q$.

Now we refine our question further:

Can we modify a given generator system ${\bH}$ of ${\bI}_P$ into ${\bH}'$ so
that ${\bH}'$ stays being a leading generator system of ${\bI}_Q$ for
any closed point $Q \in C \cap {\fm}\text{-}\on{Spec}\hskip.02in R \subset
U_P \cap {\fm}\text{-}\on{Spec}\hskip.02in R$ where $C = \{Q \in U_P \mid
\sigma(P) = \sigma(Q)\}$ ?

The main goal of the next subsection is to give an affirmative answer to
this last question (adding one extra condition of the point $Q$ being 
on the support $\on{Supp}(\bI)$ of the idealistic filtration), and also 
to give an explicit description of how we make the
modification.  We say we modify the given leading generator system into one
which is ``{\it uniformly pure}'' (along $C$ intersected 
with $\on{Supp}(\bI)$).

\end{subsection}

\begin{subsection}{Modification of a given leading generator system into one
which is uniformly pure.}\label{1.3.3}

\begin{defn}\label{1.3.3.1} Let ${\bH}$ be a leading generator system of the
localization ${\bI}_P$ of the ${\fD}$-saturated idealistic filtration $\bI$
over $R$ at a closed
point $P
\in {\fm}\text{-}\on{Spec}\hskip.02in R$.  We say ${\bH}$ is uniformly pure
(in a neighborhood $U_P$ of $P$ along the local maximum locus $C$ of the
invariant $\sigma$ intersected with the support $\on{Supp}(\bI)$ 
of the idealistic filtration) if there exists an open neighborhood 
$U_P$ of $P$ such that the following conditions are satisfied:

\begin{center} 
\begin{tabular}{ll}
(1) &${\bH} \subset {\bI}_Q \hskip.1in \forall Q \in U_P$,\\
(2) &$\sigma(P)$ is the maximum of the invariant $\sigma$ over $U_P$, i.e.,
$\sigma(P) \geq \sigma(Q) \hskip.1in \forall Q \in U_P$\\
(3) &$C = \{Q \in U_P \mid \sigma(P) = \sigma(Q)\}$ is a closed subset of
$U_P$, and \\
(4) &${\bH}$ is a leading generator system of ${\bI}_Q$ for any $Q \in C
\cap \on{Supp}({\bI}) \cap {\fm}\text{-}\on{Spec}\hskip.02in R$.\\
\end{tabular} 
\end{center} 

(For the definition of the support $\on{Supp}({\bI})$ of the idealistic
filtration ${\bI}$, we refer the reader to Definition 2.1.1.1 in Part I.)
\end{defn}

\begin{rem}\label{1.3.3.2} We remark that in condition (4) of Definition
\ref{1.3.3.1}, in order for ${\bH}$ to be uniformly pure, we require ${\bH}$
is a leading generator system of ${\bI}_Q$ for
any closed point ``$Q
\in C
\cap
\on{Supp}({\bI})
\cap {\fm}\text{-}\on{Spec}\hskip.02in R$'' (i.e., we only consider those
closed points on the support $\on{Supp}({\bI})$ of the idealistic filtration
${\bI}$),
where in the last form of the basic question in \ref{1.3.2} we merely wrote
``$Q \in C \cap {\fm}\text{-}\on{Spec}\hskip.02in R$''.  The reason to add
this extra condition on $Q$ (as mentioned in the last 
paragraph of \ref{1.3.2}) is as follows:

Consider the case when $\sigma(P) = {\mathbb O}$.  
(Recall that the symbol ${\mathbb O} = (0, \cdots, 0, \cdots)$ 
represents the absolute minimum in the value set of the invariant $\sigma$.)

By the upper semi-continuity of the invariant $\sigma$, for a sufficiently
small open neighborhood $U_P$ of $P$, we have $\sigma(Q) =
\sigma(P) = {\mathbb O}$ for any closed point 
$Q \in U_P \cap {\fm}\text{-}\on{Spec}\hskip.02in R$ 
and hence we have $C \cap
{\fm}\text{-}\on{Spec}\hskip.02in R = U_P \cap
{\fm}\text{-}\on{Spec}\hskip.02in R$.  On the other hand, 
the condition $\sigma(P) = {\mathbb O}$ implies that, given any leading
generator system ${\bH}$ of ${\bI}_P$, the elements $\{h_{ij}\}$ are
generators of the maximal ideal ${\fm}_P$ with
$\#\{ij\} = d$.  (Note that, in this case, all the elements of a leading
generator system are concentrated at level $1$, i.e., $1 = i = M$ and $0 =
e_1 = e_i = e_M$.)  Therefore,
${\bH}$ can not be a leading generator system of
${\bI}_Q$ for a closed point
$Q
\in U_P
\cap {\fm}\text{-}\on{Spec}\hskip.02in R$ if $Q \neq P$.  That is to say, it
would not satisfy the condition described in the last form of the basic
question.  However, in this case, we have either
$U_P
\cap
\on{Supp}({\bI}) =
\emptyset$ or
$U_P
\cap
\on{Supp}({\bI}) = \{P\}$ (if we take $U_P$ sufficiently small).  Therefore,
condition (4) in Definition \ref{1.3.3.1} is automatically satisfied.

Consider the case when $\sigma(P) \neq {\mathbb O}$.

In this case, we have $C \cap {\fm}\text{-}\on{Spec}\hskip.02in R = C \cap
\on{Supp}({\bI}) \cap {\fm}\text{-}\on{Spec}\hskip.02in R$, since any closed
point $Q \in C \cap
{\fm}\text{-}\on{Spec}\hskip.02in R$ (i.e., we have $\sigma(Q) = \sigma(P)
\neq {\mathbb O}$ is necessarily in the support
$\on{Supp}({\bI})$ of the idealistic filtration (cf. 
Lemma \ref{1.1.2.1}).  Therefore, there is no difference
between the condition in the last form of the basic question and condition
(4) in Definition \ref{1.3.3.1}.

In other words, the extra condition for $Q$ to be in the support
$\on{Supp}({\bI})$ is introduced so that we can avoid 
the ``obvious'' counter example to an affirmative answer to the last form 
of the basic question in the special case $\sigma(P) = {\mathbb O}$.

\end{rem}

\begin{prop}\label{1.3.3.3} Let ${\bH} = 
\left\{(h_{ij},p^{e_i})\right\}_{i = 1, \cdots, M}^j$ be a leading
generator system of the localization ${\bI}_P$ of the ${\fD}$-saturated
idealistic filtration $\bI$ over $R$ at a closed point 
$P\in {\fm}\text{-}\on{Spec}\hskip.02in R$, with nonnegative
integers $0 \leq e_1 < \cdots < e_M$ attached.  
Then ${\bH}$ can be modified into another leading generator system
$\bH'$ which is uniformly pure.

\vskip.1in

More precisely, there exists $\left\{g_{ijB}\right\} \subset {\fm}_P$, 
where the subscript $B$ ranges over the set
$$\on{Mix}_{{\bH},i} = \left\{B = (b_{\alpha\beta}) \in ({\bZ}_{\geq
0})^{\#{\bH}} \mid |[B]| = p^{e_i},
\hskip.05in b_{\alpha\beta} = 0 \hskip.05in \on{if} \hskip.05in \alpha \geq
i, \hskip.05in
\on{and}
\hskip.05in p^{e_{\alpha}}b_{\alpha\beta} \neq p^{e_i} \hskip.05in \forall
\alpha\beta\right\},$$
such that, setting $h_{ij}' = h_{ij} - \sum g_{ijB}H^B$, the modified set
${\bH}' =
\left\{(h_{ij}',p^{e_i})\right\}_{i = 1, \cdots, M}^j$ is 
a leading generator system of ${\bI}_P$ which is uniformly pure.

\end{prop}

\begin{proof} It suffices to prove that there exists an affine open
neighborhood $U_P = \on{Spec}\hskip.02in R_f$ of $P$, where $R_f$ is the
localization of $R$ by an element $f \in R$, such that
the following conditions are satisfied:

\begin{center} 
\begin{enumerate}
\item 
${\bH} \subset {\bI}_f$ (and hence ${\bH}' \subset {\bI}_f$ where
${\bH}'$ is described in condition (4)),\\
\item
$\sigma(P)$ is the maximum of the invariant $\sigma$ over $U_P$, i.e.,
$\sigma(P) \geq \sigma(Q) \hskip.1in \forall Q \in U_P$,\\
\item
$C = \left\{Q \in U_P \mid \sigma(P) = \sigma(Q)\right\}$ is 
a closed subset of $U_P$, and \\
\item
there exists $\left\{g_{ijB}\right\} \subset R_f$, where 
the subscript $B$ ranges over the set $\on{Mix}_{{\bH},i}$, 
such that $\left\{g_{ijB}\right\} \subset {\fm}_P$ and that, 
setting $h_{ij}' = h_{ij}- \sum g_{ijB}H^B$, 
the modified set 
$${\bH}' = \left\{(h_{ij}',p^{e_i})\right\}_{i = 1, \cdots, M}^j$$
is a leading generator system of ${\bI}_Q$ 
for any $Q \in C \cap \on{Supp}({\bI})
\cap {\fm}\text{-}\on{Spec}\hskip.02in R$.\\
\end{enumerate}
\end{center} 

\begin{step}{Check conditions (1), (2) and (3).} 

It is easy to choose an affine open neighborhood $U_P = \on{Spec}\hskip.02in
R_f$ of $P$ satisfying condition (1).  By the upper semi-continuity of the
invariant $\sigma$, we
may also assume condition (2) is satisfied (cf. condition (i) in Lemma
\ref{1.2.1.2}).  Then condition (3) automatically follows, since $C = U_P
\cap (\on{Spec}\hskip.02in R)_{\geq \sigma(P)}$ is
closed (cf. Definition \ref{1.2.1.1}).

We remark that in terms of the invariant $\tau$ (cf. Remark \ref{1.1.2.2})
conditions (2) and (3) are equivalent to the following
\begin{center} 
\begin{tabular}{ll}
$(2)_{\tau}$ &$\tau(P)$ is the minimum of the invariant $\tau$ over $U_P$,
i.e., $\tau(P) \leq \tau(Q) \hskip.1in \forall Q \in U_P$, and\\
$(3)_{\tau}$ &$C = \left\{Q \in U_P \mid \tau(P) = \tau(Q)\right\}$.\\
\end{tabular} 
\end{center} 
Now we have only to check, by shrinking $U_P$ if necessary, that condition
(4) is also satisfied.
\end{step}

\vskip.1in

\begin{step}{Preliminary analysis to check condition (4).}

First consider the idealistic filtration ${\bJ} = G_{R_f}({\bH})$ generated
by ${\bH}$ over $R_f$.  Note that ${\bJ} \subset {\bI}_f$ but that ${\bJ}$
may not be ${\fD}$-saturated.  In order to
distinguish the invariant $\tau$ for ${\bI}$ (or equivalently for ${\bI}_f$
over $U_P$) from the invariant $\tau$ for ${\bJ}$, we denote them by
$\tau_{\bI}$ and $\tau_{\bJ}$, respectively.

Since
$\tau_{\bJ}$ is lower semi-continuous, by shrinking $U_P$ if necessary, we
may assume 

\begin{center}
\begin{tabular}{ll}
$(2)_{\tau_{\bJ}}$ &$\tau_{\bJ}(P)$ is the minimum of the invariant
$\tau_{\bJ}$ over $U_P$, i.e., $\tau_{\bJ}(P) \leq \tau_{\bJ}(Q) \hskip.1in
\forall Q \in U_P$.
\end{tabular} 
\end{center} 

For any closed point $Q \in C \cap \on{Supp}({\bI}) \cap
{\fm}\text{-}\on{Spec}\hskip.02in R$, we compute

\begin{displaymath}
\tau_{\bI}(P) = \tau_{\bJ}(P) \leq \tau_{\bJ}(Q) \leq \tau_{\bI}(Q) =
\tau_{\bI}(P).
\end{displaymath}

We remark that the first equality is a consequence of the fact 
that the set 
$\left\{\overline{h_{ij,P}} = \right.$
$\left.(h_{ij}\bmod {\fm}_P^{p^{e_i}+1})\right\}_{i = 1, \cdots, M}^j$ 
generates both
$L({\bI}_P)$ and
$L({\bJ}_P)$ as $k$-algebras, the second inequality is a consequence of
$(2)_{\tau_{\bJ}}$, the third inequality is a consequence of the inclusion
${\bJ} \subset {\bI}_f$, and that the last equality follows from the
definition of the closed subset $C$.

Therefore, we see that
$$\tau_{\bI}(P) = \tau_{\bJ}(P) = \tau_{\bJ}(Q) = \tau_{\bI}(Q) \hskip.1in
\forall Q \in C \cap \on{Supp}({\bI}) \cap {\fm}\text{-}\on{Spec}\hskip.02in
R.$$

\end{step}

\vskip.1in

\begin{step}{Some conclusions of the equality 
$\tau_{\bI}(P) = \tau_{\bJ}(P) = \tau_{\bJ}(Q) = \tau_{\bI}(Q)$ 
for any 
$Q \in C \cap \on{Supp}({\bI}) \cap {\fm}\text{-}\on{Spec}\hskip.02in R.$}

The equality obtained at the end of Step 2 leads to a few conclusions that 
we list below:

\begin{center}
\begin{tabular}{ll}
$(a)$ &The set $\left\{\overline{h_{ij,Q}} = (h_{ij} \bmod
{\fm}_Q^{p^{e_i}+1})\right\}_{i = 1, \cdots, M}^j$ generates 
$L({\bI}_Q)$ as a $k$-algebra 
\\
&for any $Q \in C \cap \on{Supp}({\bI}) \cap 
{\fm}\text{-}\on{Spec}\hskip.02in R$. \\
&Moreover $\left\{\overline{H_Q}^B \mid B = (b_{ij}), |[B]| = p^e, \hskip.02in 
\on{and} \hskip.02in b_{ij} = 0 \hskip.02in \on{if} \hskip.02in e_i >
e\right\}$ \\
&forms a basis of $L({\bI}_Q)_{p^e}$ as a
$k$-vector space, \\
&since it obviously generates $L({\bI}_Q)_{p^e}$ and since\\
&$\#\left\{\overline{H_Q}^B \mid B = (b_{ij}), 
|[B]| = p^e, \hskip.02in \on{and} 
\hskip.02in b_{ij} = 0 \hskip.02in \on{if} \hskip.02in e_i >
e\right\}$ \\
&$= \#\left\{\overline{H_P}^B \mid B = (b_{ij}), |[B]| = p^e, 
\hskip.02in \on{and} 
\hskip.02in b_{ij} = 0 \hskip.02in \on{if} \hskip.02in e_i >
e\right\}$\\
&$= l_{p^e}(P) = l_{p^e}(Q) = \dim_kL({\bI}_Q)_{p^e}$.\\
$(b)$ & There exist nonnegative integers $0 \leq e_1 < \cdots < e_M$, \\
&independent of $Q \in C \cap \on{Supp}({\bI}) \cap
{\fm}\text{-}\on{Spec}\hskip.02in R$, such that the jumping of\\
&the dimension of the pure part only occurs at these numbers, i.e.,\\ 
\end{tabular}
\end{center}

\begin{displaymath}
\begin{array}{ll}
0 &= l_{p^0}^{\on{pure}}(Q) = \cdots = l_{p^{e_1-1}}^{\on{pure}}(Q) \\
&< l_{p^{e_1}}^{\on{pure}}(Q) = \cdots = l_{p^{e_2-1}}^{\on{pure}}(Q) \\
& \cdots \\
&< l_{p^{e_M}}^{\on{pure}}(Q) = \cdots , \\
\end{array}
\end{displaymath}

\noindent
as $l_{p^e}^{\on{pure}}(Q) = l_{p^e}^{\on{pure}}(P)$ for any $Q \in C \cap 
\on{Supp}({\bI}) \cap
{\fm}\text{-}\on{Spec}\hskip.02in R$ and $e \in {\bZ}_{\geq 0}$ (cf. Remark 
\ref{1.1.2.2}).

Applying Lemma 3.1.2.1 in Part I to $L({\bI}_Q)$ for $Q \in C \cap 
\on{Supp}({\bI}) \cap {\fm}\text{-}\on{Spec}\hskip.02in R$, we see that we 
can take $V_{1,Q} \sqcup \cdots \sqcup V_{M,Q}
\subset G_{1,Q} = {\fm}_Q/{\fm}_Q^2$ with $V_{i,Q} =
\left\{v_{ij,Q}\right\}_j$, where $1 \leq j \leq l_{p^{e_i}}^{\on{pure}}(Q) - 
l_{p^{e_{i-1}}}^{\on{pure}}(Q)$, satisfying the following conditions

\begin{center} 
\begin{tabular}{ll} 
(i) & $F^{e_i}(V_{i,Q}) \subset L({\bI}_Q)_{p^{e_i}}^{\on{pure}}$ for $1 
\leq i \leq M$,
\\ 
(ii) & $\bigsqcup_{e_i \leq e}F^e(V_{i,Q})$ is a $k$-basis of 
$L({\bI}_Q)_{p^e}^{\on{pure}}$ for any $e \in {\bZ}_{\geq 0}$.
\\ 
\end{tabular} 
\end{center} 

Since $L({\bI}_Q)^{\on{pure}}$ generates $L({\bI}_Q)$, we have $L({\bI}_Q) = 
k[\bigsqcup_{i = 1}^MF^{e_i}(V_{i,Q})]$.

\medskip

Using this information, we also conclude the following.

\begin{center}
\begin{tabular}{ll}
$(c)$ &The $\overline{h_{ij,Q}}$ are all pure when $i = 1$, i.e., 
$\overline{h_{1j,Q}} \in L({\bI}_Q)_{p^{e_1}}^{\on{pure}}$, and we take 
$v_{1j,Q}' \in G_{1,Q}$ \\
& so that $F^{e_1}(v_{1j,Q}') = \overline{h_{1j,Q}}$ for $j = 1, \ldots,
l_{p^{e_1}}^{\on{pure}}(Q)$.\\
&As can be seen by induction on $i = 1, \ldots, M$, for each $ij$, \\ 
&there exists {\it uniquely} a set $\left\{c_{ijB,Q}\right\}_{B \in 
\on{Mix}_{{\bH},i}} \subset k$ such that \\
&$\overline{h_{ij,Q}} - \sum_{B \in 
\on{Mix}_{{\bH},i}}c_{ijB,Q}\overline{H_Q}^B$ is pure, i.e., 
$\overline{h_{ij,Q}} - \sum_{B \in 
\on{Mix}_{{\bH},i}}c_{ijB,Q}\overline{H_Q}^B \in
L({\bI}_Q)_{p^{e_i}}^{\on{pure}}$.
\\ &We take $v_{ij,Q}' \in G_{1,Q}$ such that $F^{e_i}(v_{ij,Q}') = 
\overline{h_{ij,Q}} - \sum_{B \in \on{Mix}_{{\bH},i}}
c_{ijB,Q}\overline{H_Q}^B$.\\ &Setting $V_{i,Q}' = 
\left\{v_{ij,Q}'\right\}_j$, we see 
that we can replace $V_{1,Q} \sqcup \cdots  \sqcup V_{M,Q}$ with \\
&$V_{1,Q}' \sqcup \cdots  \sqcup V_{M,Q}'$, i.e., \\
&\hskip.2in (i) $F^{e_i}(V_{i,Q}') \subset L({\bI}_Q)_{p^{e_i}}^{\on{pure}}$ 
for $1 \leq i \leq M$, \\
&\hskip.2in (ii) $\bigsqcup_{e_i \leq e}F^e(V_{i,Q}')$ is a $k$-basis of 
$L({\bI}_Q)_{p^e}^{\on{pure}}$ for any $e \in {\bZ}_{\geq 0}$ \\
&We also have $L({\bI}_Q) = k[\bigsqcup_{i = 1}^MF^{e_i}(V_{i,Q}')]$.\\
\end{tabular}
\end{center}

In fact, we prove below conclusion $(c)$, claiming the existence and 
uniqueness of such $\left\{c_{ijB,Q}\right\} \subset k$ as described 
above, showing 
simultaneously by induction on
$i$ that we can replace
$V_{1,Q} \sqcup \cdots  \sqcup V_{M,Q}$ with $V_{1,Q}' \sqcup \cdots \sqcup 
V_{i,Q}' \sqcup V_{i+1,Q} \sqcup \cdots \sqcup V_{M,Q}$ in the assertions of 
Lemma 3.1.2.1 in Part I, and hence
ultimately with $V_{1,Q}' \sqcup \cdots  \sqcup V_{M,Q}'$.

(Existence) By inductional hypothesis, we may replace
$V_{1,Q} \sqcup \cdots  \sqcup V_{M,Q}$ with 
$V_{1,Q}' \sqcup \cdots \sqcup 
V_{i-1,Q}' \sqcup V_{i,Q} \sqcup \cdots \sqcup V_{M,Q}$ in the assertions of 
Lemma 3.1.2.1 in Part I.  Expressing
$\overline{h_{ij,Q}}$ as a degree $p^{e_i}$ homogeneous polynomial 
in terms of 
$F^{e_1}(V_{1,Q}') \sqcup \cdots \sqcup 
F^{e_{i-1}}(V_{i-1,Q}') \sqcup F^{e_i}(V_{i,Q})$, we see that there exists 
$\left\{a_{ijB,Q}\right\}_{B \in \on{Mix}_{{\bH},i}} \subset k$ 
such that $\overline{h_{ij,Q}} - \sum_{B \in 
\on{Mix}_{{\bH},i}}a_{ijB,Q}F^*(V_Q')^B \in
L({\bI}_Q)_{p^{e_i}}^{\on{pure}}$, where
$F^*(V_Q') = (F^{e_{\alpha}}(v_{\alpha\beta,Q}'))$.  Note that, although 
$v_{\alpha\beta,Q}'$ has yet to be defined if $\alpha \geq i$, since 
$b_{\alpha\beta} = 0$ if $\alpha \geq i$ for $B =
(b_{\alpha\beta}) \in \on{Mix}_{{\bH},i}$, the expression 
$\overline{h_{ij,Q}} - \sum_{B \in \on{Mix}_{{\bH},i}}a_{ijB,Q}F^*(V_Q')^B$ 
is well-defined.  By substituting 
$$F^{e_{\alpha}}(v_{\alpha\beta,Q}') =
\overline{h_{\alpha\beta,Q}} - \sum_{B \in \on{Mix}_{{\bH},\alpha}} 
c_{\alpha\beta B,Q}\overline{H_Q}^B \hskip.1in \on{for} \hskip.1in 
\alpha <  i,$$
we see that there exists $\left\{c_{ijB,Q}\right\}_{B \in
\on{Mix}_{{\bH},i}}
\subset k$ such that 
$$\overline{h_{ij,Q}} - \sum_{B \in
\on{Mix}_{{\bH},i}}c_{ijB,Q}\overline{H_Q}^B \in 
L({\bI}_Q)_{p^{e_i}}^{\on{pure}}.$$

We remark that the set
$$\left\{F^{e_i - e_{\alpha}}\left(\overline{h_{\alpha\beta,Q}} - \sum_{B 
\in
\on{Mix}_{{\bH},\alpha}}c_{\alpha\beta 
B,Q}\overline{H_Q}^B\right)\right\}_{\alpha = 1, \ldots, i, \hskip.05in 
\beta = 1, \ldots, l_{p^{e_{\alpha}}}^{\on{pure}}(Q) -
l_{p^{e_{\alpha-1}}}^{\on{pure}}(Q)} \subset 
L({\bI}_Q)_{p^{e_i}}^{\on{pure}}$$
is linearly independent, since 
$$\left\{\overline{H_Q}^B \mid B = (b_{\alpha\beta}), |[B]| = p^{e_i}, 
\hskip.02in \on{and} \hskip.02in b_{\alpha\beta} = 0 \text{\ if\ }
e_{\alpha} > e_i\right\}$$ 
is linearly independent (cf. conclusion (a) above), and that 
its cardinality 
$\sum_{\alpha = 1}^i
\left(l_{p^{e_{\alpha}}}^{\on{pure}}\right.(Q)$
$\left.- l_{p^{e_{\alpha-1}}}^{\on{pure}}(Q)\right)$
is equal to
$l_{p^{e_i}}^{\on{pure}}(Q)$.  Therefore, we conclude that the above set 
forms a basis of $L({\bI}_Q)_{p^{e_i}}^{\on{pure}}$.

(Uniqueness) Suppose there exists another set 
$\left\{c_{ijB,Q}'\right\}_{B \in \on{Mix}_{{\bH},i}} \subset k$ such
that 
$\overline{h_{ij,Q}} - \sum_{B \in \on{Mix}_{{\bH},i}}c_{ijB,Q}'
\overline{H_Q}^B \in 
L({\bI}_Q)_{p^{e_i}}^{\on{pure}}$.  Then 
$$\sum_{B \in \on{Mix}_{{\bH},i}}c_{ijB,Q}\overline{H_Q}^B - \sum_{B \in 
\on{Mix}_{{\bH},i}}c_{ijB,Q}'\overline{H_Q}^B = \sum_{B \in 
\on{Mix}_{{\bH},i}}\left(c_{ijB,Q} -
c_{ijB,Q}'\right)\overline{H_Q}^B
\in L({\bI}_Q)_{p^{e_i}}^{\on{pure}}.$$
From the conclusion at the end of the argument for (Existence) it follows 
that there exists 
$$\left\{\gamma_{\alpha\beta}\right\}_{\alpha = 1, \ldots, i,
 \hskip.05in \beta = 1, \ldots, l_{p^{e_{\alpha}}}^{\on{pure}}(Q) -
l_{p^{e_{\alpha-1}}}^{\on{pure}}(Q)} \subset k$$
such that
\begin{eqnarray*}
& &\sum_{B \in \on{Mix}_{{\bH},i}}\left(c_{ijB,Q} - 
c_{ijB,Q}'\right)\overline{H_Q}^B \\
& & \hskip.1in = \sum_{\alpha = 1, \ldots, i, \hskip.05in \beta = 1, \ldots, 
l_{p^{e_{\alpha}}}^{\on{pure}}(Q) -
l_{p^{e_{\alpha-1}}}^{\on{pure}}(Q)} \gamma_{\alpha\beta}F^{e_i - 
e_{\alpha}}\left(\overline{h_{\alpha\beta,Q}} - \sum_{B \in
\on{Mix}_{{\bH},\alpha}}c_{\alpha\beta B,Q}\overline{H_Q}^B\right).\\
\end{eqnarray*}

Again since $\left\{\overline{H_Q}^B \mid B = (b_{\alpha\beta}), 
|[B]| = p^{e_i}, 
\hskip.02in \on{and} \hskip.02in b_{\alpha\beta} = 0 \hskip.02in \on{if} 
\hskip.02in e_{\alpha} >
e_i\right\}$ is linearly independent, we conclude that 
$\gamma_{\alpha\beta} = 0 
\hskip.05in \forall \alpha\beta$ and hence that 
$$c_{ijB,Q} - c_{ijB,Q}' = 0 \hskip.05in \forall B \in \on{Mix}_{{\bH},i}.$$
This finishes the proof of the uniqueness.

\vskip.1in

Now take $v_{ij,Q}' \in G_{1,Q}$ such that $F^{e_i}(v_{ij,Q}') = 
\overline{h_{ij,Q}} - \sum_{B \in \on{Mix}_{{\bH},i}}
c_{ijB,Q}\overline{H_Q}^B$.  

Setting $V_{i,Q}' = \left\{v_{ij,Q}'\right\}_j$, 
we see that we can replace
$V_{1,Q} \sqcup \cdots  \sqcup V_{M,Q}$ with 
$V_{1,Q}' \sqcup \cdots $
$\sqcup V_{i,Q}' \sqcup V_{i+1,Q} \sqcup \cdots \sqcup V_{M,Q}$ 
in the assertions of Lemma 3.1.2.1 in Part I.
\end{step}

This completes the proof for conclusion (c) by induction on $i$.

\vskip.1in

\begin{step}{Finishing argument to check condition (4).}

In order to check condition (4), it suffices to show that there exists 
$$\left\{g_{ijB}\right\}_{B \in \on{Mix}_{{\bH},i}} \subset R_f$$ 
such that 
$$g_{ijB}(Q) = c_{ijB,Q} \hskip.1in \forall Q \in C \cap \on{Supp}({\bI}) 
\cap {\fm}\text{-}\on{Spec}\hskip.02in R.$$
Fix a regular system of parameters $(x_1, \ldots, x_d)$ at $P$.  By 
shrinking $U_P$ if necessary, we may assume that $(x_1, \ldots, x_d)$ is a 
regular system of parameters over $U_P$, i.e., 
$(x_1 - x_1(Q),\ldots,x_d - x_d(Q))$ is a regular system of 
parameters at $Q$ for any 
$Q \in U_P \cap {\fm}\text{-}\on{Spec}\hskip.02in R$.

Now we analyze the condition of $\overline{h_{ij,Q}} - \sum_{B \in 
\on{Mix}_{{\bH},i}}c_{ijB,Q}\overline{H_Q}^B$ being pure, i.e., 
$$(\heartsuit) \hskip.1in \overline{h_{ij,Q}} - \sum_{B \in 
\on{Mix}_{{\bH},i}}c_{ijB,Q}\overline{H_Q}^B \in
L({\bI}_Q)_{p^{e_i}}^{\on{pure}}.$$
This happens if and only if, when we compute the power series expansions of 
$h_{ij}$ and $\sum_{B \in \on{Mix}_{{\bH},i}}c_{ijB,Q}H_Q^B$ with respect to 
the regular system of parameters 
$(x_1 - x_1(Q),\ldots, x_d - x_d(Q))$ 
and when we compare the degree $p^{e_i}$ terms, their 
mixed parts coincide (even though their pure parts may well not coincide).  
Since the coefficients of (the mixed parts
of) the power series can be computed using the partial derivatives with 
respect to $X = (x_1, \ldots, x_d)$, we conclude that condition 
$(\heartsuit)$ is equivalent to the following linear
equation 
$$(\heartsuit\heartsuit) \hskip.1in \left[\partial_{X^I}H^B(Q)\right]_{I \in 
\on{Mix}_{X,i}}^{B \in \on{Mix}_{{\bH},i}}\left[c_{ijB,Q}\right]_{B \in 
\on{Mix}_{{\bH},i}} =
\left[\partial_{X^I}h_{ij}(Q)\right]_{I \in \on{Mix}_{X,i}},$$
where
$$\on{Mix}_{X,i} = \left\{I = (i_1, \ldots, i_d) \mid |I| = p^{e_i}, i_l \neq 
p^{e_i} \hskip.05in \forall l = 1, \ldots, d\right\}$$
and where
\begin{displaymath}
\begin{array}{ll}
\left[\partial_{X^I}H^B(Q)\right]_{I \in \on{Mix}_{X,i}}^{B \in 
\on{Mix}_{{\bH},i}}& \on{is\ a\ matrix\ of\ size} \hskip.05in (\# 
\on{Mix}_{X,i}) \times (\# \on{Mix}_{{\bH},i}) \\
\left[c_{ijB,Q}\right]_{B \in \on{Mix}_{{\bH},i}}& \on{is\ a\ matrix\ of\ 
size} \hskip.05in (\# \on{Mix}_{{\bH},i}) \times 1, \hskip.05in \on{and} \\
\left[\partial_{X^I}h_{ij}(Q)\right]_{I \in \on{Mix}_{X,i}}& \on{is\ a\ 
matrix\ of\ size} \hskip.05in (\# \on{Mix}_{X,i}) \times 1.\\
\end{array}
\end{displaymath}

In particular, at the closed point $P$, we have the following linear 
equation
$$\left[\partial_{X^I}H^B(P)\right]_{I \in \on{Mix}_{X,i}}^{B \in 
\on{Mix}_{{\bH},i}}\left[c_{ijB,P}\right]_{B \in \on{Mix}_{{\bH},i}} =
\left[\partial_{X^I}h_{ij}(P)\right]_{I \in \on{Mix}_{X,i}}.$$
Since the solution $\left[c_{ijB,P}\right]_{B \in \on{Mix}_{{\bH},i}}$  
{\it uniquely} exists (cf. conclusion (c)), we conclude that the coefficient 
matrix of the linear equation has full rank,
i.e., 
$$\on{rank} \left[\partial_{X^I}H^B(P)\right]_{I \in \on{Mix}_{X,i}}^{B \in
\on{Mix}_{{\bH},i}} = \#\on{Mix}_{{\bH},i}.$$
Therefore, there exists a subset $S \subset \on{Mix}_{X,i}$ with $\# S = 
\#\on{Mix}_{{\bH},i}$ such that the corresponding minor has a nonzero 
determinant, i.e.,  
$$\det \left[\partial_{X^I}H^B(P)\right]_{I \in S}^{B \in
\on{Mix}_{{\bH},i}} \in k^{\times}.$$
Then the solution $\left[c_{ijB,P}\right]_{B \in \on{Mix}_{{\bH},i}}$ can be 
expressed as follows
$$\left[c_{ijB,P}\right]_{B \in \on{Mix}_{{\bH},i}} = 
\left(\left[\partial_{X^I}H^B(P)\right]_{I \in S}^{B \in
\on{Mix}_{{\bH},i}}\right)^{-1}\left[\partial_{X^I}h_{ij}(P)\right]_{I \in 
S}.$$
(Note that actually the matrix $\left[c_{ijB,P}\right]_{B 
\in \on{Mix}_{{\bH},i}}$ as well as the matrix 
$\left[\partial_{X^I}h_{ij}(P)\right]_{I \in 
S}$ is a zero matrix.)  By shrinking $U_P$ if necessary, we may assume
$$\det \left[\partial_{X^I}H^B\right]_{I \in S}^{B \in
\on{Mix}_{{\bH},i}} \in (R_f)^{\times}$$
and hence that
$$\det \left[\partial_{X^I}H^B(Q)\right]_{I \in S}^{B \in
\on{Mix}_{{\bH},i}} \in k^{\times} \hskip.1in \forall Q \in C \cap 
\on{Supp}({\bI}) \cap {\fm}\text{-}\on{Spec}\hskip.02in R.$$
Then the solution $\left[c_{ijB,Q}\right]_{B \in \on{Mix}_{{\bH},i}}$ for 
$(\heartsuit\heartsuit)$ can be expressed as follows
$$\left[c_{ijB,Q}\right]_{B \in \on{Mix}_{{\bH},i}} = 
\left(\left[\partial_{X^I}H^B(Q)\right]_{I \in S}^{B \in
\on{Mix}_{{\bH},i}}\right)^{-1}\left[\partial_{X^I}h_{ij}(Q)\right]_{I \in 
S}.$$
It follows immediately from this that, if we define the set $\{g_{ijB}\}_{B 
\in \on{Mix}_{{\bH},i}}$ by the formula
$$\left[g_{ijB}\right]_{B \in \on{Mix}_{{\bH},i}} = 
\left(\left[\partial_{X^I}H^B\right]_{I \in S}^{B \in
\on{Mix}_{{\bH},i}}\right)^{-1}\left[\partial_{X^I}h_{ij}\right]_{I \in 
S},$$
then it satisfies the desired condition
$$g_{ijB}(Q) = c_{ijB,Q} \hskip.1in \forall Q \in C \cap \on{Supp}({\bI}) 
\cap {\fm}\text{-}\on{Spec}\hskip.02in R.$$
Finally, by shrinking $U_P$ if necessary so that the above argument is valid 
for any element $h_{ij}$ taken from the given leading generator system 
${\bH}$, we see that condition (4) is satisfied.
\end{step}

This completes the proof of Proposition \ref{1.3.3.3}.

\end{proof}

\end{subsection}

\end{section}

\end{chapter} 

\begin{chapter}{Power series expansion}
As in Chapter 1, we denote by $R$ the coordinate ring of an affine open
subset $\on{Spec}\hskip.02in R$ of a nonsingular variety $W$ of $\dim W = d$
over an
algebraically closed  field $k$ of positive characteristic ${\on{char}}(k) =
p$ or of characteristic zero ${\on{char}}(k) = 0$, where in the latter case
we 
formally set 
$p = \infty$ (cf. 0.2.3.2.1 and Definition 3.1.1.1 (2) in Part I).

We fix a closed point $P \in W$.

Let ${\bI}_P$ be a ${\fD}$-saturated idealistic filtration over $R_P =
{\mathcal O}_{W,P}$, the local ring at the closed point, with ${\fm}_P$
being its maximal ideal.

Let ${\bH} = \left\{(h_l,p^{e_l})\right\}_{l = 1}^N$ be a leading
generator system of ${\bI}_P$.

In characteristic zero, the elements in the leading generator system are all
concentrated at level 1, i.e., $e_l = 0$ and $p^{e_l} = 1$ for $l = 1,
\ldots,
N$ (cf. Chapter 3 in Part I).  This implies by definition of a leading
generator system that the set of the elements $H = (h_l \mid l = 1, 
\ldots, N)$ forms (a part of) a regular system of parameters
$(x_1,
\ldots, x_d)$.  (Say $h_l = x_l$ for
$l = 1,
\ldots, N$.)  In positive characteristic, this is no longer the case.  However,
we can still regard the notion of a leading generator system as a
generalization of
the notion of a regular system of parameters, and we may expect some similar
properties between the two notions.

Now any element $f
\in R_P$ (or more generally any element $f \in \widehat{R_P}$) can be
expressed as a power series with respect to the regular system of
parameters and hence with respect to the leading generator system as above
in characteristic zero.  That is to say, we can write
$$f = \sum_{I \in ({\bZ}_{\geq 0})^d}c_IX^I = \sum_{B \in ({\bZ}_{\geq
0})^N} a_BH^B$$
where $c_I \in k$ and where $a_B$ is a power series in terms of the
remainder $(x_{N+1}, \ldots, x_d)$ of the regular system of
parameters.  

Chapter 2 is devoted to the study of the power series expansion with respect
to the elements in a
leading generator system (and its (weakly-)associated regular system of
parameters), one of the expected similar properties mentioned above, which is
valid both
in characteristic zero and in positive characteristic.

\begin{section}{Existence and uniqueness.}\label{2.1}
\end{section} 
\begin{subsection}{Setting for the power series expansion.}\label{2.1.1}
First we describe the setting for 
\linebreak[4]
Chapter 2, which is slightly more
general than just dealing with a leading generator system.  Actually, until
we reach \S 2.2, our argument does {\it not} involve the notion of an idealistic
filtration.

Let ${\mathcal H} = \{h_1, \ldots, h_N\} \subset R_P$ be a subset consisting
of $N$ elements, and 
\text{$0\leq e_1\leq\cdots \leq e_N$} 
nonnegative integers attached to these elements,
satisfying the following
conditions (cf. 4.1.1 in Part I):
\begin{center} 
\begin{tabular}{ll}
(i) & $h_l \in {\fm}_P^{p^{e_l}}$ and $\overline{h_l} = (h_l 
\bmod  {\fm}_P^{p^{e_l}+1}) = v_l^{p^{e_l}}$ with
$v_l
\in {\fm}_P/{\fm}_P^2$ for $l = 1, \ldots, N$,
\\ 
(ii) & $\left\{v_l \mid l = 1, \ldots, N\right\} 
\subset {\fm}_P/{\fm}_P^2$ consists of
$N$-distinct and $k$-linearly independent \\
&elements in the
$k$-vector space ${\fm}_P/{\fm}_P^2$.
\\ 
\end{tabular} 
\end{center} 
We also take a regular system of parameters $(x_1, \ldots, x_d)$ such that
$$(\on{asc}) \hskip.1in v_l = \overline{x_l} = (x_l \bmod 
{\fm}_P^2) \text{\ for\ } l = 1, \ldots, N.$$
We say $(x_1, \ldots, x_d)$ is associated to $H = (h_1, \ldots, h_N)$ if
the above condition $(\on{asc})$ is satisfied.
\end{subsection}

\begin{subsection}{Existence and uniqueness of the power series
expansion.}\label{2.1.2}
\begin{lem}\label{2.1.2.1} Let the setting be as described in \ref{2.1.1}.
Then any element $f \in \widehat{R_P}$ has a power series expansion, with
respect to $H = (h_1, \ldots, h_N)$ and 
its associated regular system of parameters
$(x_1,
\ldots, x_d)$, of the form
$$(\star) \hskip.1in f = \sum_{B \in ({\bZ}_{\geq 0})^N}a_BH^B \text{\ 
where\ }a_B = \sum_{K \in ({\bZ}_{\geq 0})^d}b_{B,K}X^K,$$
with $b_{B,K}$ being a power series in terms of the remainder
$(x_{N+1},
\ldots, x_d)$ of the regular system of parameters, and with $K = (k_1,
\ldots, k_d)$ varying in the range satisfying the condition
$$0 \leq k_l \leq p^{e_l} - 1 \text{\ for\ }l = 1, \ldots, N \hskip.05in
\on{and} \hskip.05in k_l = 0 \text{\ for\ }l = N+1, \ldots, d.$$
Moreover, the power series expansion of the form $(\star)$ is unique.
\end{lem}

\begin{proof} ($\mathbf{Existence}$) We construct a sequence 
$\left\{f_r\right\}_{r \in {\bZ}_{\geq 0}} \subset R_P$ 
in the following manner.

\begin{case}{Construction of $f_0$.}

In this case, choose $f_0 = a_{{\mathbb O},0} \in k$ such that

\begin{center}
\begin{tabular}{ll}
(i) & $f - f_0 \in \widehat{{\fm}_P}^1$, \\
(ii) & $f_0 = \sum_{|[B]| \leq 0}a_{B,0}H^B$. \\
\end{tabular}
\end{center}
\end{case}

\vskip.1in

\begin{case}{Construction of $f_{r+1}$ assuming that of $f_r$.}

Suppose inductively that we have constructed $f_r$ satisfying the following
conditions:

\begin{center}
\begin{tabular}{ll}
(i) & $f - f_r \in \widehat{{\fm}_P}^{r+1}$, \\
(ii) & $f_r = \sum_{|[B]| \leq r}a_{B,r}H^B$ where $a_{B,r} = \sum
b_{B,K,r}X^K$ \\
& with $b_{B,K,r}$ being a polynomial in $(x_{N+1}, \ldots, x_d)$, i.e.,\\
& $b_{B,K,r} = \sum_{|[B]| + |K| + |J| \leq r} c_{B,K,J}X^J$ where
$c_{B,K,J} \in k$ and \\
&where $J = (j_1, \ldots, j_d)$ with $j_l = 0$ for $l = 1, \ldots, N$,\\
&and $K$ varying in the range specified above, and \\
&satisfying the condition $|[B]| + |K| + |J| \leq r$. \\
\end{tabular}
\end{center}

Now express $f - f_r = \sum c_{I,r}X^I$ with $c_{I,r} \in k$ as a power
series expansion in terms of the regular system of parameters $X =
(x_1,
\ldots, x_d)$.  

Given $I = (i_1, \ldots, i_d)$ with $|I| = r+1$, determine

\begin{displaymath}
\left\{\begin{array}{ll}
B &= (b_1, \ldots, b_N), \\
K &= (k_1, \ldots, k_N,\ \ 0,\ \ \ldots\ ,\ \ 0) \in ({\bZ}_{\geq 0})^d, \\
J &= (0,\ \ \ldots,\ \ 0, j_{N+1}, \ldots, j_d) \in ({\bZ}_{\geq 0})^d \\
\end{array}\right.
\end{displaymath}
by the formulas below
\begin{displaymath}
\left\{\begin{array}{lll}
i_l &= b_lp^{e_l} + k_l \hskip.1in \on{with} \hskip.1in b_l \in {\bZ}_{\geq
0} \hskip.1in \on{and} \hskip.1in 0 \leq k_l \leq p^{e_l} - 1
&\on{for} \hskip.1in l = 1,
\ldots, N
\\
i_l &= j_l &\on{for} \hskip.1in l = N+1, \ldots, d.\\
\end{array}\right.
\end{displaymath}
Then it is straightforward to see, after renaming $c_{I,r}$ as $c_{B,K,J}$,
that the following equality holds
$$\sum_{|I| = r+1}c_{I,r}X^I = \sum_{|[B]| + |K| + |J| =
r+1}c_{B,K,J}X^JX^KH^B \bmod  {\fm}_P^{r+2}.$$
Set
$$\left\{\begin{array}{ll}
b_{B,K,r+1} &= \sum_{|[B]| + |K| + |J| \leq r+1}c_{B,K,J}X^J \\
a_{B,r+1} &= \sum b_{B,K,r+1}X^K \\
f_{r+1} &= \sum_{|[B]| \leq r+1}a_{B,r+1}H^B
\end{array}\right.$$

Then $f_{r+1}$ clearly satisfies conditions (i) and (ii).

This finishes the inductive construction of the sequence 
$\left\{f_r\right\}_{r \in {\bZ}_{\geq 0}} \subset R_P$.
\end{case}

\vskip.1in

Now set
\begin{displaymath}
\left\{\begin{array}{ll}
b_{B,K} &= \lim_{r \rightarrow \infty}b_{B,K,r} = \sum c_{B,K,J}X^J \\
a_B &= \lim_{r \rightarrow \infty}a_{B,r} = \sum b_{B,K}X^K, \\
\end{array}\right.
\end{displaymath}
where each of the above limits exists by condition (ii).

Then condition (i) implies
$$f = \lim_{r \rightarrow \infty}f_r = \lim_{r \rightarrow
\infty}\sum_{|[B]| \leq r}a_{B,r}H^B = \sum a_BH^B,$$
proving the existence of a power series expansion of the form $(\star)$.

\vskip.1in

($\mathbf{Uniqueness}$) In order to show the uniqueness of the power 
series expansion of the form $(\star)$, we have only to verify
$$0 = \sum_{B \in ({\bZ}_{\geq 0})^N}a_BH^B \hskip.1in \on{of\ the\ 
form\ }(\star) \Longleftrightarrow a_B = 0 \hskip.1in \forall B \in 
({\bZ}_{\geq 0})^N.$$
As the implication $(\Longleftarrow)$ is obvious, we show the opposite
implication $(\Longrightarrow)$ in what follows.

Suppose $0 = \sum_{B \in ({\bZ}_{\geq 0})^N}a_BH^B$.

Assume that there exists $B \in ({\bZ}_{\geq 0})^N$ such that $a_B \neq 0$.

Set $s = \min\left\{\on{ord}\left(a_BH^B\right) \mid a_B \neq 0\right\}$.

Write
$$a_B = \sum_{K \in ({\bZ}_{\geq 0})^d}b_{B,K}X^K \hskip.1in \on{and}
\hskip.1in b_{B,K} = \sum_{J \in ({\bZ}_{\geq 0})^d} c_{B,K,J}X^J
\hskip.1in \on{with} \hskip.1in c_{B,K,J} \in k,$$
where $K = (k_1, \ldots, k_d)$ varies in the range satisfying the condition
$$0 \leq k_l \leq p^{e_l} - 1 \text{\ for\ }l = 1, \ldots, N \hskip.05in
\on{and} \hskip.05in k_l = 0 \text{\ for\ }l = N+1, \ldots, d,$$
and where $J = (j_1, \ldots, j_d)$ varies in the range satisfying the
condition
$$j_l = 0 \hskip.05in \text{\ for\ } \hskip.05in l = 1, \ldots, N.$$
Then we have
\begin{eqnarray*}
0 &=& \sum_Ba_BH^B = \sum_B\sum_K\left(\sum_J c_{B,K,J}X^JX^K\right)H^B\\
&=& \sum_{|[B]| + |K| + |J| = s} c_{B,K,J}X^JX^K\left(\prod_{l = 1}^N
x_l^{p^{e_l}b_l}\right) \hskip.3in \bmod 
\widehat{{\fm}_P}^{s+1}.
\end{eqnarray*}

On the other hand, we observe that the set of 
all the monomials of degree $s$ 
\linebreak[4]
$\{X^JX^K\left(\prod_{l = 1}^N
x_l^{p^{e_l}b_l}\right)\}_{|[B]| + |K| + |J| = s} =
\{X^I\}_{|I| = s}$ obviously forms a basis of the vector space
$\widehat{{\fm}_P}^s/\widehat{{\fm}_P}^{s+1}$, and that $c_{B,K,J} \neq 0$
for
some
$B, K, J$ with
$|[B]| + |K| + |J| = s$ by the assumption and by the choice of $s$.

This is a contradiction !

Therefore, we conclude that $a_B = 0 \hskip.1in \forall B \in ({\bZ}_{\geq
0})^N$.

This finishes the proof of the implication $(\Longrightarrow)$, and hence
the proof of the uniqueness of the power series expansion of the form
$(\star)$.

This completes the proof of Lemma \ref{2.1.2.1}.
\end{proof}

\begin{rem}\label{2.1.2.2}

\begin{enumerate}

\item It follows immediately from the argument to show the existence and
uniqueness of the power series expansion $f = \sum a_BH^B$ of the form
$(\star)$ that
$$\on{ord}(f) = \min\left\{\on{ord}\left(a_BH^B\right)\right\} =
\min\left\{\on{ord}\left(a_B\right) + |[B]|\right\}$$
and hence that
$$\on{ord}(a_B) \geq \on{ord}(f) - |[B]| \hskip.1in \forall B \in
({\bZ}_{\geq 0})^N.$$

\item In the setting \ref{2.1.1}, we defined the notion of a regular system
associated to $H = (h_1, \ldots, h_N)$.
We say that a regular system of parameters $(x_1, \ldots, x_d)$ is 
{\it weakly-associated} to $H = (h_1, \ldots, h_N)$, if the following
condition holds:
$$\det\left[\partial_{x_i^{p^e}}(h_l^{p^{e - e_l}})\right]_{i = 1, \ldots,
L_e}^{l = 1, \ldots, L_e} \in R_P^{\times} \hskip.05in \on{for} \hskip.05in
e = e_1,
\ldots, e_N \hskip.05in \on{where}
\hskip.05in L_e =
\#\left\{l
\mid e_l
\leq e\right\}.$$

All the assertions of Lemma \ref{2.1.2.1} hold, even if we only require a
regular system of parameters $(x_1, \ldots, x_d)$ to be 
weakly-associated to $H$, instead of associated to $H$.
\end{enumerate}
\end{rem}

\end{subsection}

\begin{section}{Formal coefficient lemma.}\label{2.2}

\begin{subsection}{Setting for the formal coefficient lemma.}\label{2.2.1}
As we can see from the description of the
setting \ref{2.1.1}, our discussion on the power series expansion of the
form
$(\star)$ (cf. Lemma \ref{2.1.2.1}) so far does not involve the notion 
of an idealistic filtration.  However, the most interesting and 
important result of Chapter 2 is obtained as
Lemma \ref{2.2.2.1} below, which we call the formal coefficient lemma, when
we get the notion of an idealistic filtration involved and impose an extra
condition related to it as follows:

Let ${\mathcal H} = \{h_1, \ldots, h_N\} \subset R_P$ be a subset consisting
of $N$ elements, and \text{$0 \leq e_1 \leq \cdots \leq e_N$}
nonnegative integers attached to these elements, satisfying conditions (i)
and (ii), as described in the setting \ref{2.1.1}.  Let $X = (x_1,
\ldots, x_d)$ be a regular system of parameters associated to $H = (h_1,
\ldots, h_N)$ with $h_l = x_l^{p^{e_l}}\ \on{mod}\ {\fm}_P^{p^{e_l} 
+  1}$ for $l = 1, \cdots, N$.  

Let ${\bI}_P$ be a ${\fD}$-saturated idealistic filtration over $R_P$.

We impose the following extra condition

\begin{center} 
\begin{tabular}{ll}
(iii) & $(h_l, p^{e_l}) \in {\bI}_P$ for $l = 1, \ldots, N$.
\\ 
\end{tabular} 
\end{center} 
\end{subsection}

\begin{subsection}{Statement of the formal coefficient lemma and its
proof.}\label{2.2.2} Now our assertion is that, under the setting of
\ref{2.2.1} and given an element in the idealistic filtration, the
coefficients of the power series expansion of the form $(\star)$, with
``appropriate'' levels attached, belongs to (the completion of) 
the idealistic filtration.  We formulate this
assertion as the following formal coefficient lemma.

\begin{lem}\label{2.2.2.1} 
Let the setting be as described in \ref{2.2.1}.
Let $\widehat{{\bI}_P}$ be the completion of 
the idealistic filtration ${\bI}_P$ (cf. \S 2.4 in Part I).

Take an element $(f,a) \in \widehat{{\bI}_P}$.  

Let $f = \sum_{B \in ({\bZ}_{\geq 0})^N}a_BH^B$ be the power series
expansion of the form $(\star)$ (cf. Lemma \ref{2.1.2.1}).  Then we have 
$$(a_B,a - |[B]|) \in \widehat{{\bI}_P} \hskip.1in \forall B \in
({\bZ}_{\geq 0})^N.$$
\end{lem}

\begin{proof} We will derive a contradiction assuming
$$(a_B,a - |[B]|) \not\in \widehat{{\bI}_P} \hskip.1in \on{for\ some}
\hskip.1in B \in ({\bZ}_{\geq 0})^N.$$
Note that, under the assumption, there should exist $B \in
({\bZ}_{\geq
0})^N$ with $B \neq {\mathbb O}$ such that 
$(a_B,a - |[B]|) \not\in \widehat{{\bI}_P}$.
(In fact, suppose $(a_B,a - |[B]|) \in 
\widehat{{\bI}_P} \hskip.1in \forall B \neq {\mathbb O}$.  
Then the equality $a_{\mathbb O} = f - \sum_{B \neq {\mathbb
O}}a_BH^B$ and the inclusions $(f,a) \in \widehat{{\bI}_P}$ and $(a_BH^B,a)
\in \widehat{{\bI}_P} \hskip.05in \forall B \neq {\mathbb O}$, would
imply
$(a_{\mathbb O},a) = (a_{\mathbb O},a - |[{\mathbb O}]|) \in
\widehat{{\bI}_P}$, which is against the assumption.)

\vskip.1in

We introduce the following notations:

\vskip.05in

$$\left\{\begin{array}{ll}
l_B &= |[B]| + \sup\left\{n \in {\bZ}_{\geq 0} \mid a_B \in
(\widehat{{\bI}_P})_{a - |[B]|} + \widehat{{\fm}_P}^n\right\}  \hskip.1in
\on{for\ } B \in ({\bZ}_{\geq0})^N \setminus\{{\mathbb O}\},\\
l &= \min_{B \in ({\bZ}_{\geq0})^N, B \neq {\mathbb O}}\left\{l_B\right\}, \\
\Gamma_B &= (\widehat{{\bI}_P})_{a - |[B]|} + \widehat{{\fm}_P}^{l - |[B]| +
1} \hskip.1in \on{for\ } B \in ({\bZ}_{\geq0})^N, \\
L_B &= \max\left\{B + K \mid a_B \in \Gamma_B + \sum_{K
\leq M}\widehat{{\fm}_P}^{l - |[B + M]|}H^M\right\} 
\on{\ for\ }B \in ({\bZ}_{\geq0})^N\setminus\{{\mathbb O}\}, l_B = l,\\
L &= \min_{B \in ({\bZ}_{\geq
0})^N, B \neq {\mathbb O}, l_B = l}\left\{L_B\right\}, \\
B_o &= \max_{B \in ({\bZ}_{\geq0})^N, B \neq {\mathbb O}, 
l_B = l, L_B = L}\left\{B\right\} \\
\Lambda_B &= \Gamma_B + \sum_{L < B + M} \widehat{{\fm}_P}^{l - |[B +
M]|}H^M \on{\ for\ }B \in ({\bZ}_{\geq0})^N. \\
\end{array}
\right.$$

\vskip.05in

Note that $l < \infty$ by the assumption $a_B \not\in (\widehat{{\bI}_P})_{a
- |[B]|}$ for some $B \neq {\mathbb O}$.  Note that the maximum of 
$B + K$, the minimum of
$L_B$, and the maximum of $B$ are all taken with respect to the
lexicographical order on
$({\bZ}_{\geq 0})^N$.  Note that, if $l_B = l$, then $[B] < a$.  This
guarantees the existence of the maximum of $B \in ({\bZ}_{\geq 0})^N$ with
$B \neq {\mathbb O},
l_B = l, L_B = L$.  We remark that, when
$r
\leq 0$, we understand by convention
$\widehat{{\fm}_P}^r$  represents
$\widehat{R_P}$. 

We claim, for $B, K \in ({\bZ}_{\geq 0})^N$,

\hskip1.7in (i) $H^K\Lambda_{B+K}\subset \Lambda_B$,

\hskip1.7in (ii) $\partial_{[K]}(\Lambda_B) \subset \Lambda_{B + K}$.

\vskip.1in

(We remark that we identify $[K]$, for 
$K = (k_1, \ldots, k_N) \in ({\bZ}_{\geq0})^N$, 
with 
$(p^{e_1}k_1, \ldots,$ 
$p^{e_N}k_N,0, \ldots, 0) \in({\bZ}_{\geq0})^d$, 
and hence that we understand
$\partial_{[K]}$ denotes  $\partial_{X^{[K]}} =
\partial_{x_1^{p^{e_1}k_1} \cdots x_N^{p^{e_N}k_N}}$ in claim (ii).)

In fact, since $(H^K,|[K]|) \in \widehat{{\bI}_P}$ and since $H^K \in
\widehat{{\fm}_P}^{|[K]|}$, we see
\begin{eqnarray*}
H^K\Lambda_{B+K} &=& H^K\left(\Gamma_{B+K} + \sum_{L < B + K + M}
\widehat{{\fm}_P}^{l - |[B + K + M]|}H^M\right) \\
&=& H^K\left((\widehat{{\bI}_P})_{a - |[B + K]|} + \widehat{{\fm}_P}^{l -
|[B + K]| + 1} + \sum_{L < B + K + M} \widehat{{\fm}_P}^{l - |[B + K
+ M]|}H^M\right) \\
&\subset& (\widehat{{\bI}_P})_{a - |[B]|} + \widehat{{\fm}_P}^{l - |[B]| +
1} + \sum_{L < B + M} \widehat{{\fm}_P}^{l - |[B + M]|}H^M \\
&&\hskip1.6in (\on{by\
replacing\ old\ }M + K \on{\ with\ new\ }M)\\
&=& \Gamma_B + \sum_{L < B + M} \widehat{{\fm}_P}^{l - |[B + M]|}H^M =
\Lambda_B, \\
\end{eqnarray*}
checking claim (i).

\newpage

In order to see claim (ii), observe

\begin{itemize}

\item $\partial_{[K]}\left((\widehat{{\bI}_P})_{a - |[B]|}\right) \subset
(\widehat{{\bI}_P})_{a - |[B + K]|}$, since $\widehat{{\bI}_P}$ is
${\fD}$-saturated,

\item $\partial_{[K]}\left(\widehat{{\fm}_P}^{l - |[B]| + 1}\right) \subset
\widehat{{\fm}_P}^{l - |[B + K]| + 1}$, and

$\partial_{[K] - I}\left(\widehat{{\fm}_P}^{l -
|[B + M]|}\right) \subset \widehat{{\fm}_P}^{l - |[B + K + M]| + |I|}$ 
for $I$ with $I \leq [K]$,

\item $\partial_I(H^M) \subset \binom{[M]}{I}H^{M - I} +
\widehat{{\fm}_P}^{|[M]| - |I| + 1}$, and

$\binom{[M]}{I} = 0$ unless $I = [J]$ for some $J \in ({\bZ}_{\geq 0})^N$.
\end{itemize}

Using these observations, we compute

\begin{eqnarray*}
\lefteqn{
\partial_{[K]}\left(\Lambda_B\right) = \partial_{[K]}\left(\Gamma_B +
\sum_{L < B + M} \widehat{{\fm}_P}^{l - |[B + M]|}H^M\right)}\\
&=& \partial_{[K]}\left((\widehat{{\bI}_P})_{a - |[B]|} +
\widehat{{\fm}_P}^{l - |[B]| + 1} + \sum_{L < B + M} \widehat{{\fm}_P}^{l -
|[B + M]|}H^M\right) \\
&=& \partial_{[K]}\left((\widehat{{\bI}_P})_{a - |[B]|}\right) +
\partial_{[K]}\left(\widehat{{\fm}_P}^{l - |[B]| + 1}\right) + \sum_{L < B +
M}
\partial_{[K]}\left(\widehat{{\fm}_P}^{l - |[B + M]|}H^M\right) \\
&=& \partial_{[K]}\left((\widehat{{\bI}_P})_{a - |[B]|}\right) +
\partial_{[K]}\left(\widehat{{\fm}_P}^{l - |[B]| + 1}\right) \\
&&+ \sum_{L < B +
M}\left[\sum_{I \leq [K]}\partial_{[K] - I}\left(\widehat{{\fm}_P}^{l - |[B
+ M]|}\right)\partial_I\left(H^M\right)\right] \\
&&\hskip.5in (\on{by\ the\ generalized\ product\ rule\ 
(cf.\ Lemma\ 1.2.1.2\ (3)\ in\ Part\ I)}) \\
&=& \partial_{[K]}\left((\widehat{{\bI}_P})_{a - |[B]|}\right) +
\partial_{[K]}\left(\widehat{{\fm}_P}^{l - |[B]| + 1}\right) + \sum_{L < B +
M}\\
&&\left[\sum_{I = [J], I \leq [K]} \partial_{[K] -
I}\left(\widehat{{\fm}_P}^{l - |[B + M]|}\right)\partial_I\left(H^M\right) +
\sum_{I \neq
[J], I \leq [K]}\partial_{[K] - I}\left(\widehat{{\fm}_P}^{l - |[B +
M]|}\right)\partial_I\left(H^M\right)\right]\\
&\subset& (\widehat{{\bI}_P})_{a - |[B + K]|} + \widehat{{\fm}_P}^{l - |[B +
K]| + 1}     
+ \sum_{L < B + M} \left[\sum_{I = [J], J \leq K, J \leq M}
\widehat{{\fm}_P}^{l - |[B + M + K - J]|}H^{M-J}\right]  \\
&=& \Gamma_{B + K} + \sum_{L < B + M + (K - J) = B + K + (M - J), 
J \leq K, J \leq M}
\widehat{{\fm}_P}^{l - |[B + K + M -
J]|}H^{M-J}
\\
&\subset& \Gamma_{B + K} + \sum_{L < B + K + M}\widehat{{\fm}_P}^{l - |[B +
K
+ M]|}H^M = \Lambda_{B + K} \\
&& \hskip2in (\on{by\
replacing\ old\ }M - J \on{\ with\ new\ }M),\\
\end{eqnarray*}
checking claim (ii).

\vskip.1in

Now by definition, for each $B \in ({\bZ}_{\geq 0})^N$ with $B \neq {\mathbb
O}, l_B = l, L_B = L$, we can choose $b_B \in
\widehat{{\fm}_P}^{l - |[L]|}$ such that $a_B - b_BH^{L-B} \in \Lambda_B$.
For each $B \in ({\bZ}_{\geq 0})^N$ with $B
\neq {\mathbb O}$ but $l_B \neq l$ or $L_B \neq L$, we set $b_B = 0$ and have
$a_B - b_BH^{L-B} \in \Lambda_B$.

Therefore, we have, for each $B \in ({\bZ}_{\geq 0})^N$ with $B
\neq {\mathbb O}$,
$$a_B - b_BH^{L-B} \in \Lambda_B$$
and\ hence by claim (i) (with $B, K \in ({\bZ}_{\geq 0})^N$ there being
equal to ${\mathbb O}, B$ below, respectively)
$$(a_B - b_BH^{L-B})H^B \in \Lambda_{\mathbb O}.$$


Now we compute (with the symbol ``$\equiv$'' denoting the equality modulo
$\Lambda_{B_o}$): 
\begin{eqnarray*}
\lefteqn{
\partial_{[B_o]}f = \partial_{[B_o]}\left(\sum a_BH^B\right) =
\partial_{[B_o]}\left(\sum_{B \neq {\mathbb O}} a_BH^B\right) 
\equiv \partial_{[B_o]}\left(\sum_{B \neq {\mathbb O}}b_BH^L\right)
}\\
&& \hskip.3in (\on{since\ }\sum_{B \neq {\mathbb O}}a_BH^B - 
\sum_{B \neq {\mathbb O}}b_BH^L \in \Lambda_{\mathbb O} \hskip.1in 
\on{and\ by\ claim\ (ii)}) \\ &=& \sum_{B \neq {\mathbb O}, l_B = l, 
L_B = L}\partial_{[B_o]}\left(b_BH^{L-B}H^B\right) \\
&=& \sum_{B \neq {\mathbb O},l_B = l, L_B = L}\left[\sum_{I \leq
[B_o]}\partial_I\left(b_BH^{L-B}\right)\partial_{[B_o]-I}\left(H^B\right)
\right]  
\\ 
&& \hskip1.5in (\on{by\ the\ generalized\ product\ rule}) \\
&\equiv& \sum_{B \neq {\mathbb O},l_B = l, L_B = L} \left[\sum_{I \leq
[B_o], I = [K]}\partial_{[K]}\left(b_BH^{L-B}\right)\partial_{[B_o -
K]}\left(H^B\right)\right]
\\
&& \hskip.5in (\on{refer\ to\ the\ last\ observation\ 
used\ to\ see\ claim\ (ii)})
\\
&\equiv& \sum_{B \neq {\mathbb O},l_B = l, L_B = L}
b_BH^{L-B}\partial_{[B_o]}\left(H^B\right)
\\
&& \hskip.3in (\on{since\ for\ }K \neq {\mathbb O}\on{\ we\ have\ }
\partial_{[K]}\left(b_BH^{L-B}\right) = \partial_{[K]}\left(- \left(a_B -
b_BH^{L -
B}\right)\right)
\in \Lambda_{B + K} \\
&& \hskip.3in \on{and\ hence\
}\partial_{[K]}\left(b_BH^{L-B}\right)\partial_{[B_o - K]}\left(H^B\right)
\in \Lambda_{B + K}\partial_{[B_o -
K]}\left(H^B\right) \subset \Lambda_{B_o})\\ &\equiv& \sum_{B \neq {\mathbb
O},l_B = l, L_B = L}\binom{B}{B_o}b_BH^{L - B_o} \\
&=& b_{B_o}H^{L - B_o} \hskip1.3in (\on{by\ the\ maximality\ of\ }B_o).\\
\end{eqnarray*}
Note that the inclusion $\Lambda_{B + K}\partial_{[B_o -
K]}\left(H^B\right) \subset \Lambda_{B_o}$ used above is verified as
follows:

\begin{eqnarray*}
\lefteqn{
\Lambda_{B + K}\partial_{[B_o - K]}\left(H^B\right)} \\
&=& \left(\Gamma_{B + K} + \sum_{L < B + K + M}\widehat{{\fm}_P}^{l - |[B +
K + M]|}H^M\right)\partial_{[B_o -
K]}\left(H^B\right) \\
&=& \left((\widehat{{\bI}_P})_{a - |[B + K]|} + \widehat{{\fm}_P}^{l - |[B +
K]| + 1} + \sum_{L < B + K + M}\widehat{{\fm}_P}^{l - |[B + K +
M]|}H^M\right)\partial_{[B_o -
K]}\left(H^B\right) \\
&\subset& (\widehat{{\bI}_P})_{a - |[B_o]|} + \widehat{{\fm}_P}^{l - |[B_o]|
+ 1} \\
&&+ \left(\sum_{L < B + K + M}\widehat{{\fm}_P}^{l - |[B + K +
M]|}H^M\right)\left(\binom{[B]}{[B_o - K]}H^{B - (B_o - K)} +
\widehat{{\fm}_P}^{[B] - [B_o - K] + 1}\right) \\
&& (\on{since\ }\partial_{[B_o - K]}\left(H^B\right) \in
(\widehat{{\bI}_P})_{|[B]| - |[B_o - K]|} \on{\ and\ since\ }\partial_{[B_o
- K]}\left(H^B\right) \in
\widehat{{\fm}_P}^{|[B]| - |[B_o - K]|})
\\
&&(\on{refer\ also\ to\ the\ last\ observation\ used\ to\ see\ claim\ (ii)})
\\
&\subset& (\widehat{{\bI}_P})_{a - |[B_o]|} + \widehat{{\fm}_P}^{l - |[B_o]|
+ 1} \\
&&+ \sum_{L < B_o + (M + B + K - B_o)}\widehat{{\fm}_P}^{l - |[B_o + (M + B
+ K -
B_o)]|}H^{M + B + K - B_o} \\
\end{eqnarray*}

\begin{eqnarray*}
&\subset& (\widehat{{\bI}_P})_{a - |[B_o]|} + \widehat{{\fm}_P}^{l - |[B_o]|
+ 1} + \sum_{L < B_o + M}\widehat{{\fm}_P}^{l - |[B_o + M]|}H^M \\
&=& \Gamma_{B_o} + \sum_{L < B_o + M}\widehat{{\fm}_P}^{l - |[B_o + M]|}H^M
= \Lambda_{B_o}.\\
\end{eqnarray*} 

However, since $\partial_{[B_o]}f \in (\widehat{{\bI}_P})_{a - |[B_o]|}
\subset \Lambda_{B_o}$, we conclude
$$b_{B_o}H^{L - B_o} \in \Lambda_{B_o} \hskip.1in \on{and\ hence} \hskip.1in
a_{B_o} \in \Lambda_{B_o},$$
which contradicts the choice of $B_o$ with $L_{B_o} = L$.

This finishes the proof of Lemma \ref{2.2.2.1}.

\end{proof} 

We would like to remark that the same idea of the proof by
contradiction above actually leads to an explicit construction 
of the coefficients
through differential operators and taking limits.  We present such a
construction, which is of interest on its own and which is slightly
different from the direct
translation of the proof by contradiction above, as an alternative 
proof.  This alternative proof leads to an improved version of the 
formal coefficient lemma, which we will state and use in the 
subsequent papers to show that the new invariant $\widetilde{\nu}$ 
we introduce (cf. \ref{0.3.1}) is well-defined.

\begin{proof}[Alternative Proof]

\begin{step}{We show the statement when $B = {\mathbb O}$, i.e., we show
$(a_{\mathbb O},a - |[{\mathbb O}]|) = (a_{\mathbb O},a) \in 
\widehat{{\bI}_P}$.}

Let $g = \sum_{B \in ({\bZ}_{\geq 0})^N}a_{B,g}H^B$ be the power series
expansion of the form $(\star)$ for $g \in \widehat{R_P}$.  Note that we
add the subscript ``${}_g$'' to the notation for the coefficients $a_{B,g}$,
to emphasize their dependence on $g$ and in particular to distinguish
them from the coefficients $a_B = a_{B,f}$ for $f$.

Observe first (cf. Remark \ref{2.1.2.2} (1)) that, for any $B \in
({\bZ}_{\geq 0})^N$, we
have
$$\on{ord}(a_{B,g}H^B) \geq \on{ord}(g) \hskip.05in \on{and\ hence}
\hskip.05in \on{ord}(a_{B,g}) \geq \on{ord}(g) - |[B]|.$$
Given $g \in
\widehat{R_P}$, we define the invariant $\eta(g)$ by the formula
$$\eta(g) := \left(\on{ord}(g), \min\left\{B \mid \on{ord}(a_{B,g}H^B) =
\on{ord}(g)\right\}\right) \in {\bZ}_{\geq 0} \times ({\bZ}_{\geq 0})^N,$$
where the minimum in the second factor is taken with respect to the
lexicographical order on $({\bZ}_{\geq 0})^N$.  The values of the invariant
$\eta$ are ordered according to the lexicographical order given to
${\bZ}_{\geq 0} \times ({\bZ}_{\geq 0})^N$.

We will construct a sequence $\left\{g_n\right\}_{n \in {\bZ}_{\geq 0}} 
\subset \widehat{R_P}$ inductively, satisfying the following conditions:
\begin{center} 
\begin{tabular}{ll}
$(0)_n$ & $g_n \in ({\mathcal H})$,\\
$(1)_n$ & $(a_{\mathbb O} + g_n, a) \in \widehat{{\bI}_P}$,
\\ 
$(2)_n$ & $\eta(g_{n-1}) < \eta(g_n)$.\\
\end{tabular} 
\end{center} 

The construction of such a sequence is sufficient to prove $(a_{\mathbb O},a)
\in \widehat{{\bI}_P}$.

In fact, since there are only finitely many $B$'s with $\on{ord}(H^B) \leq
\nu$ for a fixed $\nu \in {\bZ}_{\geq 0}$, and since
$$\left\{B \mid \on{ord}(H^B) \leq \nu\right\} \supset \left\{B \mid 
\on{ord}(a_{B,g}H^B) = \on{ord}(g)\right\}$$
for any $g \in \widehat{R_P}$ with $\on{ord}(g) = \nu$, we conclude by
condition $(2)_n$ that
$$\lim_{n \rightarrow \infty}\on{ord}(g_n) = \infty \hskip.1in 
\on{and\ hence} \hskip.1in \lim_{n \rightarrow \infty}g_n = 0.$$
This implies by condition $(1)_n$ that
$$(a_{\mathbb O},a) = (\lim_{n \rightarrow \infty}(a_{\mathbb O} + g_n),a) =
\lim_{n \rightarrow \infty}(a_{\mathbb O} + g_n,a) \in \widehat{{\bI}_P},$$
since $(\widehat{{\bI}_P})_a = \widehat{({\bI}_P)_a}$ is complete.

\vskip.1in

\begin{case}{Construction of $g_0$.}

Set 
$$g_0 = f - a_{\mathbb O}.$$

Then we check
\begin{center} 
\begin{tabular}{ll}
$(0)_0$ & $g_0 = f - a_{\mathbb O} = \sum_{B \neq {\mathbb O}}a_BH^B \in
({\mathcal H})$,\\
$(1)_0$ & $(a_{\mathbb O} + g_0, a) = (f,a) \in \widehat{{\bI}_P}$,
\\ 
$(2)_0$ & the condition $(2)_n$ is void when $n = 0$.\\
\end{tabular} 
\end{center} 
\end{case}

\begin{case}{Construction of $g_n$ assuming that of $g_{n-1}$.}

\vskip.1in

We look at the power series expansion of the form $(\star)$
$$g_{n-1} = \sum_{B \in ({\bZ}_{\geq 0})^N}a_{B,g_{n-1}}H^B.$$
Set
\begin{displaymath}
\left\{\begin{array}{ll}
\nu &= \on{ord}(g_{n-1}) \\
B_o &= \min\left\{B \mid \on{ord}(a_{B,g_{n-1}}H^B) = \nu\right\} \\
B_o &= (b_{o1}, b_{o2}, \ldots, b_{oN}) \\
\partial_{[B_o]} &= \partial_{X^{[B_o]}} =
\partial_{(x_1^{p^{e_1}})^{b_{o1}}(x_2^{p^{e_2}})^{b_{o2}} \cdots
(x_N^{p^{e_N}})^{b_{oN}}}.    \\
\end{array}\right.
\end{displaymath}
Note that $B_o \neq {\mathbb O}$, which follows from condition $(0)_{n-1}$.

We set
$$g_n = \left(1 - H^{B_o}\partial_{[B_o]}\right)g_{n-1}.$$

We check conditions $(0)_n, (1)_n$ and $(2)_n$ in the following.

\vskip.1in

$\boxed{\on{Condition\ }(0)_n}$

\vskip.1in

We compute
$$g_n = \left(1 - H^{B_o}\partial_{[B_o]}\right)g_{n-1} = g_{n-1} -
H^{B_o}\partial_{[B_o]}g_{n-1},$$
where 
\begin{center} 
\begin{tabular}{ll}
&$g_{n-1} \in ({\mathcal H})$ by condition $(0)_{n-1}$, and \\
&$H^{B_o}\partial_{[B_o]}g_{n-1} \in ({\mathcal H})$ since $B_o \neq {\mathbb
O}$. \\
\end{tabular}
\end{center}
Therefore, we conclude
$$g_n \in ({\mathcal H}),$$
checking condition $(0)_n$.

\vskip.1in

$\boxed{\on{Condition\ }(1)_n}$

\vskip.1in

By inductional hypothesis, condition $(1)_{n-1}$ holds, 
i.e., we have the first inclusion
$$(a_{\mathbb O} + g_{n-1}, a) \in \widehat{{\bI}_P}.$$
Since ${\bI}_P$ is ${\fD}$-saturated, so is $\widehat{{\bI}_P}$ (cf.
compatibility of completion with ${\fD}$-saturation, Proposition 2.4.2.1 (2)
in Part I). 
Therefore, the first inclusion implies the second inclusion
$$(\partial_{[B_o]}(a_{\mathbb O} + g_{n-1}),a - |[B_o]|) \in
\widehat{{\bI}_P}.$$
The second inclusion combined with the third inclusion below
$$(H^{B_o}, |[B_o]|) \in \widehat{{\bI}_P}$$
implies the fourth inclusion
$$(H^{B_o}\partial_{[B_o]}(a_{\mathbb O} 
+ g_{n-1}),a) \in \widehat{{\bI}_P}.$$
Therefore, we conclude
\begin{displaymath}
\begin{array}{ll}
(a_{\mathbb O} + g_n,a) &= (a_{\mathbb O} + \left(1 -
H^{B_o}\partial_{[B_o]}\right)g_{n-1},a) \\
&= (\left(1 - H^{B_o}\partial_{[B_o]}\right)(a_{\mathbb O} + g_{n-1}),a)
\hskip.1in (\on{since\ }H^{B_o}\partial_{[B_o]}a_{\mathbb O} = 0 \on{\ as\
}B_o
\neq 0)\\ &= ((a_{\mathbb O} + g_{n-1}) - H^{B_o}
\partial_{[B_o]}(a_{\mathbb O} +
g_{n-1}),a) \in \widehat{{\bI}_P} \\
&\hfill (\on{by\ condition\ }(1)_{n-1} \on{\ and\
by\ the\ fourth\ inclusion\ above}).\\
\end{array}
\end{displaymath}
That is to say, we have
$$(a_{\mathbb O} + g_n,a) \in \widehat{{\bI}_P},$$
checking condition $(1)_n$.

\vskip.1in

$\boxed{\on{Condition\ }(2)_n}$

\vskip.1in

Observe that, for any $B \in ({\bZ}_{\geq 0})^N$ with
$\on{ord}\left(a_{B,g_{n-1}}H^B\right) = \nu$, we have

\begin{eqnarray*}
\left(1 - H^{B_o}\partial_{[B_o]}\right)\left(a_{B,g_{n-1}}H^B\right) &=&
a_{B,g_{n-1}}H^B - H^{B_o}\left(\binom{B}{B_o}a_{B,g_{n-1}}H^{B -
B_o} + s_B\right)\\
&=& \left(1 - \binom{B}{B_o}\right)a_{B,g_{n-1}}H^B + r_B \\
\end{eqnarray*}
where $s_B$ and $r_B$ are elements in $\widehat{R_P}$ with $\on{ord}(s_B) >
\nu - |[B_o]|$ and $\on{ord}(r_B) > \nu$, respectively.

Therefore, we compute
\begin{eqnarray*}
\lefteqn{
g_n = \left(1 - H^{B_o}\partial_{[B_o]}\right)g_{n-1} 
= \left(1 - H^{B_o}\partial_{[B_o]}\right)\left(\sum
a_{B,g_{n-1}}H^B\right) }\\
&=& \left(1 - H^{B_o}\partial_{[B_o]}\right)\left(\sum_{B \mid
\on{ord}\left(a_{B,g_{n-1}}H^B\right) = \nu}
a_{B,g_{n-1}}H^B + \sum_{B \mid \on{ord}\left(a_{B,g_{n-1}}H^B\right) > \nu}
a_{B,g_{n-1}}H^B\right) \\
&=& \sum_{B \mid \on{ord}\left(a_{B,g_{n-1}}H^B\right) = \nu} \left(1 -
\binom{B}{B_o}\right)a_{B,g_{n-1}}H^B + r 
\\
&=& \sum_{B \mid \on{ord}\left(a_{B,g_{n-1}}H^B\right) = \nu, 
\hskip.02in B> B_o} \left(1 - \binom{B}{B_o}\right)a_{B,g_{n-1}}H^B + r \\
\end{eqnarray*}
where $r$ is an element in $\widehat{R_P}$ with $\on{ord}(r) > \nu$.

{}From this computation it immediately follows that
$$\on{ord}(g_n) \geq \nu = \on{ord}(g_{n-1})$$
and that, if $\on{ord}(g_n) = \nu = \on{ord}(g_{n-1})$, then

\begin{displaymath}
\begin{array}{rcl}
\on{the\ 2nd\ factor\ in\ }\eta(g_n) &=& \min\left\{B \mid 
\on{ord}(a_{B,g_n}H^B)= \nu\right\} \\
&>& B_o \\ 
&=& \on{the\ 2nd\ factor\ in\ }\eta(g_{n-1}). \\
\end{array}
\end{displaymath}

Thus we conclude
$$\eta(g_{n-1}) < \eta(g_n),$$
checking condition $(2)_n$.
\end{case}

This completes the inductive construction of the sequence 
$\left\{g_n\right\}_{n \in
{\bZ}_{\geq 0}} \subset \widehat{R_P}$ satisfying conditions $(0)_n,
(1)_n$ and $(2)_n$.

This completes the argument in Step 1, showing $(a_{\mathbb O},a) \in
\widehat{{\bI}_P}$.
\end{step}

\vskip.1in

\begin{step}{We show the statement in the general case, i.e., we show
$(a_B,a - |[B]|) \in \widehat{{\bI}_P}$ for any $B \in
({\bZ}_{\geq 0})^N$.}

We will construct a sequence $\left\{g_n\right\}_{n 
\in {\bZ}_{\geq 0}} \subset \widehat{R_P}$ inductively, 
satisfying the following conditions:
\begin{center} 
\begin{tabular}{ll}
$(0)_n$ & $(g_n, a) \in \widehat{{\bI}_P}$,
\\ 
$(1)_n$ & $(a_{B,f - g_n},a - |[B]|) \in \widehat{{\bI}_P}$ for any $B \in
({\bZ}_{\geq 0})^N$,\\
$(2)_n$ & $\eta(g_{n-1}) < \eta(g_n)$.\\
\end{tabular} 
\end{center} 

The construction of such a sequence is sufficient to prove the statement
$$(a_B,a - |[B]|) \in \widehat{{\bI}_P} \hskip.1in \on{for\ any} \hskip.1in
B \in ({\bZ}_{\geq 0})^N.$$

In fact, since there are only finitely many $B$'s with $\on{ord}(H^B) \leq
\nu$ for a fixed $\nu \in {\bZ}_{\geq 0}$, and since
$$\left\{B \mid \on{ord}(H^B) \leq \nu\right\} \supset 
\left\{B \mid \on{ord}(a_{B,g}H^B) =
\on{ord}(g)\right\}$$
for any $g \in \widehat{R_P}$ with $\on{ord}(g) = \nu$, we conclude by
condition $(2)_n$ that
$$\lim_{n \rightarrow \infty}\on{ord}(g_n) = \infty \hskip.1in 
\on{and\ hence} \hskip.1in \lim_{n \rightarrow \infty}g_n = 0.$$
This implies by condition $(1)_n$ that, for any $B \in ({\bZ}_{\geq 0})^N$,
$$(a_B,a - |[B]|) = (a_{B,f},a - |[B]|) = (\lim_{n \rightarrow
\infty}(a_{B,f - g_n}),a - |[B]|) = \lim_{n \rightarrow \infty}(a_{B,f -
g_n},a)
\in
\widehat{{\bI}_P},$$ since $(\widehat{{\bI}_P})_a = \widehat{({\bI}_P)_a}$
is complete.

\vskip.1in

\setcounter{case}0\begin{case}{Construction of $g_0$.}

Set 
$$g_0 = f.$$

Then we check
\begin{center} 
\begin{tabular}{ll}
$(0)_0$ & $(g_0,a) = (f,a) \in {\bI}_P \subset \widehat{{\bI}_P}$,\\
$(1)_0$ & $(a_{B,f - g_0}, a - |[B]|) = (0,a - |[B]|) \in \widehat{{\bI}_P}$
for any $B \in ({\bZ}_{\geq 0})^N$,
\\ 
$(2)_0$ & the condition $(2)_n$ is void when $n = 0$.\\
\end{tabular} 
\end{center} 
\end{case}

\newpage

\begin{case}{Construction of $g_n$ assuming that of $g_{n-1}$.}

\vskip.1in

We look at the power series expansion of the form $(\star)$
$$g_{n-1} = \sum_{B \in ({\bZ}_{\geq 0})^N}a_{B,g_{n-1}}H^B.$$
Set
\begin{displaymath}
\left\{\begin{array}{ll}
\nu &= \on{ord}(g_{n-1}) \\
B_o &= \min\left\{B \mid \on{ord}(a_{B,g_{n-1}}H^B) = \nu\right\} \\
B_o &= (b_{o1}, b_{o2}, \ldots, b_{oN}) \\
\partial_{[B_o]} &= \partial_{X^{[B_o]}} = 
\partial_{(x_1^{p^{e_1}})^{b_{o1}}(x_2^{p^{e_2}})^{b_{o2}} \cdots 
(x_N^{p^{e_N}})^{b_{oN}}}.    \\
\end{array}\right.
\end{displaymath}

We set
$$g_n = \left(1 - H^{B_o}a_{{\mathbb O},*}\partial_{[B_o]}\right)g_{n-1} = 
g_{n-1} - a_{{\mathbb O},\partial_{[B_o]}g_{n-1}}H^{B_o},$$
where $a_{{\mathbb O},*}$ denotes the operator such that $a_{{\mathbb O},*}g = 
a_{{\mathbb O},g}$ for any $g \in \widehat{R_P}$.

We check conditions $(0)_n, (1)_n$ and $(2)_n$ in the following.

\vskip.1in

$\boxed{\on{Condition\ }(0)_n}$

\vskip.1in

By inductional hypothesis, condition $(0)_{n-1}$ holds, 

i.e., we have the first inclusion
$$(g_{n-1},a) \in \widehat{{\bI}_P}.$$
Since ${\bI}_P$ is ${\fD}$-saturated, so is $\widehat{{\bI}_P}$ (cf. 
compatibility of completion with ${\fD}$-saturation, Proposition 2.4.2.1 (2) 
in Part I). 
Therefore, the first inclusion implies the second inclusion
$$(\partial_{[B_o]}g_{n-1},a - |[B_o]|) \in \widehat{{\bI}_P}.$$
By Step 1 the second inclusion implies the third
$$(a_{{\mathbb O},\partial_{[B_o]}g_{n-1}},a - |[B_o]|) \in 
\widehat{{\bI}_P}.$$
The third inclusion combined with the fourth inclusion below
$$(H^{B_o},|[B_o]|) \in \widehat{{\bI}_P}$$
implies the fifth inclusion
$$(a_{{\mathbb O},\partial_{[B_o]}g_{n-1}}H^{B_o},a) \in \widehat{{\bI}_P}.$$
Therefore, we conclude
$$(g_n,a) = (g_{n-1} - a_{{\mathbb O},\partial_{[B_o]}g_{n-1}}H^{B_o},a) \in 
\widehat{{\bI}_P},$$
checking condition $(0)_n$.

\vskip.1in

$\boxed{\on{Condition\ }(1)_n}$

\vskip.1in

When $B \neq B_o$, we have
$$a_{B,f - g_n} = a_{B,f - g_{n-1} - a_{{\mathbb 
O},\partial_{[B_o]}g_{n-1}}H^{B_o}} = a_{B,f - g_{n-1}}.$$
Therefore, by condition $(1)_{n-1}$, we conclude
$$(a_{B,f - g_n},a - |[B]|) = (a_{B,f - g_{n-1}},a - |[B]|) \in 
\widehat{{\bI}_P}.$$
When $B = B_o$, we have
$$a_{B_o,f - g_n} = a_{B_o,f - g_{n-1} - a_{{\mathbb 
O},\partial_{[B_o]}g_{n-1}}H^{B_o}} = a_{B_o,f - g_{n-1}} - a_{{\mathbb
O},\partial_{[B_o]}g_{n-1}}.$$
Therefore, by condition $(1)_{n-1}$ and Step 1, we conclude
$$(a_{B_o,f - g_n},a - |[B_o]|) = (a_{B_o,f - g_{n-1}} - a_{{\mathbb 
O},\partial_{[B_o]}g_{n-1}},a - |[B]|) \in \widehat{{\bI}_P}.$$
This checks condition $(1)_n$.

\vskip.1in

$\boxed{\on{Condition\ }(2)_n}$

\vskip.1in

Observe that, for any $B \in ({\bZ}_{\geq 0})^N$ with 
$\on{ord}\left(a_{B,g_{n-1}}H^B\right) = \nu$, we have

\begin{eqnarray*}
\left(1 - H^{B_o}a_{{\mathbb 
O},*}\partial_{[B_o]}\right)\left(a_{B,g_{n-1}}H^B\right) &=& 
a_{B,g_{n-1}}H^B -
H^{B_o}a_{{\mathbb O},*}\left(\binom{B}{B_o}a_{B,g_{n-1}}H^{B - B_o} + 
s_B\right)\\
&=& \left(1 - \delta_{BB_o}\right)a_{B,g_{n-1}}H^B + r_B' \\
\end{eqnarray*}
where $s_B$ and $r_B'$ are elements in $\widehat{R_P}$ with $\on{ord}(s_B) > 
\nu - |[B_o]|$ and $\on{ord}(r_B') > \nu$, respectively, and where
$\delta_{BB_o}$ denotes the Kronecker delta.

Therefore, we compute
\begin{eqnarray*}
\lefteqn{
g_n = \left(1 - H^{B_o}a_{{\mathbb O},*}\partial_{[B_o]}\right)g_{n-1} 
= \left(1 - H^{B_o}a_{{\mathbb O},*}\partial_{[B_o]}\right)\left(\sum 
a_{B,g_{n-1}}H^B\right)} \\
&=& \left(1 - H^{B_o}a_{{\mathbb O},*}
\partial_{[B_o]}\right)\left(\sum_{B \mid 
\on{ord}\left(a_{B,g_{n-1}}H^B\right) = \nu}
a_{B,g_{n-1}}H^B + \sum_{B \mid \on{ord}\left(a_{B,g_{n-1}}H^B\right) > \nu}
a_{B,g_{n-1}}H^B\right) \\
&=& \sum_{B \mid \on{ord}\left(a_{B,g_{n-1}}H^B\right) = \nu} \left(1 - 
\delta_{BB_o}\right)a_{B,g_{n-1}}H^B + r' 
= \sum_{B \mid \on{ord}\left(a_{B,g_{n-1}}H^B\right) = \nu, \hskip.02in B 
> B_o} a_{B,g_{n-1}}H^B + r' \\
\end{eqnarray*}
where $r'$ is an element in $\widehat{R_P}$ with $\on{ord}(r') > \nu$.

{}From this computation it immediately follows that
$$\on{ord}(g_n) \geq \nu = \on{ord}(g_{n-1})$$
and that, if $\on{ord}(g_n) = \nu = \on{ord}(g_{n-1})$, then

\begin{displaymath}
\begin{array}{rcl}
\on{the\ 2nd\ factor\ in\ }\eta(g_n) &=& \min\left\{B \mid 
\on{ord}(a_{B,g_n}H^B) = \nu\right\} \\
&>& B_o \\ 
&=& \on{the\ 2nd\ factor\ in\ }\eta(g_{n-1}). \\
\end{array}
\end{displaymath}

Thus we conclude
$$\eta(g_{n-1}) < \eta(g_n),$$
checking condition $(2)_n$.
\end{case}

This completes the inductive construction of the sequence 
$\left\{g_n\right\}_{n \in 
{\bZ}_{\geq 0}} \subset \widehat{R_P}$ satisfying conditions $(0)_n,
(1)_n$ and $(2)_n$.

This completes the argument in Step 2, showing $(a_B,a - |[B]|) \in 
\widehat{{\bI}_P}$.
\end{step}

\vskip.1in

This completes the alternative proof of Lemma \ref{2.2.2.1}.
\end{proof}

\begin{rem}\label{2.2.2.2}

\begin{enumerate}

\item The basic strategy of the above proof can be seen in a more 
transparent way, if we consider the following special case:  
Suppose that we can take the associated regular system of 
parameters $(x_1, \ldots, x_d)$ in such a way that $h_l = x_l$ for
$l = 1, \ldots, N$. 

Given $(f,a) \in {\bI}_P$, with $f = \sum_Ba_BH^B$ being the power series 
expansion of the form $(\star)$, we proceed as follows.

\begin{itemize}

\item[Step 1.] We compute
$$\prod_{0 < |[B]| < a}\left(1 - H^B\partial_{H^B}\right)f = a_{{\mathbb O}} - 
\sum_{|[B]| \geq a}c_BH^B$$
for some $c_B \in \widehat{R_P}$.  (Note that, since the operators $\left(1 
- H^B\partial_{H^B}\right)$ do not commute, the product symbol
$\prod_{0 < |[B]| < a}$ is understood to align the factors from right to left 
according to the lexicographical order among $(|[B]|,B)$'s.)

Since ${\bI}_P$ is ${\fD}$-saturated, we see $(\prod_{0 < |[B]| < a}\left(1 
- H^B\partial_{H^B}\right)f,a) \in {\bI}_P$, while we have obviously
$(\sum_{|[B]| \geq a}c_BH^B,a) \in \widehat{{\bI}_P}$.  Therefore, we 
conclude
$$(a_{\mathbb O},a) \in \widehat{{\bI}_P}.$$

\item[Step 2.] For $B \in ({\bZ}_{\geq 0})^N$, 
set $g = \partial_{H^B}f$.  Let $g = 
\sum_{\beta}a_{\beta,g}H^{\beta}$
be the power series expansion of the form $(\star)$.  Observe $a_{{\mathbb 
O},g} = a_B$.  Since ${\bI}_P$ is ${\fD}$-saturated, we see $(g,a
- |[B]|)
\in {\bI}_P$.  By the previous step, we conclude
$$(a_B,a - |[B]|) = (a_{{\mathbb O},g},a - |[B]|) \in \widehat{{\bI}_P}.$$

\end{itemize}

In the general case,  since the set ${\mathcal H}$ does not 
coincide with a part of the associated regular system of parameters, 
we can not follow the steps of the special case above {\it literally}.  
However, by substituting $\partial_{X^{[B]}}$ for $\partial_{H^B}$ 
and by filling the gap of the substitution through the process 
of taking the limits, we can try to follow the steps of the 
special case {\it in spirit}.  That is the basic strategy 
of the above proof.

\item In Chapter 3, we derive Lemma 4.1.4.1 (Coefficient Lemma) 
in Part I as a corollary to the formal coefficient lemma above. 
\end{enumerate}

\end{rem}

\end{subsection}

\end{section} 
\end{chapter} 

\begin{chapter}{Invariant $\widetilde{\mu}$}

The purpose of this chapter is to study the basic properties of the
invariant $\widetilde{\mu}$.  As the unit for the strand of invariants 
in our algorithm is a triplet of numbers
$(\sigma,\widetilde{\mu},s)$ (or a quadruplet $(\sigma,\widetilde{\mu},
\widetilde{\nu},s)$ (cf. \ref{0.3.1})), we also study the behavior of 
the pair $(\sigma,\widetilde{\mu})$ endowed with the lexicographical 
order.  The discussion of the invariant $\widetilde{\mu}$ in this 
chapter is restricted to and concentrated on the case where there 
are no exceptional divisors involved, and hence can only be applied 
directly to the process {\it at year 0} of our algorithm.  We will 
postpone the general discussion, involving the exceptional divisors 
and hence applicable to the process {\it after year 0} of our 
algorithm, to Part III or Part IV (cf. \ref{0.2.3}).

The setting for this chapter is identical to that of Chapter 1.

Namely, $R$ represents the coordinate ring of an affine open subset
$\on{Spec}\hskip.02in R$ of a nonsingular variety $W$ of $\dim W = d$ over
an algebraically
closed  field $k$ of positive characteristic ${\on{char}}(k) = p$ or of
characteristic zero ${\on{char}}(k) = 0$, where in the latter case we
formally set 
$p = \infty$ (cf. 0.2.3.2.1 and Definition 3.1.1.1 (2) in Part I).

Let ${\bI}$ be an idealistic filtration over $R$.  We assume that ${\bI}$ is
${\fD}$-saturated.  We remark
that then, by compatibility of localization with ${\fD}$-saturation (cf.
Proposition 2.4.2.1 (2) in Part I), the localization ${\bI}_P$ is 
also ${\fD}$-saturated for any closed point $P \in {\on{Spec}}\ R$.

\begin{section}{Definition of $\widetilde{\mu}$.}\label{3.1}

\begin{subsection}{Definition of $\widetilde{\mu}$ as 
${\mu}_{\mathcal H}$.}\label{3.1.1} We fix a closed point 
$P \in \on{Spec}\hskip.02in R \subset W$.  
Take a leading generator system
${\bH} =
\left\{(h_l,p^{e_l})\right\}_{l = 1, \ldots, N}$ 
with associated nonnegative integers
$0 \leq e_1 \leq \cdots \leq e_N$ for the ${\fD}$-saturated idealistic
filtration
${\bI}_P$.  Let
${\mathcal H} =
\left\{h_l\right\}_{l = 1,
\ldots, N}$ be the set of its elements in $R_P$, and $({\mathcal H}) 
\subset R_P$ the ideal generated by ${\mathcal H}$.

\begin{defn}\label{3.1.1.1} First we recall a few definitions given 
in \S 3.2 in Part I.  For $f \in R_P$ (or more generally for $f \in
\widehat{R_P}$), we define its multiplicity (or order) modulo
$({\mathcal H})$, denoted by
$\on{ord}_{\mathcal H}(f)$, to be
$$\begin{array}{lcl}
\on{ord}_{\mathcal H}(f) &=& \sup\left\{n \in {\bZ}_{\geq 0} \mid f \in
{\fm}_P^n + ({\mathcal H})\right\} \\
&(\text{or}& \sup\left\{n \in {\bZ}_{\geq 0} \mid f \in \widehat{{\fm}_P}^n
+ ({\mathcal H})\right\}). \\
\end{array}$$
Note that we set $\on{ord}_{\mathcal H}(0) = \infty$ by definition.  
We also define 
$${\mu}_{\mathcal H}({\bI}_P) := \inf\left\{\mu_{\mathcal H}(f,a) =
\frac{\on{ord}_{\mathcal H}(f)}{a} \mid (f,a) \in {\bI}_P, a >
0\right\}.$$
(We remark that ${\mu}_{\mathcal H}(\widehat{{\bI}_P})$ is defined 
in a similar manner.)

Finally the invariant $\widetilde{\mu}$ at $P$, which we denote by
$\widetilde{\mu}(P)$, is defined by the formula
$$\widetilde{\mu}(P) = {\mu}_{\mathcal H}({\bI}_P).$$
\end{defn}

In order to justify the definition, we should show that ${\mu}_{\mathcal
H}({\bI}_P)$ is independent of the choice of ${\mathcal H}$, 
i.e.,
independent of the choice of a leading generator system ${\bH}$ for
${\bI}_P$.  We will show this
independence in the next subsection.

\begin{rem}\label{3.1.1.2} (1) The usual order is multiplicative, 
i.e., we have an equality
$$\on{ord}(fg) = \on{ord}(f) + \on{ord}(g) \hskip.1in \forall f, g \in
R_P.$$
The order modulo $({\mathcal H})$ is also multiplicative if $e_1 =
\cdots = e_N = 0$.  However, in general, we can only expect that 
the order modulo $({\mathcal H})$ is only weakly multiplicative, 
i.e., we have only an inequality
$$\on{ord}(fg) \geq \on{ord}(f) + \on{ord}(g) \hskip.1in \forall f, g \in
R_P.$$
In fact, if $e_l > 0$ for some $l = 1, \ldots, N$, then it is 
easy to see (cf. Remark \ref{3.2.1.2} (1)) that we indeed have 
a strict inequality for some $f, g \in R_P$, i.e.,
$$\on{ord}(fg) > \on{ord}(f) + \on{ord}(g) \hskip.1in \text{\ for\ some\ }
f, g \in R_P.$$

(2) Assume further that the idealistic filtration ${\bI}$ is of r.f.g.~type
(cf. Definition 2.1.1.1 (4) and \S 2.3 in Part
I).  Then the invariant $\widetilde{\mu}$ takes the rational values with
some bounded denominator $\delta$ (independent of $P$).

In fact, take a
finite set of generators $T$ for ${\bI} = G_R(T)$ of the form
$$T = \left\{(f_{\lambda},a_{\lambda})
\right\}_{\lambda \in \Lambda} \subset R \times
{\bQ}_{> 0}, \# \Lambda < \infty \hskip.03in \on{with} \hskip.03in
a_{\lambda} =
\frac{p_{\lambda}}{q_{\lambda}} \hskip.03in \on{where} \hskip.03in
p_{\lambda}, q_{\lambda} \in {\bZ}_{> 0}.$$

Set $\delta = \prod_{\lambda \in \Lambda} p_{\lambda}$.

Then

\begin{eqnarray*}
\widetilde{\mu}(P) &=& {\mu}_{\mathcal H}({\bI}_P) =
\inf\left\{\mu_{\mathcal H}(f,a) = \frac{\on{ord}_{\mathcal H}(f)}{a} \mid
(f,a) \in {\bI}_P, a >
0\right\}\\
&=& \min\left\{\mu_{\mathcal H}(f_{\lambda},a_{\lambda}) =
\frac{\on{ord}_{\mathcal H}(f_{\lambda})}{a_{\lambda}} =
\frac{\on{ord}_{\mathcal
H}(f_{\lambda}) \cdot q_{\lambda}}{p_{\lambda}}\right\}\\
&& (\text{cf.\ Lemma\ 2.2.1.2\ (1)\ in\ Part\ I\ 
and\ Remark\ \ref{3.1.1.2}\ (1)\ above})\\
&\in& \frac{1}{\delta}{\bZ}_{\geq 0} \cup \{\infty\}.\\
\end{eqnarray*}
\end{rem}

\end{subsection}

\begin{subsection}{Invariant ${\mu}_{\mathcal H}$ is independent of
${\mathcal H}$.}\label{3.1.2} We show that ${\mu}_{\mathcal H}({\bI}_P)$
is independent of the choice of ${\mathcal H}$.

\begin{prop}\label{3.1.2.1} Let the setting be as described in \ref{3.1.1}.

Then ${\mu}_{\mathcal H}({\bI}_P)$ is independent of the choice of
${\mathcal H}$, 
i.e., independent of the choice of a
leading generator system ${\bH}$ for ${\bI}_P$.
\end{prop}

\begin{proof} Suppose
$$\mu_P({\mathbb I}) = \inf\left\{\mu_P(f,a) = \frac{\on{ord}_P(f)}{a} \mid
(f,a) \in {\bI}_P, a > 0\right\} < 1.$$
Then, since
${\bI}_P$ is ${\fD}$-saturated, we have ${\bI}_P = 
R_P \times {\mathbb R}$ (cf. Lemma 1.1.2.1 $\mathbf{Case: 
P \not\in \on{Supp}(\bI)}$).  We conclude that the set of 
elements ${\mathcal H}$ in any leading generator system
${\bH}$ for ${\bI}_P$ is a regular system of parameters
$\{x_1, \ldots , x_d\}$ for $R_P$, where $d = \dim W$.  
Accordingly, we have $${\mu}_{\mathcal H}({\mathbb I}_P) = 0,$$
independent of the choice of
${\mathcal H}$.

Therefore, in the following, we may assume $1 \leq \mu_P({\bI})$ and hence
that $1 \leq \mu_P({\bI}) \leq {\mu}_{\mathcal H}({\bI}_P)$ 
for any choice of ${\mathcal H}$.

\vskip.1in

\begin{case}{${\mu}_{\mathcal H}({\bI}_P) = 1$ for any choice of ${\mathcal
H}$.}

In this case, ${\mu}_{\mathcal H}({\bI}_P) = 1$ is obviously independent of
the choice of ${\mathcal H}$ by the case assumption.
\end{case}

\vskip.1in

\begin{case}{${\mu}_{\mathcal H}({\bI}_P) > 1$ for some choice of ${\mathcal
H}$.}

In this case, fixing the set of elements ${\mathcal H}$ of a leading
generator system ${\bH}$ for ${\bI}_P$ with ${\mu}_{\mathcal
H}({\bI}_P) > 1$, we show
$$(*) \hskip.1in {\mu}_{\mathcal H'}({\bI}_P) \geq {\mu}_{\mathcal
H}({\bI}_P) \hskip.1in (> 1)$$
where ${\mathcal H'}$ is the set of elements of another leading
generator system ${\bH'}$ for ${\bI}_P$.

This is actually sufficient to show the required independence, since by
switching the roles of ${\mathcal H}$ and
${\mathcal H'}$, we conclude ${\mu}_{\mathcal H}({\bI}_P) \geq
{\mu}_{\mathcal
H'}({\bI}_P)$ and hence ${\mu}_{\mathcal H}({\bI}_P) = {\mu}_{\mathcal
H'}({\bI}_P)$.

\vskip.1in

First we make the following two easy observations:

\begin{enumerate}

\item Let ${\mathcal H''} = \left\{h_l''\right\}_{l = 1, \ldots, N}$ 
be another set of elements in $R_P$ obtained from ${\mathcal H'}$ 
by a linear transformation, i.e., for each
$e
\in {\bZ}_{\geq 0}$ we have
$$\left[h_l''^{p^{e - e_l}} \mid e_l \leq e\right] = \left[h_l'^{p^{e -
e_l}} \mid e_l \leq e\right]g_e \hskip.05in \text{for\ some\ } g_e
\in \on{GL}\left(\#\left\{e_l \mid e_l \leq e\right\},k\right).$$
Then ${\mathcal H''}$ is the set of elements of a leading generator system
${\bH''}=\left\{(h_l'',p^{e_l})\right\}_{l=1}^N$ for ${\bI}_P$, and we
have
$${\mu}_{\mathcal H'}({\bI}_P) = {\mu}_{\mathcal H''}({\bI}_P).$$

Going back to our situation, we see that there is ${\mathcal H''}$, 
obtained from ${\mathcal H'}$ by a linear
transformation, such that ${\mathcal H''}$ and ${\mathcal H}$ share the same
leading terms.

Therefore, in order to show the inequality $(*)$, by replacing ${\mathcal
H'}$ with ${\mathcal H''}$ we may assume that ${\mathcal H}$ and
${\mathcal H'}$ share the same leading terms, 
i.e.,
$$h_l \equiv h_l' \bmod {\fm}_P^{p^{e_l} + 1}
\hskip.1in \on{for} \hskip.1in l = 1, \ldots, N.$$

\item Assume that ${\mathcal H}$ and ${\mathcal H'}$ share the same leading
terms.  Then we have a sequence of the sets of elements of leading
generator systems for
${\bI}_P$ 
$${\mathcal H} = {\mathcal H}_0, {\mathcal H}_1, \ldots, {\mathcal H}_N =
{\mathcal H'}$$
where the adjacent sets share all but one elements in common.  We have only
to show
$${\mu}_{{\mathcal H}_l}({\bI}_P) \geq {\mu}_{{\mathcal H}_{l-1}}({\bI}_P)
\hskip.1in \on{for} \hskip.1in l = 1, \ldots, N$$
in order to verify the inequality $(*)$.
\end{enumerate}

\vskip.1in

According to the observations above, therefore, we have only to show 
the inequality $(*)$ under the following extra assumptions:

\begin{enumerate}

\item ${\mathcal H}$ and ${\mathcal H'}$ share the same leading terms, 
i.e.,
$$h_l \equiv h_l' \bmod  {\fm}_P^{p^{e_l} + 1}
\hskip.1in \on{for} \hskip.1in l = 1, \ldots, N.$$

\item ${\mathcal H}$ and ${\mathcal H'}$ share all but one element in
common, 
i.e.,
$$h_l = h_l' \hskip.1in \on{for} \hskip.1in l = 1, \ldots, N \hskip.1in
\on{except} \hskip.1in l = l_o.$$
\end{enumerate}

\vskip.1in

In order to ease the notation, we set
$$h = h_{l_o}, h' = h_{l_o}', {\mathcal G} = {\mathcal H} \setminus
\{h_{l_o}\} = {\mathcal H'} \setminus
\{h_{l_o}'\}.$$ 
Let $\nu$ be any positive number such that $1 < \nu < {\mu}_{\mathcal
H}({\bI}_P)$.

Since $(h,p^{e_{l_o}}), (h',p^{e_{l_o}}) \in {\bI}_P$, we have $(h -
h',p^{e_{l_o}}) \in {\bI}_P$.  Therefore, by definition of ${\mu}_{\mathcal
H}({\bI}_P)$ and by the inequality $1 < \nu < {\mu}_{\mathcal H}({\bI}_P)$,
we have
$$h - h' \in {\fm}_P^{\lceil \nu p^{e_{l_o}} \rceil} + ({\mathcal H})
\hskip.1in \on{
i.e.,} \hskip.1in h - h' = f_1 + f_2 \hskip.05in \on{with}
\hskip.05in f_1 \in {\fm}_P^{\lceil \nu p^{e_{l_o}} \rceil}, f_2 \in
({\mathcal H}).$$
On the other hand, by extra assumption (1), we have
$$h - h' \in {\fm}_P^{p^{e_{l_o}} + 1}.$$
Observing $\lceil \nu p^{e_{l_o}} \rceil \geq p^{e_{l_o}} + 1$ (recall $\nu
> 1$), we thus conclude that
$$f_2 = (h - h') - f_1 \in ({\mathcal H}) \cap {\fm}_P^{p^{e_{l_o}} + 1}$$
and hence that
$$
\begin{array}{ll}
h - h' = f_1 + f_2 &\in \hskip.05in {\fm}_P^{\lceil \nu p^{e_{l_o}} 
\rceil} + ({\mathcal H}) \cap {\fm}_P^{p^{e_{l_o}} + 1} \\
&\subset \hskip.05in {\fm}_P^{\lceil \nu p^{e_{l_o}} \rceil} + h{\fm}_P +
({\mathcal G}) \cap {\fm}_P^{p^{e_{l_o}} + 1},
\\
\end{array}$$
where the second inclusion follows from Lemma 4.1.2.3 in Part I.

That is to say, we have
$$h - h' = g_1 + hr + g_2 \hskip.1in \on{with} \hskip.1in g_1 \in
{\fm}_P^{\lceil \nu p^{e_{l_o}} \rceil}, r \in {\fm}_P, g_2 \in ({\mathcal
G}) \cap
{\fm}_P^{p^{e_{l_o}} + 1}.$$
Therefore, we have
$$(1 - r)h = g_1 + h' + g_2.$$
Since $u = 1 - r$ is a unit in $R_P$, we conclude
$$h = u^{-1}g_1 + u^{-1}h' + u^{-1}g_2 \in {\fm}_P^{\lceil \nu p^{e_{l_o}}
\rceil} + ({\mathcal H'}).$$

Given an element $(f,a) \in {\bI}_P \hskip.1in (a > 0)$, we hence have

\begin{eqnarray*}
f &\in& \sum_B ({\bI}_P)_{a - |[B]|}'H^B \hskip.2in 
(\on{by\ Coefficient\ Lemma\ in\ Part\ I,\ where\ }({\bI}_P)_t' =
({\bI}_P)_t \cap {\fm}_P^{\lceil \nu t \rceil}) \\
&=& \sum_{b = b_{l_o}, C = (b_1, \ldots, b_{l_o - 1}, 0, b_{l_o + 1},
\ldots, b_N)}({\bI}_P)_{a - |[C]| - bp^{e_{l_o}}}'h^bH^C \hskip.1in 
\subset \hskip.1in \sum_b({\bI}_P)_{a - bp^{e_{l_o}}}'h^b + ({\mathcal G}) \\
&\subset& \sum_b({\bI}_P)_{a - bp^{e_{l_o}}}'{\fm}_P^{b \lceil \nu
p^{e_{l_o}} \rceil} + ({\mathcal H'}) \hskip.2in (\on{since\ }h 
\in {\fm}_P^{\lceil \nu p^{e_{l_o}} \rceil} + ({\mathcal
H'})). \\
\end{eqnarray*}

Therefore, we compute

\begin{eqnarray*}
\on{ord}_{\mathcal H'}(f) &\geq& \min_b\left\{\on{ord}_P\left(({\bI}_P)_{a -
bp^{e_{l_o}}}'{\fm}_P^{b \lceil \nu p^{e_{l_o}} \rceil}\right)\right\} \\
&\geq& \min_b\left\{\lceil \nu (a - bp^{e_{l_o}}) \rceil + b \lceil \nu
p^{e_{l_o}} \rceil\right\} \geq \nu a.\\
\end{eqnarray*}

This implies
$${\mu}_{\mathcal H'}(f,a) = \frac{\on{ord}_{\mathcal H'}(f)}{a} \geq \nu.$$
Since this inequality holds for any positive number with $1 < \nu <
{\mu}_{\mathcal H}({\bI}_P)$, we conclude
$${\mu}_{\mathcal H'}(f,a) \geq {\mu}_{\mathcal H}({\bI}_P).$$
Since $(f,a) \in {\bI}_P \hskip.1in (a > 0)$ is arbitrary, we finally
conclude
$${\mu}_{\mathcal H'}({\bI}_P) \geq {\mu}_{\mathcal H}({\bI}_P).$$

This completes the proof of the inequality $(*)$, and hence the proof of
Proposition \ref{3.1.2.1}.
\end{case}

\end{proof}

\end{subsection}

\end{section}

\begin{section}{Interpretation of $\widetilde{\mu}$ in terms of the power
series expansion.}\label{3.2}

The purpose of this section is to give an interpretation of 
the invariant $\widetilde{\mu} = 
{\mu}_{\mathcal H}$ in terms of the power series expansion of the form
$(\star)$ discussed in Chapter 2.

\begin{subsection}{The order $\on{ord}_{\mathcal H}(f)$ of $f$ modulo
$({\mathcal H})$ is equal to the order $\on{ord}(a_{\mathbb O})$ of the
constant
term of the power series expansion for $f$.}\label{3.2.1}

\begin{lem}\label{3.2.1.1} Let the setting be as described in \ref{3.1.1}.

Then we have
$$\on{ord}_{\mathcal H}(f) = \on{ord}(a_{\mathbb O}),$$
where $a_{\mathbb O}$ is the ``constant'' term of the power series expansion $f
= \sum a_BH^B$ of the form $(\star)$ as described in Lemma \ref{2.1.2.1}.
\end{lem}

\begin{proof} Since $f \equiv a_{\mathbb O} \bmod 
({\mathcal H})$, we obviously have
$$\on{ord}_{\mathcal H}(f) = \on{ord}_{\mathcal H}(a_{\mathbb O}) \geq
\on{ord}(a_{\mathbb O}).$$
Suppose $\on{ord}_{\mathcal H}(f) > \on{ord}(a_{\mathbb O}) = r$.  Then by
definition we can write
$$f = f_1 + f_2 \hskip.1in \on{with} \hskip.1in f_1 \in \widehat{{\fm}_P}^{r
+ 1}, f_2 \in ({\mathcal H}).$$
Therefore, we have
$$f_1 = f - f_2 = \sum a_BH^B - f_2.$$
Since $f_2 \in ({\mathcal H})$, we conclude by the uniqueness of the power
series expansion (for $f_1$) of the form $(\star)$ that the constant
term $a_{\mathbb O} = a_{{\mathbb O},f}$ for $f$ is also the constant term
$a_{{\mathbb O},f_1}$ for $f_1$, 
i.e.,
$$a_{\mathbb O} = a_{{\mathbb O},f_1}.$$
On the other hand, we have by Remark \ref{2.1.2.2} (1)
$$r = \on{ord}(a_{\mathbb O}) = \on{ord}(a_{{\mathbb O},f_1}) \geq \on{ord}(f_1)
\geq r + 1,$$
a contradiction !

Therefore, we have
$$\on{ord}_{\mathcal H}(f) = \on{ord}(a_{\mathbb O}).$$
This completes the proof of Lemma \ref{3.2.1.1}.
\end{proof}

\begin{rem}\label{3.2.1.2} (1) We give the justification of Remark
\ref{3.1.1.2} (1), using Lemma \ref{3.2.1.1}. 

Suppose $e_l > 0$ for some $l
= 1,
\ldots, N$.  

Take $a, b \in {\bZ}_{> 0}$ such that $a + b = p^{e_l}$.  Set
$f = x_l^a$ and $g = x_l^b$.  

Then $a_{{\mathbb O},f} = x_l^a$ and $a_{{\mathbb
O},g} = x_l^b$.  Therefore, we have by Lemma \ref{3.2.1.1}
$$\on{ord}_{\mathcal H}(f) + \on{ord}_{\mathcal H}(g) = \on{ord}(a_{{\mathbb
O},f}) + \on{ord}(a_{{\mathbb O},g}) = a + b.$$
On the other hand, we observe
$$fg = x_l^ax_l^b = x_l^{p^{e_l}} \in {\fm}_P^{p^{e_l}+1} + (h_l) \subset
{\fm}_P^{p^{e_l}+1} + ({\mathcal H}),$$
which implies
$$\on{ord}_{\mathcal H}(fg) \geq p^{e_l} + 1 > p^{e_l} = a + b =
\on{ord}_{\mathcal H}(f) + \on{ord}_{\mathcal H}(g).$$

(2) We remark that the above interpretation of $\on{ord}_{\mathcal H}(f)$ is
still valid, even if we consider the power series
expansion of the form $(\star)$ with respect to $H = (h_1, \ldots, h_N)$ and
a regular system of parameters only weakly-associated to $H$
(cf. Remark
\ref{2.1.2.2} (2), instead of the power series expansion of the form
$(\star)$ with respect to $H$ and a regular system of parameters
associated to
$H$ as described in Lemma
\ref{2.1.2.1}.
\end{rem}

\end{subsection}

\begin{subsection}{Alternative proof to Coefficient Lemma.}\label{3.2.2} The
interpretation given in \ref{3.2.1} allows us to derive
Coefficient Lemma (Lemma 4.1.4.1 in Part I) as a corollary to the formal
coefficient lemma (Lemma
\ref{2.2.2.1} in Part II).

\begin{cor}{$(= \on{Coefficient\ Lemma})$}\label{3.2.2.1} Let $\nu \in
{\bR}_{\geq 0}$ be a nonnegative number such that $\nu < {\mu}_{\mathcal
H}({\bI}_p)$.  Set
$$({\bI}_P)_t' = ({\bI}_P)_t \cap {\fm}_P^{\lceil \nu t\rceil},$$
where we use the convention that ${\fm}_P^n = R_P$ for $n \leq 0$.  Then for
any $a \in {\bR}$, we have
$$({\bI}_P)_a = \sum_B ({\bI}_P)_{a - |[B]|}'H^B.$$
\end{cor}

\begin{proof} Note that we already gave a proof to Coefficient Lemma in Part
I.  Here we present a different proof based upon the formal coefficient
lemma, although both proofs share some common spirit.

Since $H^B \in ({\bI}_P)_{|[B]|}$, we clearly have the inclusion
$$({\bI}_P)_a \supset \sum_B ({\bI}_P)_{a - |[B]|}'H^B.$$
Therefore, we have only to show the opposite inclusion
$$({\bI}_P)_a \subset \sum_B ({\bI}_P)_{a - |[B]|}'H^B.$$
Now, as observed in Remark 4.1.4.2 (2) in Part I, we have
$$\sum_B ({\bI}_P)_{a - |[B]|}'H^B = \sum_{|[B]| < a + p^{e_N}}({\bI}_P)_{a
- |[B]|}'H^B.$$
Therefore, actually we have only to show
$$({\bI}_P)_a \subset \sum_{|[B]| < a + p^{e_N}}({\bI}_P)_{a - |[B]|}'H^B.$$
Since $\widehat{R_P}$ is faithfully flat over $R_P$, we have only to prove
this inclusion at the level of completion.  That is to say, we have only
to show
$$(\widehat{{\bI}_P})_a \subset \sum_{|[B]| < a +
p^{e_N}}(\widehat{{\bI}_P})_{a - |[B]|}'H^B,$$
noting
$$\left\{
\begin{array}{ll}
({\bI}_P)_t \otimes_{R_P}\widehat{R_P} &= (\widehat{{\bI}_P})_t \\
({\bI}_P)_t' \otimes_{R_P}\widehat{R_P} &= \left(({\bI}_P)_t \cap
{\fm}_P^{\lceil \nu t\rceil}\right) \otimes_R\widehat{R_P} =
(\widehat{{\bI}_P})_t
\cap \widehat{{\fm}_P}^{\lceil \nu t\rceil} = (\widehat{{\bI}_P})_t'.\\
\end{array}
\right.$$

Take $f \in (\widehat{{\bI}_P})_a$.

Let $f = \sum_B a_BH^B$ be the power series expansion of the form $(\star)$
as described in Lemma \ref{2.1.2.1}.

Observe that, for each $C \in ({\bZ}_{\geq 0})^N$ with $|[C]| \geq a +
p^{e_N}$, there exists $B_C \in ({\bZ}_{\geq 0})^N$ with $a \leq |[B_C]| < a
+
p^{e_N}$ such that $B_C < C$ (cf. Remark 4.1.4.2 (2) in Part I).  We choose
one such $B_C$ and call it $\phi(C)$.

For each $B \in ({\bZ}_{\geq 0})^N$ with $a \leq |[B]| < a + p^{e_N}$, we
set
$$a_B' = a_B + \sum_{C \on{with} \phi(C) = B}a_CH^{C - B}.$$
We see then
$$a_B' \in \widehat{R_P} = (\widehat{{\bI}_P})_{a - |[B]|}',$$
since $a - |[B]| \leq 0$.

On the other hand, for each $B \in ({\bZ}_{\geq 0})^N$ with $|[B]| < a$, we
have by the formal coefficient lemma
$$a_B \in (\widehat{{\bI}_P})_{a - |[B]|}.$$
We also have by Lemma \ref{3.2.1.1}
$$\on{ord}(a_B) = \on{ord}_{\mathcal H}(a_B) \geq \lceil {\mu}_{\mathcal
H}({\bI}_P) (a - |[B]) \rceil \geq \lceil \nu (a - |[B]|) \rceil.$$
Therefore, we see
$$a_B \in (\widehat{{\bI}_P})_{a - |[B]|} \cap \widehat{{\fm}_P}^{\lceil \nu
(a - |[B]|) \rceil} = (\widehat{{\bI}_P})_{a - |[B]|}'.$$
We conclude

\begin{eqnarray*}
f &=& \sum_B a_BH^B = \sum_{|[B]| < a}a_BH^B + \sum_{a \leq |[B]| < a +
p^{e_N}}a_B'H^B \\
&\subset& \sum_{|[B]| < a + p^{e_N}}(\widehat{{\bI}_P})_{a - |[B]|}'H^B.\\
\end{eqnarray*}

Since $f \in (\widehat{{\bI}_P})_a$ is arbitrary, we finally conclude
$$(\widehat{{\bI}_P})_a \subset \sum_{|[B]| < a +
p^{e_N}}(\widehat{{\bI}_P})_{a - |[B]|}'H^B.$$

This completes the ``different'' proof of Coefficient Lemma.
\end{proof}

\end{subsection}

\begin{subsection}{Alternative proof to Proposition
\ref{3.1.2.1}.}\label{3.2.3}

The interpretation given in \ref{3.2.1} also allows us to provide an
alternative proof to Proposition \ref{3.1.2.1} via the formal
coefficient lemma (cf. Lemma \ref{2.2.2.1}).

\begin{cor}{$(= \on{Proposition\ }\ref{3.1.2.1})$}\label{3.2.3.1} Let the
setting be as described in \ref{3.1.1}.

Then ${\mu}_{\mathcal H}({\bI}_P)$ is independent of the choice of
${\mathcal H}$, 
i.e., independent of the choice of a leading generator
system
${\bH}$ for ${\bI}_P$.
\end{cor}

\begin{proof}[Alternative Proof] Let ${\mathcal H'}$ 
be the set of elements of another leading generator system 
${\mathbb H'}$.  We want to show ${\mu}_{\mathcal H'}(\bI_P) 
= {\mu}_{\mathcal H}(\bI_P)$.  By the same argument as in the 
proof of Proposition \ref{3.1.2.1}, we may assume that 
${\mathcal H}$ and ${\mathcal H'}$ share the same leading terms, 
i.e.,
$$h_l \equiv h_l' \bmod {\fm}_P^{p^{e_l} + 1}
\hskip.1in \on{for} \hskip.1in l = 1, \ldots, N.$$
Since ${\mathcal H}$ and ${\mathcal H'}$ share the same 
leading terms, we can take a regular system parameters 
$(x_1, \ldots, x_d)$ associated both to $H$ and to $H'$  
{\it simultaneously}.  In the following, when we consider 
the power series expansion of the form $(\star)$, we understand 
that it is with respect to $H$ and $(x_1, \ldots, x_d)$ or with 
respect to $H'$ and $(x_1, \ldots, x_d)$.

Now since ${\mu}_{\mathcal H'}(\bI_P) = {\mu}_{\mathcal H'}
(\widehat{\bI_P})$ and since ${\mu}_{\mathcal H}(\bI_P) = 
{\mu}_{\mathcal H}(\widehat{\bI_P})$, we have only to show 
$${\mu}_{\mathcal H'}(\widehat{\bI_P}) = {\mu}_{\mathcal H}
(\widehat{\bI_P}).$$
We observe that
$$\begin{array}{ll}
\mu_{\mathcal H'}(\widehat{\bI_P}) &= \inf\left\{\mu_{\mathcal H'}(f,a) 
= \frac{\on{ord}_{\mathcal H'}(f)}{a} \mid
(f,a) \in \widehat{\bI_P}, a > 0\right\} \\
&= \inf\left\{\frac{\on{ord}(a_{{\mathbb O},f}')}{a} \mid
(f,a) \in \widehat{\bI_P}, a > 0, f = \sum a_{B,f}'{H'}^B\right\} \\
& \hskip1in (\text{by\ the\ interpretation\ given\ in\ \ref{3.2.1}}) \\
&= \inf\left\{\frac{\on{ord}(f)}{a} \mid
(f,a) \in \widehat{\bI_P}, a > 0, f = a_{{\mathbb O},f}'\right\} \\
& \hskip1in (\text{by\ the\ formal\ coefficient\ lemma}) \\
\end{array}$$
and similarly that
$$\begin{array}{ll}
\mu_{\mathcal H}(\widehat{\bI_P}) &= \inf\left\{\mu_{\mathcal H}(f,a) 
= \frac{\on{ord}_{\mathcal H}(f)}{a} \mid
(f,a) \in \widehat{\bI_P}, a > 0\right\} \\
&= \inf\left\{\frac{\on{ord}(a_{{\mathbb O},f})}{a} \mid
(f,a) \in \widehat{\bI_P}, a > 0, f = \sum a_{B,f}H^B\right\} \\
& \hskip1in (\text{by\ the\ interpretation\ given\ in\ \ref{3.2.1}}) \\
&= \inf\left\{\frac{\on{ord}(f)}{a} \mid
(f,a) \in \widehat{\bI_P}, a > 0, f = a_{{\mathbb O},f}\right\} \\
& \hskip1in (\text{by\ the\ formal\ coefficient\ lemma}). \\
\end{array}$$

On the other hand, the condition $f = a_{{\mathbb O},f}'$ is 
equivalent to saying that $f$, as a power series in terms of 
$(x_1, \ldots, x_N, x_{N+1}, \ldots, x_d)$ is of the form $f 
= \sum b_KX^K$, with $b_K$ being a power series in terms of the remainder
$(x_{N+1},
\ldots, x_d)$ of the regular system of parameters, and with $K = (k_1,
\ldots, k_d)$ varying in the range satisfying the condition
$$0 \leq k_l \leq p^{e_l} - 1 \text{\ for\ }l = 1, \ldots, N \hskip.05in
\on{and} \hskip.05in k_l = 0 \text{\ for\ }l = N+1, \ldots, d.$$
Since the regular system of parameters $(x_1, \ldots, x_N, 
x_{N+1}, \ldots, x_d)$ is associated both to $H$ and $H'$ 
simultaneously, this condition is no different from the 
condition $f = a_{{\mathbb O},f}$.  That is to say, we have
$$f = a_{{\mathbb O},f}' \Longleftrightarrow f = a_{{\mathbb O},f}.$$
Therefore, by looking at the last expressions for $\mu_{\mathcal H'}
(\widehat{\bI_P})$ and $\mu_{\mathcal H}(\widehat{\bI_P})$ above, 
we conclude
$${\mu}_{\mathcal H'}(\widehat{\bI_P}) 
= {\mu}_{\mathcal H}(\widehat{\bI_P}).$$

This completes the presentation of the alternative proof.
\end{proof}

\end{subsection}

\end{section}

\begin{section}{Upper semi-continuity of
$(\sigma,\widetilde{\mu})$.}\label{3.3}

The purpose of this section is to establish the upper 
semi-continuity of $(\sigma,\widetilde{\mu})$, where the pair 
is endowed with the lexicographical order.

Recall that we have a ${\fD}$-saturated idealistic filtration ${\bI}$ over
$R$.

\begin{subsection}{Statement of the upper semi-continuity of
$(\sigma,\widetilde{\mu})$ and its proof.}\label{3.3.1} 

\begin{prop}\label{3.3.1.1} The function
$$(\sigma,\widetilde{\mu}): X = {\fm}\text{-}\on{Spec}\hskip.02in R
\rightarrow \left(\prod_{e \in {\bZ}_{\geq 0}}{\bZ}_{\geq 0}\right) \times
\left({\bR}_{\geq 0}\cup\{\infty\}\right)$$
is upper semi-continuous with respect to the lexicographical order on 
$\left(\prod_{e \in {\bZ}_{\geq 0}}{\bZ}_{\geq 0}\right) 
\times ({\bR}_{\geq 0}\cup$
$\{\infty\})$.
That is to say, for any $(\alpha,\beta) \in \left(\prod_{e \in {\bZ}_{\geq
0}}{\bZ}_{\geq 0}\right) \times \left({\bR}_{\geq 0} \cup
\{\infty\}\right)$, the locus $X_{\geq (\alpha,\beta)}$ is closed 
(cf.~Definition \ref{1.2.1.1}).

Assume further that the idealistic filtration ${\bI}$ is of r.f.g. type 
(cf.~Definition 2.1.1.1 (4) and \S 2.3 in Part I).  Then the invariant
$\widetilde{\mu}$ takes the rational values with some bounded denominator
$\delta$, and this upper semi-continuity allows us to extend the domain to
define the function
$$(\sigma,\widetilde{\mu}): \on{Spec}\hskip.02in R \rightarrow
\left(\prod_{e \in {\bZ}_{\geq 0}}{\bZ}_{\geq 0}\right) \times
\left({\bR}_{\geq
0}
\cup
\{\infty\}\right),$$
where for $Q \in \on{Spec}\hskip.02in R$ we have by definition
$$(\sigma,\widetilde{\mu})(Q) = \min\left\{(\sigma,\mu)(P) =
(\sigma(P),\widetilde{\mu}(P)) \mid P \in {\fm}\text{-}\on{Spec}\hskip.02in
R, P \in \overline{Q}\right\}$$
or equivalently $(\sigma,\widetilde{\mu})(Q)$ is equal to
$(\sigma,\widetilde{\mu})(P)$ with $P$ being a general closed point on
$\overline{Q}$.  The
function $(\sigma,\widetilde{\mu})$ with the extended domain is upper
semi-continuous.  

Moreover, since $\on{Spec}\hskip.02in R$ is noetherian and since
$\left(\prod_{e \in {\bZ}_{\geq 0}}{\bZ}_{\geq 0}\right) \times
\left({\bR}_{\geq
0}
\cup
\{\infty\}\right)$ can be replaced with a well-ordered set $T$ 
(e.g., $T$ can be obtained by replacing the
first factor $\prod_{e \in {\bZ}_{\geq 0}}{\bZ}_{\geq 0}$ 
with the well-ordered set as described in the proof of Corollary
\ref{1.2.1.3} and the second factor \text{${\bR}_{\geq
0}
\cup
\{\infty\}$} with $\frac{1}{\delta}{\bZ}_{\geq
0}
\cup
\{\infty\}$), conditions (i) and (ii) in Lemma \ref{1.2.1.2}, as well as the
assertions in Corollary \ref{1.2.1.4} hold for the upper semi-continuous
function $(\sigma,\widetilde{\mu}): \on{Spec}\hskip.02in R \rightarrow T$.
\end{prop}

\begin{proof} First we show the upper semi-continuity of the function
$$(\sigma,\widetilde{\mu}): X = {\fm}\text{-}\on{Spec}\hskip.02in R
\rightarrow \left(\prod_{e \in {\bZ}_{\geq 0}}{\bZ}_{\geq 0}\right) \times
\left({\bR}_{\geq 0}
\cup
\{\infty\}\right).$$

We have only to show that, for any $(\alpha,\beta) \in \left(\prod_{e \in
{\bZ}_{\geq 0}}{\bZ}_{\geq 0}\right) \times \left({\bR}_{\geq 0} \cup
\{\infty\}\right)$, the locus $X_{\geq (\alpha,\beta)}$ is closed.

\vskip.1in

\begin{step}{Reduction to the (local) situation where 
$X = {\fm}\text{-}\on{Spec}\hskip.02in R$ is an affine open 
neighborhood of a fixed point $P$, $\alpha = \sigma(P) \neq 
{\mathbb O}$ is the maximum of the invariant $\sigma$, and where 
a leading generator system ${\bH}$ of $\bI_P$ is uniformly pure
along the (local) maximum locus
$C$ of the invariant
$\sigma$.}

Observe $X_{\geq (\alpha,\beta)} = X_{> \alpha} \cup (X_{\leq \alpha} \cap
X_{\geq (\alpha,\beta)})$.  Since $X_{> \alpha}$ is a closed subset by the
upper semi-continuity of the invariant $\sigma$ (cf. Corollary \ref{1.2.1.3}
and Proposition \ref{1.2.2.1}), we have only to show $X_{\geq
(\alpha,\beta)}$ is closed inside of the open subset $X_{\leq \alpha}$, or
equivalently inside of any affine open subset $U$ contained in $X_{\leq
\alpha}$.  By replacing $X$ with $U$, we may assume that the invariant 
$\sigma$ never exceeds $\alpha$ in $X$.  Then again by the upper
semi-continuity of the invariant
$\sigma$, the maximum locus $C = \{Q \in X \mid \sigma(Q) = \alpha\} (=
X_{\geq \alpha})$ of the invariant
$\sigma$ is a closed subset.  Since $X_{\geq (\alpha,\beta)} \subset C$, 
we have only to show that, for any point $P \in C$, there exists an 
affine open neighborhood $U_P$ of $P$ such that $U_P \cap X_{\geq 
(\alpha,\beta)}$ is closed.

Suppose $\alpha = \sigma(P) = {\mathbb O}$.  Then, taking $U_P$ 
sufficiently small, we have $U_P \cap \on{Supp}(\bI) = \emptyset$ 
or $\{P\}$.  Therefore, we conclude that $U_P \cap X_{\geq 
(\alpha,\beta)} = U_P, \{P\}$ or $\emptyset$, and hence is closed.  
(Note that, for a point $Q \in U_P$, we have 
$(\sigma,\widetilde{\mu})(Q) = ({\mathbb O},0)$ if 
$Q \notin \on{Supp}(\bI)$, and $(\sigma,\widetilde{\mu})(Q) 
= ({\mathbb O},\infty)$ if $Q \in \on{Supp}(\bI)$.)

Therefore, in the following, we may concentrate on the case 
where $\alpha = \sigma(P) \neq {\mathbb O}$.  We take a leading 
generator system $\bH$ of ${\bI}_P$.  By Proposition
\ref{1.3.3.3} and by shrinking $U_P$ if necessary, we may assume 
that $\bH$ is uniformly pure along $C$.  Note that $C = C \cap 
\on{Supp}(\bI)$, due to the condition $\sigma(P) \neq {\mathbb O}$ 
(cf. Remark \ref{1.3.3.2}).

Finally by replacing $X$ with $U_P$, we are reduced to the 
(local) situation as described in Step 1.

We may also assume by shrinking $U_P$ if necessary, after taking a regular
system of parameters $(x_1, \ldots, x_d)$ 
associated to $H = (h_1, \ldots, h_N)$
at
$P$, that we have a regular system of parameters $(x_1, \ldots, x_d)$ over
$\on{Spec}\hskip.02in R$ such that the matrices
$$\left[\partial_{x_i^{p^e}}(h_j^{p^{e - e_j}})\right]_{i = 1, \ldots,
L_e}^{j = 1, \ldots, L_e} \hskip.05in \on{for} \hskip.05in e = e_1,
\ldots, e_N \hskip.05in \on{where}
\hskip.05in L_e =
\#\{l
\mid e_l
\leq e\}$$
are all invertible, and hence that the conditions described in the setting
4.1.1 of Part I for the supporting lemmas to
hold are satisfied (at any point in $C$).  

(We would like to bring the attention of the reader to the difference in 
notation between here in Part II and there in 4.1.1 of Part I.  The 
symbol ``$R$'' here denotes the coordinate ring of an affine open subset 
$\on{Spec}\hskip.02in R$ in $W$ (cf. the beginning of Chapter 3), while 
the symbol ``$R$'' there denotes the local ring at a closed point.)

\end{step}

\vskip.1in

\begin{step}{Reduction to statement $(\spadesuit)$, which is further reduced
to statement $(\heartsuit)$.}

We observe that, in order to provide an argument for the upper
semi-continuity, it suffices to prove the following slightly more general
statement
$(\spadesuit)$ (which does not involve any idealistic filtration):

\vskip.1in

$(\spadesuit)$ \hskip.1in Let $C \subset {\fm}\text{-}\on{Spec}\hskip.02in
R$ be a closed subset.

Let ${\mathcal H} = \{h_1, \ldots, h_N\} \subset R$ be a subset consisting
of $N$ elements, and 
$0\leq e_1\leq\cdots \leq e_N$ 
nonnegative integers attached to these elements,
satisfying the following
conditions at each point $P \in C$ (cf. 4.1.1 in Part I):
\begin{center} 
\begin{tabular}{ll}
(i) & $h_l \in {\fm}_P^{p^{e_l}}$ and $\overline{h_l} = (h_l \bmod  
{\fm}_P^{p^{e_l}+1}) = v_l^{p^{e_l}}$ with
$v_l
\in {\fm}_P/{\fm}_P^2$ for $l = 1, \ldots, N$,
\\ 
(ii) & $\{v_l \mid l = 1, \ldots, N\} \subset {\fm}_P/{\fm}_P^2$ consists of
$N$-distinct and $k$-linearly independent \\
&elements in the
$k$-vector space ${\fm}_P/{\fm}_P^2$.
\\ 
\end{tabular} 
\end{center} 

We also have a regular system of parameters $(x_1, \ldots, x_d)$ over
$\on{Spec}\hskip.02in R$ such that the matrices
$$\left[\partial_{x_i^{p^e}}(h_j^{p^{e - e_j}})\right]_{i = 1, \ldots,
L_e}^{j = 1, \ldots, L_e} \hskip.05in \on{for} \hskip.05in e = e_1,
\ldots, e_N \hskip.05in \on{where}
\hskip.05in L_e =
\#\{l
\mid e_l
\leq e\}$$
are all invertible.

Then for any $f \in R$ and $r \in {\bZ}_{\geq 0}$ the locus
$$V_r(f,{\mathcal H}) := \left\{P \in C \mid f \in {\fm}_P^rR_P 
+ ({\mathcal H})R_P\right\} = \left\{P \in C \mid 
\on{ord}_{\mathcal H}(f)(P) \geq r\right\}$$
is a closed subset.

\vskip.1in

In fact, if we prove statement $(\spadesuit)$, then
$$X_{(\alpha,\beta)} = \bigcap_{(f,a) \in {\bI}, a > 0}V_{\lceil \beta a
\rceil}(f,{\mathcal H})$$
is closed for any $(\alpha,\beta) \in \left(\prod_{e \in {\bZ}_{\geq
0}}{\bZ}_{\geq 0}\right) \times \left({\bR}_{\geq 0} \cup
\{\infty\}\right)$, and hence we have the required upper semi-continuity of
the function $(\sigma,\widetilde{\mu})$.

Furthermore, in order to prove statement $(\spadesuit)$, we have only to
show the following statement $(\heartsuit)$ for any $f \in R$ and $r \in
{\bZ}_{\geq 0}$:

$(\heartsuit)$ \hskip.1in There exist $\omega_l \in R \hskip.03in (l = 1,
\ldots, N)$ such that
$$V_r(f,{\mathcal H}) = \left\{P \in C \mid f - \sum_{l = 1}^N \omega_lh_l \in
{\fm}_P^rR_P\right\}.$$

In fact, if we show statement $(\heartsuit)$, then $V_r(f,{\mathcal H})$
being a closed set follows from the usual upper semi-continuity of the
order function for $f - \sum_{l = 1}^N \omega_lh_l$, and hence we have
statement $(\spadesuit)$.

\end{step}

\vskip.1in

\begin{step}{Show statement $(\heartsuit)$ by induction on $r$.}

Step 3 is dedicated to showing statement $(\heartsuit)$ by induction on $r$.

\vskip.1in 

We set

$$\left\{
\begin{array}{ll}
e &:= e_1 = \min\left\{e_l \mid l = 1, \ldots, N\right\}, \\
L &:= \max\left\{l \mid l = 1, \ldots, N, e_l = e\right\} 
= \#\left\{l \mid 1, \ldots, N,
e_l = e\right\}, \\
e' &= e_{L+1} \hskip.1in (\on{if\ }L = N, \on{\ then\ we\ set\ }e' =
\infty),\\
\chi &= \#\left\{e_1, \ldots, e_N\right\}.\\
\end{array}\right.
$$

\vskip.1in

\begin{case}{$r \leq p^e$}

In this case, we have only to set $\omega_l = 0$ \hskip.02in ($l = 1, 
\ldots, l$) in order to see statement $(\heartsuit)$.
\end{case}

\begin{case}{$r > p^e$}

Observing $V_r(f,{\mathcal H}) \subset V_{r-1}(f,{\mathcal H})$ and 
replacing $f$ with $f - \sum_{l = 1}^N \omega_lh_l$ via application of 
statement $(\heartsuit)$ for $r-1$ by
induction, we may assume
$$f \in {\fm}_P^{r-1}R_P \hskip.1in \forall P \in V_r(f,{\mathcal H}).$$
We also observe then, by Supporting Lemma 3 in Part I 
(cf. Lemma 4.1.2.3 in Part I), that, at each point 
$P \in V_r(f,{\mathcal H})$, there exist
$\beta_{l,P} \in {\fm}_P^{r-1-p^{e_l}}R_P$ such that
$$f - \sum_{l = 1}\beta_{l,P}h_l \in {\fm}_P^rR_P.$$
Now we use the induction on the pair $(\chi, L)$.

\vskip.1in

$\boxed{\on{Case}: \chi = 1 \hskip.05in (L = N, e' = \infty)}$

\vskip.1in

In this case, by applying Supporting Lemma 2 in Part I (cf. Lemma 4.1.2.2 in 
Part I) with $v = r, s = r-1$ and $\alpha = - f$, we see

\begin{eqnarray*}
(*) \hskip.1in \beta_L &\in& F_v(-f) + \sum_{l = 1}^{N - 
1}(F_v\beta_{l,P})h_l + \left(h_L^r\right) + {\fm}_P^{r - p^e}R_P \\
&\subset& F_v(-f) + \sum_{l = 1}^{N - 1}(F_v\beta_{l,P})h_l + 
{\fm}_P^{r - p^e}R_P.\\
\end{eqnarray*}

See Supporting Lemmas 1 and 2 in Part I (cf. Lemma 4.1.2.1 and Lemma 4.1.2.2 
in Part I) for the definition of the differential operator $F_v$.  We
would like to emphasize that, even though Supporting Lemma 3 is a local 
statement at $P$, the differential operator $F_v$ is defined globally 
over $\on{Spec}\hskip.02in R$ and
hence that $F_v(-f) \in R$.

{}From $(*)$, we conclude the following.

When $N = 1$, we have only to set $\omega_1 = F_v(- f)$ in order to see 
statement $(\heartsuit)$.

When $N > 1$, we observe 
$$V_r(f,{\mathcal H}) = V_r(f, \{h_1, \ldots, h_{N-1}, h_N\}) = V_r(f - 
F_v(-f)h_N, \{h_1, \ldots, h_{N-1}\}).$$  
Now statement $(\heartsuit)$ for $f$ and $r$ with respect
to
${\mathcal H} = \{h_1, \ldots, h_{N-1}, h_N\}$ follows from statement 
$(\heartsuit)$ for $f - F_v(-f)h_N$ and $r$ with respect to $\{h_1, 
\ldots,
h_{N-1}\}$, which holds by induction on $(\chi,L) = (1,N-1)$.

\vskip.1in

$\boxed{\on{Case}: \chi > 1}$

\vskip.1in

In this case, by applying Supporting Lemma 2 in Part I (cf. Lemma 4.1.2.2 in 
Part I) with $v = p^{e'-e} - 1, s = r-1$ and $\alpha = - f$, we see

\begin{eqnarray*}
(*) \hskip.1in \beta_L &\in& F_v(-f) + \sum_{1 \leq l \leq N, l \neq 
L}(F_v\beta_{l,P})h_l + \left(h_L^v\right) + {\fm}_P^{r - p^e}R_P \\
&\subset& F_v(-f) + \sum_{1 \leq l \leq N, l \neq L}h_lR_P + + 
\left(h_L^v\right) + {\fm}_P^{r - p^e}R_P.\\
\end{eqnarray*}

{}From $(*)$, we conclude that
$$f - F_v(-f)h_L \in \sum_{1 \leq l \leq N, l \neq L}h_lR_P + 
h_L^{p^{e'-e}}R_P + {\fm}_P^rR_P$$
and hence that

\begin{eqnarray*}
V_r(f,{\mathcal H}) &=& V_r(f,\{h_1, \ldots, h_{L-1}, h_L, h_{L+1}, \ldots, 
h_N\}) \\
&=& V_r(f - F_v(-f), \{h_1, \ldots, h_{L-1}, h_L^{p^{e'-e}},
h_{L+1}, \ldots, h_N\}).\\
\end{eqnarray*}

Now statement $(\heartsuit)$ for $f$ and $r$ with respect to
${\mathcal H} = \{h_1, \ldots, h_{L-1}, h_L, h_{L+1}, \ldots, h_N\}$ 
follows from statement 
$(\heartsuit)$ for $f - F_v(-f)h_N$ and $r$ with respect 
to $\{h_1, \ldots, h_{L-1}, h_L^{p^{e'-e}}, h_{L+1},$
$\ldots, h_N\}$, 
which holds by induction on $(\chi,L)$.  (In fact, 
if originally $L = 1$, then the invariant $\chi$ drops by $1$, and if 
originally $L > 1$, then the
invariant $\chi$ remains the same but the invariant $L$ drops by $1$.)
\end{case}

This completes the proof of statement $(\heartsuit)$.

\end{step}

\vskip.1in

This completes the proof of the upper semi-continuity of the function
$$(\sigma,\widetilde{\mu}): X = {\fm}\text{-}\on{Spec}\hskip.02in R 
\rightarrow \left(\prod_{e \in {\bZ}_{\geq 0}}{\bZ}_{\geq 0}\right) \times 
\left({\bR}_{\geq
0}
\cup
\{\infty\}\right).$$

\vskip.1in

If we assume further that the idealistic filtration is of r.f.g. type, then, 
as shown in Remark \ref{3.1.1.2} (2), the invariant $\widetilde{\mu}$ takes
the rational values with some bounded denominator $\delta$.  
Then we may replace the target space for the function 
$$(\sigma,\widetilde{\mu}): 
{\fm}\text{-}\on{Spec}\hskip.02in R
\rightarrow \left(\prod_{e \in {\bZ}_{\geq 0}}{\bZ}_{\geq 0}\right) \times
\left({\bR}_{\geq 0}
\cup
\{\infty\}\right)$$ with a well-ordered set $T$ (e.g., $T$ can be 
obtained by replacing the first factor $\prod_{e \in {\bZ}_{\geq 0}}
{\bZ}_{\geq 0}$ with the well-ordered set as described in the proof 
of Corollary \ref{1.2.1.3} and the second factor
${\bR}_{\geq 0}\cup \{\infty\}$ with $\frac{1}{\delta}{\bZ}_{\geq0}
\cup
\{\infty\}$).  Now the assertion regarding the extension of the domain of 
the function $(\sigma,\widetilde{\mu})$ from 
${\fm}\text{-}\on{Spec}\hskip.02in R$ to
$\on{Spec}\hskip.02in R$ and the rest of the assertions in 
Proposition \ref{3.3.1.1} follow from the same 
argument as in Corollary \ref{1.2.2.2}, where we discussed the extension of 
the domain of the invariant
$\sigma$ from ${\fm}\text{-}\on{Spec}\hskip.02in R$ to
$\on{Spec}\hskip.02in R$.

\vskip.1in 

This completes the proof of Proposition \ref{3.3.1.1}.
\end{proof}

\end{subsection}

\begin{subsection}{Alternative proof to the upper semi-continuity of 
$(\sigma,\widetilde{\mu})$.}\label{3.3.2} We can give an alternative proof 
to the
upper semi-continuity of
$(\sigma,\widetilde{\mu})$, using the interpretation of $\widetilde{\mu}$ in 
terms of the power series expansion of the form $(\star)$ as presented 
in \ref{3.2}.  

\begin{proof}[Alternative Proof to the upper semi-continuity of 
$(\sigma,\widetilde{\mu})$] By the same argument as before, we are 
reduced to the (local) situation as described in Step 1 of the original proof.

Take a regular system of parameters $X_P = (x_1, \ldots, x_d) 
= (x_{1,P}, \ldots, x_{d,P})$ at $P$, which 
is associated to $H = (h_1, \ldots, h_N)$.  By shrinking 
$\on{Spec}\hskip.02in R$ if necessary, we may assume that 
$X_Q = (x_{1,Q}, \ldots, x_{d,Q}) = (x_1 
- x_1(Q), \ldots, x_d - x_d(Q))$ with $x_i(Q) \in k$ is a 
regular system of 
parameters, which is weakly-associated to $H$ at any point $Q \in C$.

\vskip.1in

By the same argument as before, we have only to show 
that, given $f \in R$ and $r \in {\bZ}_{\geq 0}$, the locus $V_r(f,{\mathcal 
H}) = \left\{Q
\in C \mid \on{ord}_{\mathcal H}(f)(Q) \geq r\right\}$ is a 
closed subset as in Step 2 of the original proof.

This is where the alternative argument using the interpretation of 
$\widetilde{\mu}$ presented in \ref{3.2} begins: Let 
$f = \sum a_{B,Q}H^B$ be the power series expansion of $f$ at $Q \in C$ 
with respect to $H$ and the regular system of parameters $X_Q$, which is
weakly-associated to $H$ at $Q$.  By Lemma \ref{3.2.1.1} and Remark 
\ref{3.2.1.2} (2), we have
$$\on{ord}_{\mathcal H}(f)(Q) = \on{ord}(a_{{\mathbb O},Q}).$$
Let 
$$a_{{\mathbb O},Q} = \sum \gamma_{I,Q}X_Q^I, \hskip.1in\gamma_{I,Q} \in k$$ 
be the power series expansion of $a_{{\mathbb O},Q}$ with respect to $X_Q$.  
Then we have
$$\on{ord}(a_{{\mathbb O},Q}) \geq r \Longleftrightarrow \gamma_{I,Q} = 0 
\hskip.1in \forall I \hskip.03in \on{with} \hskip.03in |I| < r.$$
On the other hand, since the coefficients $\gamma_{I,Q}$ can be computed 
from the coefficients of the power series expansions 

$$\left\{\begin{array}{ll}
f &= \sum c_{I,f,Q}X_Q^I \text{ with }c_{I,f,Q} \in k\\
&\hskip.76in  \text{where }c_{I,f,Q} = \partial_{X_Q^I}(f)(Q) 
= \partial_{X_P^I}(f)(Q) \\
h_l &= \sum c_{I,h_l,Q}X_Q^I \text{ with }c_{I,h_l,Q} \in k\\
&\hskip.76in  \text{where }c_{I,h_l,Q} = \partial_{X_Q^I}(h_l)(Q) 
= \partial_{X_P^I}(h_l)(Q) \\
& \hskip2in \text{for }l = 1, \ldots, N, \\
\end{array}\right.$$ 
via the invertible 
matrices appearing in the condition of $X_Q$ being weakly-associated to
$H$
$$\left[\partial_{x_{i,Q}^{p^e}}(h_j^{p^{e - e_j}})\right]_{i = 1, \ldots,
L_e}^{j = 1, \ldots, L_e}(Q) = \left[\partial_{x_{i,P}^{p^e}}
(h_j^{p^{e - e_j}})\right]_{i = 1, \ldots,
L_e}^{j = 1, \ldots, L_e}(Q) \qquad\on{\ for\ } e = e_1,
\ldots, e_N $$
where $L_e = \#\{l\mid e_l\leq e\}$, 
we conclude that, for each $I$, there exists $\gamma_I \in R$ such that 
$\gamma_I(Q) = \gamma_{I,Q} \hskip.1in \forall Q \in C$.

Finally then we conclude that
$$V_r(f,{\mathcal H}) = \left\{Q \in C \mid \on{ord}_{\mathcal H}(f)(Q) 
\geq r\right\} = \left\{Q \in C \mid \gamma_I(Q) = 0 \hskip.1in \forall 
I \hskip.05in \on{with} \hskip.05in |I| < r\right\}$$
is a closed subset.

This completes the alternative proof for the upper 
semi-continuity of $(\sigma,\widetilde{\mu})$.

\end{proof}

\end{subsection}
\end{section}
\end{chapter} 

\begin{chapter}*{Appendix}
\def\thechapter{A}
\setcounter{section}{0}
\renewcommand{\thesection}{\S\thechapter.%
\arabic{section}}
\renewcommand{\thesubsection}{\thechapter.%
\arabic{section}.\arabic{subsection}}
\renewcommand{\thesubsubsection}{\thechapter.%
\arabic{section}.\arabic{subsection}.\arabic{subsubsection}}

The purpose of this appendix is to present the new nonsingularity 
principle using only the $\fD$-saturation, as opposed to the old 
nonsingularity principle using both the $\fD$-saturation and 
$\fR$-saturation.  (The combination of the $\fD$-saturation and 
$\fR$-saturation was called the \text{$\fB$-saturation} in Part I 
(cf. 2.1.5 and 2.2.3 in Part I).)

In Part I, we emphasized the importance of the $\fR$-saturation (and of the
$\fB$-saturation) in carrying out the
IFP.  In fact, the $\fR$-saturation was
crucial in establishing the nonsingularity principle, as formulated in
Chapter 4 of Part I, which sits at the heart of constructing an algorithm.
However, the
$\fR$-saturation has also been the main culprit in our quest to complete the
algorithm, causing the following problems:

\begin{itemize}

\item By
taking the $\fR$-saturation, we may increase the denominator of the
invariant $\widetilde{\mu}$ indefinitely, and hence may not have the 
descending chain condition on the value set of the strand of 
invariants consisting of the units of the form $(\sigma, \widetilde{\mu}, s)$.
This invites the problem of termination, as we mentioned in the introduction
to Part I.  

\item If we take the $\fR$-saturation, the value of the invariant  
    $\widetilde{\mu}$ may strictly increase under blowup, even when  
    the value of the invariant $\sigma$ stays the same.  This violates  
    the principle that our strand of invariants, consisting of the  
    units of the form $(\sigma,\widetilde{\mu},s)$, should never  
    increase under blowup. 

\end{itemize}

While writing Part II, we came to realize that we
can establish the nonsingularity principle, as formulated below, with only
the $\fD$-saturation and without the
$\fR$-saturation.  This indicates that we may construct an algorithm, still
in the frame work of the IFP, without using the $\fR$-saturation, and hence
that we may avoid the problem of
termination, as well as the other technical problems, that the use of 
the $\fR$-saturation invites.

Even though we are still in the evolution process of the IFP program 
(See \ref{0.3.1} for the current status of the IFP.), we consider this 
new nonsingularity principle a substantial step forward in our quest to 
construct an algorithm for local and global resolution of singularities 
in positive characteristic.

\vskip.1in

In this appendix, $R$
represents the coordinate ring of an affine open subset
$\on{Spec}\hskip.02in R$ of a nonsingular variety $W$ of $\dim W = d$ over
an algebraically closed  field 
$k$ of positive characteristic ${\on{char}}(k) = p$ or of
characteristic zero ${\on{char}}(k) = 0$, where in the latter case we
formally set 
$p = \infty$ (cf. 0.2.3.2.1 and Definition 3.1.1.1 (2) in Part I).

\begin{section}{Nonsingularity principle with only $\fD$-saturation and
without $\fR$-saturation.}\label{A.1}

\begin{subsection}{Statement.}\label{A.1.1}

\begin{thm}\label{A.1.1.1} Let $\bI$ be an idealistic filtration over $R$.

Let $P \in \on{Spec}\hskip.02in R \subset W$ be a closed point.

\item[(1)] Assume that $\bI$ is $\fD$-saturated and that 
$\widetilde{\mu}(P) = \infty$.  

Then there exists a regular system of parameters $(x_1, \cdots, x_N,
y_{N+1},
\cdots, y_d)$ at $P$ such that

$\left\{\begin{array}{ll}
\bH &= \left\{(x_l^{p^{e_l}}, p^{e_l})\right\}_{l = 1}^N 
\hskip.1in (e_1 \leq \cdots
\leq e_N) \text{\ is\ an\ LGS\ of\ }\bI_P \\
& (\text{See\ the\ footnote\ to\ \ref{0.2.1.5}.}), \\
\bI_P &= G_{R_P}(\bH). \\
\end{array}\right.$

\item[(2)] Assume further that $\bI$ is of r.f.g. type.

Then there exists an affine open neighborhood $P \in U_P =
\on{Spec}\hskip.02in R_r$ of $P$ (Note that $R_r$ represents 
the localization of $R$ by $r \in R$.) such that $(x_1, \cdots, x_N, y_{N+1},
\cdots, y_d)$
is a regular system of parameters over $U_P$, and that

$\left\{\begin{array}{ll}
\bH &= \left\{(x_l^{p^{e_l}},p^{e_l})\right\}_{l = 1}^N \subset \bI_r, \\
\bI_r &= G_{R_r}(\bH).\\
\end{array}\right.$

In particular, we have

\begin{itemize}
\item $\on{Supp}(\bI) \cap U_P = V(x_1, \cdots, x_N)$, which is hence
nonsingular, and

\item $(\sigma(Q),\widetilde{\mu}(Q)) = (\sigma(P),\infty)$ 
for any closed point
$Q \in \on{Supp}(\bI) \cap U_P$.

\end{itemize}

\end{thm} 

\begin{rem}\label{A.1.1.2}

(1) It is straightforward to see that assertion (1) actually gives the
following characterization: An idealistic filtration $\bI_P$ over $R_P$ is
$\fD$-saturated and $\widetilde{\mu}(P) = \infty$ if
and only if there exist a regular system of parameters $(x_1, \cdots, x_N,
y_{N+1}, \cdots, y_d)$ and a subset of the form $\bH = 
\left\{(x_l^{p^{e_l}},p^{e_l})\right\}_{l =
1}^N \subset \bI_P$ such that $\bI_P =
G_{R_P}(\bH)$.  (The subset $\bH$ is then automatically an LGS of $\bI_P$.)

(2) We construct the strand of invariants in our algorithm 
(cf. 0.2.3.2.2 in the introduction to Part I), and at year 0 
it takes the following form
$$\on{inv}_{\on{new}}(P) = (\sigma_0^1, \widetilde{\mu}_0^1, 
s_0^1)(\sigma_0^2, \widetilde{\mu}_0^2, s_0^2)
\cdots (\sigma_0^{n-1}, \widetilde{\mu}_0^{n-1}, 
s_0^{n-1})(\sigma_0^n,\widetilde{\mu}_0^n,s_0^n),$$
with the last $n$-th unit $(\sigma_0^n,\widetilde{\mu}_0^n,s_0^n)$ being 
equal to $(\sigma_0^n,\infty,0)$.  The subscript ``${}_0$'' refers to 
year ``$0$'', while the superscript ``${}^{j}$'' refers to the stage 
``$j$''.  (Note that, if we insert the new invariant $\widetilde{\nu}$ 
so that the unit changes from the triplet $(\sigma,\widetilde{\mu},s)$ 
to the quadruplet $(\sigma,\widetilde{\mu},\widetilde{\nu},s)$, then 
the strand of invariants also changes accordingly (cf. \ref{0.3.1}).)

The (local) maximum locus of the strand of invariants coincides with the 
support $\on{Supp}(\bI_0^n)$ of the last $n$-th modification $\bI_0^n$. 
(Note that in year $0$ we always have $\widetilde{\mu} > 1$ and hence 
that there is no gap between the (local) maximum locus and the support 
of the modification, an anomaly observed when $\widetilde{\mu} = 1$.) 
The idealistic filtration $\bI_0^n$ is $\fD$-saturated with 
$\widetilde{\mu}(P) = \infty$.  Therefore, applying Theorem 
\ref{A.1.1.1}, we conclude that $\on{Supp}(\bI_0^n)$ is nonsingular (in 
a neighborhood of $P$).  (Note that all the idealistic filtrations we 
deal with in our algorithm are of r.f.g. type.)  Therefore, we conclude 
that the center of blowup, which is chosen to be the maximum locus of 
the strand, is nonsingular.  This is why Theorem \ref{A.1.1.1} is 
called the (new) nonsingularity principle of the center.  (After year 
$0$, we have to make several technical adjustments, including an 
adjustment to overcome the gap between the (local) maximum locus and 
the support of the last modification and another adjustment to 
introduce the $\fD_E$-saturation in the presence of the exceptional 
divisor $E$ instead of the usual $\fD$-saturation.  The basic tool for 
us to guarantee the nonsingularity of the center, however, is still 
Theorem \ref{A.1.1.1}.)

(3) If we assume further that $\bI_P$ is $\fR$-saturated, then after having
assertion (1), we immediately come to the conclusion that $e_1 = \cdots =
e_N = 0$, 
i.e., all the elements in the LGS
(and hence of any LGS) are concentrated at level $1$.  That is to say, we
obtain the old nonsingularity principle Theorem 4.2.1.1 in Part I as a 
corollary to the new nonsingularity principle Theorem \ref{A.1.1.1} of 
this appendix.

\end{rem}

\end{subsection}

\newpage

\begin{subsection}{Proof.}\label{A.1.2}
\begin{proof}[\underline{\it Proof for assertion $\on{(1)}$}]

\begin{step}{Show $\bI_P = G_{R_P}(\bH)$ for
{\it any} LGS $\bH$ of
$\bI_P$.}

First, note that, if $P \not\in \on{Supp}(\bI)$, then we would have $\bI_P =
R_P \times \bR$ (since $\bI_P$ is
$\fD$-saturated cf. Case $P \not\in \on{Supp}(\bI)$ of 
Lemma \ref{1.1.2.1}) and hence
$\widetilde{\mu}(P) = 0$.  Thus our assumption 
$\widetilde{\mu}(P) = \infty$ implies $P \in \on{Supp}(\bI)$.  
Second, we claim $\bI_P = G_{R_P}(\bH)$ for
{\it any} LGS $\bH$ of
$\bI_P$.  In order to prove this claim, we can use 
the same argument as presented in the proof of the 
nonsingularity principle formulated in Chapter
4 of Part I.  Note that this part of the proof
did not use the assumption that
$\bI_P$ is
$\fR$-saturated.  

Alternatively, we can give a proof of the claim using the formal coefficient
lemma (cf. Lemma \ref{2.2.2.1}) discussed here in Part II, without referring
to the arguments in Part I:

Take an element $f \in (\bI_P)_a \subset (\widehat{\bI_P})_a$, and let $f =
\sum_{B \in ({\bZ}_{\geq 0})^N}a_BH^B$ be the power series expansion of the
form $(\star)$ as
described in Lemma \ref{2.1.2.1}.  From the formal coefficient lemma it
follows that
$$(a_B, a - |[B]|) \in \widehat{\bI_P} \hskip.1in \forall B \in ({\bZ}_{\geq
0})^N.$$
Suppose there exists $B \in ({\bZ}_{\geq 0})^N$ with $|[B]| < a$ such that
$a_B \neq 0$.  Then we would have
$$\widetilde{\mu}(P) = \mu_{\mathcal H}(\bI_P) = \mu_{\mathcal
H}(\widehat{\bI_P}) \leq \frac{\on{ord}_{\mathcal H}(a_B)}{a - |[B]|} =
\frac{\on{ord}(a_B)}{a - |[B]|} < \infty, \hskip.1in (\on{cf.}
\text{Lemma\ } \ref{3.2.1.1})$$ a contradiction !  Therefore, we conclude
$$a_B = 0 \hskip.1in \forall B \in ({\bZ}_{\geq 0})^N \text{\ with\ } |[B]|
< a.$$
This implies 
$$f \in (H^B \mid |[B]| \geq a).$$
Since $f \in \left(\bI_P\right)_a$ is arbitrary, we conclude
$$\left(\bI_P\right)_a \subset (H^B \mid |[B]| \geq a),$$
while the opposite inclusion $\left(\bI_P\right)_a \supset (H^B \mid |[B]|
\geq a)$ is obvious.  Therefore, we finally conclude
$$\left(\bI_P\right)_a = (H^B \mid |[B]| \geq a) \hskip.1in \forall a \in
{\bR},$$
which is equivalent to saying
$$\bI_P = G_{R_P}(\bH).$$
\end{step}

\begin{step}{Inductive construction of an LGS and a regular system of
parameters of the desired form via claim $(\diamondsuit)$.}

Now, by induction, we assume that we have found an LGS $\bH =
\{(h_{ij},p^{e_i})\}$ of $\bI_P$ and a regular system of parameters
$(\{x_{ij}\}, y_{N+1}, \cdots, y_d)$ at
$P$ such that
$$\left\{\begin{array}{lcll}
h_{ij} &=& x_{ij}^{p^{e_i}} \hskip.1in &\text{\ if\ }e_i < e_u,\\
h_{ij} &=& x_{ij}^{p^{e_i}} \bmod 
{\fm_P}^{p^{e_i}+1} &\text{\ if\ }e_i \geq e_u.\\
\end{array}\right.$$
Note that we use the double subscripts $h_{ij}$ for the elements in the LGS,
where the first subscript indicates the level $p^{e_i}$ with $e_1 < \cdots <
e_M$.  So we have the total of $N$
elements at $M$ distinct levels in the LGS.  
(See \ref{1.3.1}.)  The inductive assumption means that we have 
found an LGS and a regular system of parameters of the desired 
form up to the level $e_i = e_{u-1}$.

We want to show, by replacing $h_{ij}$ and $x_{ij}$ with $e_i = e_u$ via the
use of claim $(\diamondsuit)$ which we state next, that we can also have
$$\left\{\begin{array}{lcll}
h_{ij} &=& x_{ij}^{p^{e_i}} \hskip.1in &\text{\ if\ }e_i < e_{u+1},\\
h_{ij} &=& x_{ij}^{p^{e_i}} \bmod 
{\fm_P}^{p^{e_i}+1} &\text{\ if\ }e_i \geq e_{u+1}.\\
\end{array}\right.$$
\end{step}

\vskip.1in

\begin{step}{Statement and the proof of claim $(\diamondsuit)$}

This step is devoted to proving the following claim:

\vskip.1in

$(\diamondsuit)$ \hskip.1in Set, for $l \in {\bZ}_{\geq 0}$,
$$J_l = F^{e_u}(\fm_P) + \left(X^{[C]} = \prod_{e_i < e_u}
x_{ij}^{p^{e_i}c_{ij}}\mid C = (c_{ij} \mid e_i < e_u), |[C]| \geq
p^{e_u}\right) + (\bI_P)_{p^{e_u}} \cap {\fm_P}^{p^{e_u} + 1} +
{\fm_P}^l.$$ 
(Recall that the symbol ``$F$'' represents the Frobenius map.)

Then we have the inclusion
$$(\bI_P)_{p^{e_u}} \subset J_l \hskip.1in \forall l \in {\bZ}_{\geq 0}.$$

\vskip.1in

Observe
$$\begin{array}{ll}
(\bI_P)_{p^{e_u}} &= (H^B \mid |[B]| \geq p^{e_u}) \hskip.1in 
(\text{since\ }\bI_P = G_{R_P}(\bH))\\
&= \left(X^{[C]} \mid C = (c_{ij} \mid e_i < e_u), |[C]| \geq 
p^{e_u}\right) + (h_{ij}
\mid e_i = e_u) + (h_{ij} \mid e_i > e_u) \\
&\subset \left(X^{[C]} \mid C = (c_{ij} \mid e_i < e_u), |[C]| 
\geq p^{e_u}\right) +
F^{e_u}(\fm_P) + {\fm_P}^{p^{e_u} + 1} \\
&\subset J_{p^{e_u}+1}.\\
\end{array}$$
Therefore, the required inclusion holds for $l \leq p^{e_u} + 1$.

Now assume, by induction, that the required inclusion
$$(\bI_P)_{p^{e_u}} \subset J_l$$
holds for a fixed $l \geq p^{e_u} + 1$.  We want to show
$$(\bI_P)_{p^{e_u}} \subset J_{l+1}.$$

Take an arbitrary element $f \in (\bI_P)_{p^{e_u}} \subset J_l$.

We may choose $\{\alpha_{ST} \mid S, T\} \subset k$ such that
$$f - \sum_{|(S,T)| = l}\alpha_{ST}X^SY^T \in J_{l+1}.$$
Note that then there exists $w \in \fm_P$ such that
$$(\heartsuit) \hskip.1in w^{p^{e_u}} + \sum_{|(S,T)| = l}\alpha_{ST}X^SY^T
\in ({\bI}_P)_{p^{e_u}} + {\fm_P}^{l+1}.$$

Set
$$\left\{\begin{array}{ll}
s_{ij} &= p^{e_i}s_{ij,q} + s_{ij,r} \text{\ with\ }0 \leq s_{ij,r} <
p^{e_i} \\
S_q &= (s_{ij,q}), \\
S_r &= (s_{ij,r}).\\
\end{array}\right.$$
Then we have $S = [S_q] + S_r$ and $X^SY^T = X^{[S_q]}X^{S_r}Y^T$.

We analyze the terms in
$$\sum_{|(S,T)| = l}\alpha_{ST}X^SY^T = \sum_{|(S,T)| =
l}\alpha_{ST}X^{[S_q]}X^{S_r}Y^T.$$

\pagebreak[4]

\begin{case}{$S_q = 0$.}

In this case, we write for simplicity
$$X^SY^T = X^{S_r}Y^T = Z^V$$
by setting
$$\left\{\begin{array}{ll}
(X,Y) &= (\{x_{ij}\}, y_{N+1}, \cdots, y_d) = (z_1, \cdots, z_d) = Z \\
(S,T) &= (S_r,T) = (\{s_{ij,r}\},t_{N+1}, \cdots, t_d) = (v_1, \cdots,
v_d).\\
\end{array}\right.$$

$\underline{\on{Subcase\ 1.1}}$ \hskip.02in : $p^{e_u} | \hskip.04in V$.

\vskip.1in

In this subcase, we conclude
$$\alpha_{ST}X^SY^T = \alpha_{ST}Z^V \in F^{e_u}(\fm_P) \subset J_{l+1}.$$

\vskip.1in

$\underline{\on{Subcase\ 1.2}}$ \hskip.02in : $p^{e_u} \not|  \hskip.06in
V$.

\vskip.1in

In this subcase, let $v_{\omega}$ be a factor, not divisible by $p^{e_u}$,
of $V$.

Set
$$v_{\omega} = p^{e_u}v_{\omega,q} + v_{\omega,r} \text{\ with\ }0 <
v_{\omega,r} < p^{e_u}.$$

Apply $\partial_{z_{\omega}^{v_{\omega,r}}}$ to $(\heartsuit)$ and obtain
\begin{eqnarray*}
\lefteqn{
\partial_{z_{\omega}^{v_{\omega,r}}}\left(w^{p^{e_u}} + \sum_{|(S,T)| =
l}\alpha_{ST}X^SY^T\right) 
=
\partial_{z_{\omega}^{v_{\omega,r}}}\left(\sum_{|(S,T)| =
l}\alpha_{ST}X^SY^T\right) }\\
&=& \alpha_{ST}Z^{V -
v_{\omega,r}{\mathbf e}_{\omega}} + (\text{other\ monomials\ of\ degree\ }(l
- v_{\omega,r})) \\
&\in& \left((\bI_P)_{p^{e_u} - v_{\omega,r}} 
+ {\fm_P}^{l - v_{\omega,r} + 1}\right)
\cap {\fm_P}^{l - v_{\omega,r}} \\
&=& (\bI_P)_{p^{e_u} - v_{\omega,r}} \cap {\fm_P}^{l - v_{\omega,r}} +
{\fm_P}^{l - v_{\omega,r} + 1}. 
\end{eqnarray*}
On the other hand, we observe
$$(\bI_P)_{p^{e_u} - v_{\omega,r}} \cap {\fm_P}^{l - v_{\omega,r}} \subset
\sum_{1 \leq i \leq M}h_{ij}{\fm_P}^{l - v_{\omega,r} - p^{e_i}}.$$
(We use the convention that ${\fm_P}^{l - v_{\omega,r} - p^{e_i}} = R_P$ if
$l - v_{\omega,r} - p^{e_i} \leq 0$.)

In fact, let $g \in (\bI_P)_{p^{e_u} - v_{\omega,r}} \cap {\fm_P}^{l -
v_{\omega,r}}$ be an arbitrary element, and $g = \sum a_BH^B$ the power
series expansion of the form $(\star)$ as described
in Lemma \ref{2.1.2.1}.  Then it follows 
from the condition $\widetilde{\mu}(P)
= \infty$ and $0 < p^{e_u} - v_{\omega,r}$ 
that $a_{\mathbb O} = 0$ (cf. Lemma \ref{3.2.1.1}), and 
from the construction that $\on{ord}_P(a_B)
\geq (l - v_{\omega,r}) - |[B]|$ for any $B \in ({\bZ}_{\geq 0})^N$ (cf.
Remark \ref{2.1.2.2} (1)).  Therefore, we conclude $f \in \sum_{1 \leq i
\leq M}h_{ij}{\fm_P}^{l - v_{\omega,r} -
p^{e_i}}$.  This proves the inclusion above.  (Note that the inclusion above
can also be derived using Lemma 4.1.2.3 in Part I via the fact that ${\bI}_P
= G_{R_P}(\bH)$.) 

However, this inclusion implies that any monomial of degree $l -
v_{\omega,r}$ in the power series expansion of an element in
$(\bI_P)_{p^{e_u} - v_{\omega,r}} \cap {\fm_P}^{l - v_{\omega,r}}$,
with respect to the regular system of parameters $(x_1, \cdots, x_N,
y_{N+1}, \cdots, y_d)$ should be divisible by some element in the set
$\left\{x_{ij}^{p^{e_i}} \mid 1 \leq i \leq M\right\}$, and
hence that the monomial $Z^{V - v_{\omega,r}{\mathbf e}_{\omega}}$ can not
appear as $S_q = 0$.

Therefore, in this subcase, we conclude
$$\alpha_{ST} = 0.$$
\end{case}
\pagebreak[4]
\begin{case}{$S_q \neq 0$}

$\underline{\on{Subcase\ 2.1}}$ \hskip.02in: $s_{ij,q} > 0$ for some $i \geq
u$.

\vskip.1in

In this subcase, we compute
\begin{eqnarray*}
\lefteqn{
X^SY^T \in x_{ij}^{p^{e_i}}{\fm_P}^{l - p^{e_i}} 
\subset (h_{ij} + {\fm_P}^{p^{e_i}+1}){\fm_P}^{l - p^{e_i}} }\\
&\subset& h_{ij}{\fm_P}^{l - p^{e_i}} + {\fm_P}^{l+1} 
\subset (\bI_P)_{p^{e_u}} \cap {\fm_P}^{p^{e_u}+1} + {\fm_P}^{l+1} 
\subset J_{l+1}.
\end{eqnarray*}
(Note that, in order to obtain the second last inclusion above,
we use the fact that $h_{ij} \in {\fm_P}^{p^{e_u}+1}$ if $i > u$, and use
the condition $l \geq p^{e_u} + 1$ if $i =
u$.)

Therefore, we conclude
$$\alpha_{ST}X^SY^T \in J_{l+1}.$$

\vskip.1in

$\underline{\on{Subcase\ 2.2}}$ \hskip.02in: $s_{ij,q} = 0$ for any $i \geq
u$ and $|[S_q]| \geq p^{e_u}$.

\vskip.1in

In this subcase, we conclude
\begin{eqnarray*}
\alpha_{ST} &=& \alpha_{ST}X^{[S_q]}X^{S_r}Y^T \\
&\in& \left(X^{[C]} \mid C = (c_{ij} \mid e_i < e_u), |[C]| \geq p^{e_u}\right)
\subset J_{l+1}.
\end{eqnarray*}

$\underline{\on{Subcase\ 2.3}}$ \hskip.02in: $s_{ij,q} = 0$ for any $i \geq
u$ and $|[S_q]| < p^{e_u}$.

\vskip.1in

Note that $0 < |[S_q]|$ by the case assumption.

In this subcase, apply $\partial_{X^{[S_q]}}$ to $(\heartsuit)$ and obtain

\begin{eqnarray*}
\lefteqn{\partial_{X^{[S_q]}}
\left(w^{p^{e_i}} + \sum_{|(S,T)| =l}\alpha_{ST}X^SY^T\right) 
= \partial_{X^{[S_q]}}\left(\sum_{|(S,T)|=l}\alpha_{ST}X^SY^T\right) 
}\\
&=& \alpha_{ST}X^{S_r}Y^T + (\text{other\ 
monomials\ of\ degree\ }(l -|[S_q]|)) \\
&\in& \left((\bI_P)_{p^{e_u} - |[S_q]|} + 
{\fm_P}^{l - |[S_q]| + 1}\right) \cap
{\fm_P}^{l - |[S_q]|} \\
&=& (\bI_P)_{p^{e_u} - |[S_q]|} \cap {\fm_P}^{l - |[S_q]|} + {\fm_P}^{l -
|[S_q]| + 1}. \\
\end{eqnarray*}
On the other hand, we observe
$$(\bI_P)_{p^{e_u} - |[S_q]|} \cap {\fm_P}^{l - |[S_q]|} \subset \sum_{1
\leq i \leq M}h_{ij}{\fm_P}^{l - |[S_q]| - p^{e_i}}.$$
(We use the convention that ${\fm_P}^{l - |[S_q]| - p^{e_i}} = R_P$ if $l -
|[S_q]| - p^{e_i} \leq 0$.  The inclusion follows from the same argument as
in Subcase 1.2.)

However, this inclusion implies that any monomial of degree $l - |[S_q]|$ in
the power series expansion of an element in $(\bI_P)_{p^{e_u} - |[S_q]|}
\cap {\fm_P}^{l - |[S_q]|}$, with respect to
the regular system of parameters $(x_1, \cdots, x_N, y_{N+1}, \cdots, y_d)$
should be divisible by some element in the set 
$\left\{x_{ij}^{p^{e_i}} \mid 1\leq i \leq M\right\}$, and hence that the
monomial $X^{S_r}Y^T$ can not appear.

Therefore, in this subcase, we conclude
$$\alpha_{ST} = 0.$$
\end{case}

{}From the above analysis of the terms in $\sum_{|(S,T)| =
l}\alpha_{ST}X^SY^T$, it follows that
$$f \in \sum_{|(S,T)| = l}\alpha_{ST}X^SY^T + J_{l+1} = J_{l+1}.$$
Since $f \in (\bI_P)_{p^{e_u}}$ is arbitrary, we conclude $(\bI_P)_{p^{e_u}}
\subset J_{l+1}$, completing the inductive argument for claim
$(\diamondsuit)$.
\end{step}

\vskip.1in

\begin{step}{Finishing argument for the inductive construction.}

Claim $(\diamondsuit)$ states
$$\begin{array}{ll}
(\bI_P)_{p^{e_u}} &\subset J_l = F^{e_u}(\fm_P) + \left(X^{[C]} \mid C =
(c_{ij} \mid e_i < e_u), |[C]| \geq p^{e_u}\right) \\
& \hskip1in + (\bI_P)_{p^{e_u}} \cap {\fm_P}^{p^{e_u} + 1} + {\fm_P}^l
\hskip1in
\forall l \in {\bZ}_{\geq 0}.\\
\end{array}$$

This implies
$$\begin{array}{ll}
(\bI_P)_{p^{e_u}} &\subset F^{e_u}(\fm_P) + \left(X^{[C]} \mid C = (c_{ij}
\mid e_i < e_u), |[C]| \geq p^{e_u}\right) \\
&\hskip.7in + (\bI_P)_{p^{e_u}} \cap {\fm_P}^{p^{e_u} + 1} +
F^{e_u}({\fm_P}^l)R_P
\hskip.8in
\forall l \in {\bZ}_{\geq 0}.\\
\end{array}$$

Since $R_P$ is a finite $F^{e_u}(R_P)$-module, including
$$F^{e_u}(\fm_P) + \left(X^{[C]} \mid C = (c_{ij} \mid e_i < e_u), |[C]|
\geq p^{e_u}\right) + (\bI_P)_{p^{e_u}} \cap {\fm_P}^{p^{e_u} + 1}$$
as an $F^{e_u}(R_P)$-submodule, we conclude (cf. \cite{MR879273}, page 62
last line) that

$$\begin{array}{ll}
(\bI_P)_{p^{e_u}} &\subset \bigcap_{l \in {\bZ}_{\geq
0}}\left[F^{e_u}(\fm_P) + \left(X^{[C]} \mid C = (c_{ij} \mid e_i < e_u),
|[C]| \geq p^{e_u}\right) \right.\\
&\hskip1.9in \left.+ (\bI_P)_{p^{e_u}} \cap
{\fm_P}^{p^{e_u} + 1} + F^{e_u}({\fm_P}^l)R_P\right]\\
&= F^{e_u}(\fm_P) + \left(X^{[C]} \mid C = (c_{ij} \mid e_i < e_u), |[C]|
\geq p^{e_u}\right) + (\bI_P)_{p^{e_u}} \cap {\fm_P}^{p^{e_u} + 1}.\\
\end{array}$$
Since
$$\left(X^{[C]} \mid C = (c_{ij} \mid e_i < e_u), |[C]| \geq p^{e_u}\right)
\subset (\bI_P)_{p^{e_u}},$$
we also conclude
$$(\bI_P)_{p^{e_u}} \subset F^{e_u}(\fm_P) \cap (\bI_P)_{p^{e_u}} +
\left(X^{[C]} \mid C = (c_{ij} \mid e_i < e_u), |[C]| \geq p^{e_u}\right) +
(\bI_P)_{p^{e_u}} \cap {\fm_P}^{p^{e_u} + 1}.$$
Now choose
$$\left\{h_{uj}' = {x_{uj}'}^{p^{e_u}}\right\} \subset F^{e_u}(\fm_P) \cap
(\bI_P)_{p^{e_u}}$$
such that
$$\left\{h_{uj}' \bmod  \fm_P^{p^{e_u} + 1}\right\} \cup
\left\{X^{[C]} \bmod  \fm_P^{p^{e_u} + 1} \mid C =
(c_{ij} \mid e_i < e_u), |[C]| = p^{e_u}\right\}$$
form a $k$-basis of $L(\bI_P)_{p^{e_u}}$.

In order to finish the inductive argument (cf. Step 2) 
to complete the proof for assertion (1), we have only to 
replace $\left\{h_{ij}\right\}$ and
$\left\{x_{ij}\right\}$ with 
$\left\{h_{ij}'\right\}$ and $\left\{x_{ij}'\right\}$.
\end{step}

\end{proof}

\begin{proof}[\underline{\it Proof for assertion $\on{(2)}$}] Take a regular
system of parameters $(x_1, \cdots, x_N, y_{N+1}, \cdots, y_d)$ and an LGS
$\bH$ of $\bI_P$ as described in assertion
(1).

By choosing an affine neighborhood $P \in U_P = \on{Spec}\hskip.02in R_r$ of
$P$ sufficiently small, we may assume that $(x_1, \cdots, x_N, y_{N+1},
\cdots, y_d)$ is a regular
system of parameters over $U_P$ and that $\bH =
\{(x_l^{p^{e_l}},p^{e_l})\}_{l = 1}^N \subset \bI_r$.

Now we know by assumption that $\bI$ is of r.f.g. type, 
i.e., $\bI =
G_R(\left\{(f_{\lambda},a_{\lambda})\right\}_{\lambda \in \Lambda})$
for some $\left\{(f_{\lambda},a_{\lambda})\right\}_{\lambda \in 
\Lambda} \subset R
\times {\bQ}_{\geq 0}$ with $\# \Lambda < \infty$.

Since $\bI_P = G_{R_P}(\bH)$, we can write each $f_{\lambda}$ as a finite
sum of the form $\sum g_{B,\lambda}H^B$ with $g_{B,\lambda} \in R_P$.  By
shrinking $U_P = \on{Spec}\hskip.02in R_r$ if
necessary, we may assume that the coefficients $g_{B,\lambda}$ are in $R_r$
for all $B$ and $\lambda \in \Lambda$.  Then we have
$$\bI_r = G_{R_r}(\left\{(f_{\lambda},a_{\lambda})\right\}_{\lambda \in 
\Lambda})\subset G_{R_r}(\bH).$$
Since the opposite inclusion $\bI_r \supset G_{R_r}(\bH)$ is obvious, we
conclude
$$\bI_r = G_{R_r}(\bH).$$

It follows immediately from the above conclusions that
$$\begin{array}{ll}
\on{Supp}(\bI) \cap U_P &= \on{Supp}(\bI_r) = \on{Supp}(G_{R_r}(\bH)) \\
&= \left\{Q \in U_P \mid \mu_Q(x_l^{p^{e_l}},p^{e_l}) \geq 1 
\text{\ for\ } l =1, \cdots, N\right\} \\ 
&= V(x_1,\cdots, x_N),\\
\end{array}$$
which is nonsingular.

Given any closed point $Q \in \on{Supp}(\bI) \cap U_P = V(x_1, \cdots,
x_N)$, it also follows from the above conclusions that $(x_1, \cdots, x_N)$
is a part of a regular system of parameters at
$Q$ with a subset $\bH = \left\{(x_l^{p^{e_l}},p^{e_l})
\right\}_{l = 1}^N \subset\bI_Q$ such that $\bI_Q = G_{R_Q}(\bH)$.  
This implies that $\bH$ is an LGS
of $\bI_Q$ and that $\widetilde{\mu}(Q) = \infty$.  Therefore, we conclude
$$(\sigma(Q), \widetilde{\mu}(Q)) = (\sigma(P), \infty).$$
This concludes the proof of Theorem \ref{A.1.1.1}.
\end{proof}

\end{subsection}

\end{section}

\end{chapter} 

\chapter*\bibname
\begin{quote}{\small 
The list of references for Part II is largely identical 
to the one for Part~I 
\linebreak[4]
{\cite{MR2361797}}, which we reproduce below for the convenience 
of the reader, with a few more papers and 
Part I added.
}\end{quote}

\medskip

\begin{center}{\it\noindent
References for Part I 
}\end{center}

\medskip

\begin{bibpart}{HLOQ00}
\bibitem[Abh66]{MR0217069}
Shreeram~S. Abhyankar.
\newblock {\em Resolution of singularities of embedded algebraic surfaces}.
\newblock Pure and Applied Mathematics, Vol. 24. 
 Academic Press, New York, 1966.

\bibitem[Abh77]{MR542446}
Shreeram~S. Abhyankar.
\newblock {\em Lectures on expansion techniques in algebraic geometry},
  volume~57 of {\em Tata Institute of Fundamental Research Lectures on
  Mathematics and Physics}.
\newblock Tata Institute of Fundamental Research, Bombay, 1977.
\newblock Notes by Balwant Singh.

\bibitem[Abh83]{MR713043}
Shreeram~S. Abhyankar.
\newblock Desingularization of plane curves.
\newblock In {\em Singularities, Part 1 (Arcata, Calif., 1981)}, volume~40 of
  {\em Proc. Sympos. Pure Math.}, pages 1--45. Amer. Math. Soc., Providence,
  RI, 1983.

\bibitem[AdJ97]{MR1487237}
Dan Abramovich and A.~J. de~Jong.
\newblock Smoothness, semistability, and toroidal geometry.
\newblock {\em J. Algebraic Geom.}, 6(4):789--801, 1997.

\bibitem[AHV75]{MR0444999}
Jos{\'e}~M. Aroca, Heisuke Hironaka, and Jos{\'e}~L. Vicente.
\newblock {\em The theory of the maximal contact}.
\newblock Instituto ``Jorge Juan'' de Matem\'aticas, Consejo Superior de
  Investigaciones Cientificas, Madrid, 1975.
\newblock Memorias de Matem\'atica del Instituto ``Jorge Juan'', No.~29.
  [Mathematical Memoirs of the ``Jorge Juan'' Institute, No.~29].

\bibitem[AHV77]{MR480502}
Jos{\'e}~M. Aroca, Heisuke Hironaka, and Jos{\'e}~L. Vicente.
\newblock {\em Desingularization theorems}, volume~30 of {\em Memorias de
  Matem\'atica del Instituto ``Jorge Juan'' [Mathematical Memoirs of the Jorge
  Juan Institute]}.
\newblock Consejo Superior de Investigaciones Cient\'\i ficas, Madrid, 1977.

\bibitem[Ben70]{MR0252388}
Bruce~M. Bennett.
\newblock On the characteristic functions of a local ring.
\newblock {\em Ann. of Math. (2)}, 91:25--87, 1970.

\bibitem[BEV05]{MR2174912}
Ana~M. Bravo, Santiago Encinas, and Orlando E.~Villamayor.
\newblock A simplified proof of desingularization and applications.
\newblock {\em Rev. Mat. Iberoamericana}, 21(2):349--458, 2005.

\bibitem[Bie04]{B_BIRS}
Edward Bierstone.
\newblock Resolution of singularities.
\newblock {\em A preprint for series of lectures at the ``Workshop on
  resolution of singularities, factorization of birational mappings,
  and toroidal geometry'' at Banff International Research Station,
  Banff, during 11-16 December}, 2004.

\bibitem[Bie05]{B_Harvard}
Edward Bierstone.
\newblock Functoriality in resolution of singularities.
\newblock {\em Slides for a colloquium at Harvard University on November 10},
  2005.

\bibitem[BM89]{MR1001853}
Edward Bierstone and Pierre~D. Milman.
\newblock Uniformization of analytic spaces.
\newblock {\em J. Amer. Math. Soc.}, 2(4):801--836, 1989.

\bibitem[BM90]{MR1051203}
Edward Bierstone and Pierre~D. Milman.
\newblock Local resolution of singularities.
\newblock In {\em Real analytic and algebraic geometry (Trento, 1988)}, volume
  1420 of {\em Lecture Notes in Math.}, pages 42--64. Springer, Berlin, 1990.

\bibitem[BM91]{MR1106412}
Edward Bierstone and Pierre~D. Milman.
\newblock A simple constructive proof of canonical resolution of singularities.
\newblock In {\em Effective methods in algebraic geometry (Castiglioncello,
  1990)}, volume~94 of {\em Progr. Math.}, pages 11--30. Birkh\"auser Boston,
  Boston, MA, 1991.

\bibitem[BM97]{MR1440306}
Edward Bierstone and Pierre~D. Milman.
\newblock Canonical desingularization in characteristic zero by blowing up the
  maximum strata of a local invariant.
\newblock {\em Invent. Math.}, 128(2):207--302, 1997.

\bibitem[BM99]{MR1748600}
Edward Bierstone and Pierre~D. Milman.
\newblock Resolution of singularities.
\newblock In {\em Several complex variables (Berkeley, CA, 1995--1996)},
  volume~37 of {\em Math. Sci. Res. Inst. Publ.}, pages 43--78. Cambridge Univ.
  Press, Cambridge, 1999.

\bibitem[BM03]{MR2078560}
Edward Bierstone and Pierre~D. Milman.
\newblock Desingularization algorithms. {I}. {R}ole of exceptional divisors.
\newblock {\em Mosc. Math. J.}, 3(3):751--805, 1197, 2003.
\newblock \{Dedicated to Vladimir Igorevich Arnold on the occasion of his 65th
  birthday\}.

\bibitem[BM07]{AG0702375}
Edward Bierstone and Pierre~D. Milman.
\newblock Functoriality in resolution of singularities.
\newblock {\em preprint, {\tt http://ar{X}iv.org/abs/math.AG/0702375}}, 2007.

\bibitem[Bou64]{MR0194450}
Nicolas Bourbaki.
\newblock {\em \'{E}l\'ements de math\'ematique. {F}asc. {XXX}. {A}lg\`ebre
  commutative. {C}hapitre 5: {E}ntiers. {C}hapitre 6: {V}aluations}.
\newblock Actualit\'es Scientifiques et Industrielles, No.~1308. Hermann,
  Paris, 1964.

\bibitem[BP96]{MR1397679}
Fedor~A. Bogomolov and Tony~G. Pantev.
\newblock Weak {H}ironaka theorem.
\newblock {\em Math. Res. Lett.}, 3(3):299--307, 1996.

\bibitem[BV03]{MR1971154}
Ana~M. Bravo and Orlando E.~Villamayor.
\newblock A strengthening of resolution of singularities in characteristic
  zero.
\newblock {\em Proc. London Math. Soc. (3)}, 86(2):327--357, 2003.

\bibitem[Cos87]{MR907903}
Vincent Cossart.
\newblock Forme normale d'une fonction sur un {$k$}-sch\'ema de dimension {$3$}
  et de caract\'eristique positive.
\newblock In {\em G\'eom\'etrie alg\'ebrique et applications, I (La R\'abida,
  1984)}, volume~22 of {\em Travaux en Cours}, pages 1--21. Hermann, Paris,
  1987.

\bibitem[Cos00]{MR1748622}
Vincent Cossart.
\newblock Uniformisation et d\'esingularisation
  des surfaces d'apr\`es {Z}ariski.
\newblock In {\em Resolution of singularities (Obergurgl, 1997)},
  volume 181 of {\em Progr. Math.}, pages 239--258. Birkh\"auser, Basel, 2000.

\bibitem[Cos04]{C_BIRS}
Vincent Cossart.
\newblock Towards local uniformization along a valuation in Artin-Schreier
  extensions (dimension $3$).
\newblock {\em A talk at the ``Workshop on resolution of singularities,
  factorization of birational mappings, and toroidal geometry'' at Banff
  International Research Station, Banff, during 11-16 December}, 2004.

\bibitem[Cut04]{MR2058431}
Steven~Dale Cutkosky.
\newblock {\em Resolution of singularities.}, volume~63 of {\em Graduate 
 Studies in Mathematics}.
\newblock American Mathematical Society, Providence, RI, 2004.

\bibitem[Cut06]{AG0606530}
Steven~Dale Cutkosky.
\newblock Resolution of singularities for $3$-folds in positive characteristic.
\newblock {\em preprint, {\tt http://ar{X}iv.org/abs/math.AG/0606530}}, 2006.

\bibitem[dJ96]{MR1423020}
A.~J. de~Jong.
\newblock Smoothness, semi-stability and alterations.
\newblock {\em Inst. Hautes \'Etudes Sci. Publ. Math.}, No.~83:51--93, 1996.

\bibitem[EH02]{MR1949115}
Santiago Encinas and Herwig Hauser.
\newblock Strong resolution of singularities in characteristic zero.
\newblock {\em Comment. Math. Helv.}, 77(4):821--845, 2002.

\bibitem[ENV03]{MR1974392}
Santiago Encinas, A.~Nobile, and Orlando E.~Villamayor.
\newblock On algorithmic equi-resolution and stratification of {H}ilbert
  schemes.
\newblock {\em Proc. London Math. Soc. (3)}, 86(3):607--648, 2003.

\bibitem[EV98]{MR1654779}
Santiago Encinas and Orlando E.~Villamayor.
\newblock Good points and constructive resolution of singularities.
\newblock {\em Acta Math.}, 181(1):109--158, 1998.

\bibitem[EV00]{MR1748620}
Santiago Encinas and Orlando E.~Villamayor.
\newblock A course on constructive desingularization and equivariance.
\newblock In {\em Resolution of singularities (Obergurgl, 1997)}, volume 181 of
  {\em Progr. Math.}, pages 147--227. Birkh\"auser, Basel, 2000.

\bibitem[EV03]{MR2023188}
Santiago Encinas and Orlando E.~Villamayor.
\newblock A new proof of desingularization over fields of characteristic zero.
\newblock In {\em Proceedings of the International Conference on Algebraic
  Geometry and Singularities (Spanish) (Sevilla, 2001)}, volume~19, pages
  339--353, 2003.

\bibitem[Gir74]{MR0460712}
Jean Giraud.
\newblock Sur la th\'eorie du contact maximal.
\newblock {\em Math. Z.}, 137:285--310, 1974.

\bibitem[Gir75]{MR0384799}
Jean Giraud.
\newblock Contact maximal en caract\'eristique positive.
\newblock {\em Ann. Sci. \'Ecole Norm. Sup. (4)}, 8(2):201--234, 1975.

\bibitem[Gir83]{MR734215}
Jean Giraud.
\newblock Forme normale d'une fonction sur une surface de caract\'eristique
  positive.
\newblock {\em Bull. Soc. Math. France}, 111(2):109--124, 1983.

\bibitem[Gro67]{MR0238860}
Alexander Grothendieck.
\newblock \'{E}l\'ements de g\'eom\'etrie alg\'ebrique. {IV}. \'{E}tude locale
  des sch\'emas et des morphismes de sch\'emas {IV}.
\newblock {\em Inst. Hautes \'Etudes Sci. Publ. Math.}, No.~32, 1967.

\bibitem[Hau98]{MR1652479}
Herwig Hauser.
\newblock Seventeen obstacles for resolution of singularities.
\newblock In {\em Singularities (Oberwolfach, 1996)}, volume 162 of {\em Progr.
  Math.}, pages 289--313. Birkh\"auser, Basel, 1998.

\bibitem[Hau03]{MR1978567}
Herwig Hauser.
\newblock The {H}ironaka theorem on resolution of singularities (or: {A} proof
  we always wanted to understand).
\newblock {\em Bull. Amer. Math. Soc. (N.S.)}, 40(3):323--403 (electronic),
  2003.

\bibitem[Hir63]{MR0175898}
Heisuke Hironaka.
\newblock On resolution of singularities (characteristic zero).
\newblock In {\em Proc. Internat. Congr. Mathematicians (Stockholm, 1962)},
  pages 507--521. Inst. Mittag-Leffler, Djursholm, 1963.

\bibitem[Hir64]{MR0199184}
Heisuke Hironaka.
\newblock Resolution of singularities of an algebraic variety over a field of
  characteristic zero. {I}, {II}.
\newblock {\em Ann. of Math. (2) 79 (1964), 109--203; ibid. (2)}, 79:205--326,
  1964.

\bibitem[Hir70]{MR0269658}
Heisuke Hironaka.
\newblock Additive groups associated with points of a projective space.
\newblock {\em Ann. of Math. (2)}, 92:327--334, 1970.

\bibitem[Hir72a]{MR0393555}
Heisuke Hironaka.
\newblock Gardening of infinitely near singularities.
\newblock In {\em Algebraic geometry, Oslo 1970 (Proc. Fifth Nordic Summer
  School in Math.)}, pages 315--332. Wolters-Noordhoff, Groningen, 1972.

\bibitem[Hir72b]{MR0393028}
Heisuke Hironaka.
\newblock Schemes, etc.
\newblock In {\em Algebraic geometry, Oslo 1970 (Proc. Fifth Nordic Summer
  School in Math.)}, pages 291--313. Wolters-Noordhoff, Groningen, 1972.

\bibitem[Hir77]{MR0498562}
Heisuke Hironaka.
\newblock Idealistic exponents of singularity.
\newblock In {\em Algebraic geometry (J. J. Sylvester Sympos., Johns Hopkins
  Univ., Baltimore, Md., 1976)}, pages 52--125. Johns Hopkins Univ. Press,
  Baltimore, Md., 1977.

\bibitem[Hir03]{MR1996845}
Heisuke Hironaka.
\newblock Theory of infinitely near singular points.
\newblock {\em J. Korean Math. Soc.}, 40(5):901--920, 2003.

\bibitem[Hir05]{MR2145950}
Heisuke Hironaka.
\newblock Three key theorems on infinitely near singularities.
\newblock In {\em Singularit\'es Franco-Japonaises}, volume~10 of {\em S\'emin.
  Congr.}, pages 87--126. Soc. Math. France, Paris, 2005.

\bibitem[Hir06]{H-Trieste}
Heisuke Hironaka.
\newblock A program for resolution of singularities, in all characteristics
  $p>0$ and in all dimensions.
\newblock {\em preprint for series of lectures in ``Summer School on Resolution
  of Singularities'' at International Center for Theoretical Physics, Trieste,
  during 12-30 June}, 2006.

\bibitem[HLOQ00]{MR1748614}
Herwig Hauser, Joseph Lipman, Frans Oort, and Adolfo Quir{\'o}s, editors.
\newblock {\em Resolution of singularities}, volume 181 of {\em Progress in
  Mathematics}.
\newblock Birkh\"auser Verlag, Basel, 2000.
\newblock A research textbook in tribute to Oscar Zariski, Papers from the
  Working Week held in Obergurgl, September 7--14, 1997.


\bibitem[Kol07]{MR2289519}
J{\'a}nos Koll{\'a}r.
\newblock {\em Lectures on resolution of singularities}, volume 166 of {\em
  Annals of Mathematics Studies}.
\newblock Princeton University Press, Princeton, NJ, 2007.

\bibitem[Kuh97]{Kuhlmann97}
Franz-Viktor Kuhlmann.
\newblock On local uniformization in arbitrary characteristic.
\newblock {\em The Fields Institute Preprint Series}, 1997.

\bibitem[Kuh00]{MR1748629}
Franz-Viktor Kuhlmann.
\newblock Valuation theoretic and model theoretic aspects of local
  uniformization.
\newblock In {\em Resolution of singularities (Obergurgl, 1997)}, volume
 181 of
  {\em Progr. Math.}, pages 381--456. Birkh\"auser, Basel, 2000.

\bibitem[Lip69]{MR0276239}
Joseph Lipman.
\newblock Rational singularities, with applications to algebraic surfaces and
  unique factorization.
\newblock {\em Inst. Hautes \'Etudes Sci. Publ. Math.}, No.~36:195--279, 1969.

\bibitem[Lip75]{MR0389901}
Joseph Lipman.
\newblock Introduction to resolution of singularities.
\newblock In {\em Algebraic geometry (Proc. Sympos. Pure Math., Vol. 29,
  Humboldt State Univ., Arcata, Calif., 1974)}, pages 187--230. Amer. Math.
  Soc., Providence, R.I., 1975.

\bibitem[Lip78]{MR0491722}
Joseph Lipman.
\newblock Desingularization of two-dimensional schemes.
\newblock {\em Ann. Math. (2)}, 107(1):151--207, 1978.

\bibitem[Lip83]{MR713245}
Joseph Lipman.
\newblock Quasi-ordinary singularities of surfaces in {${\bf C}\sp{3}$}.
\newblock In {\em Singularities, Part 2 (Arcata, Calif., 1981)}, volume~40 of
  {\em Proc. Sympos. Pure Math.}, pages 161--172. Amer. Math. Soc., Providence,
  RI, 1983.

\bibitem[LT74]{LT74}
Monique Lejeune-Jalabert and Bernard Teissier.
\newblock Cl\^oture int\'egrale des id\'eaux et \'equisingularit\'e.
\newblock {\em S\'eminaire sur les singularit\'e \`a l'Ecole Polytechnique},
  pages 1--66, 1974--1975.

\bibitem[Mat86]{MR879273}
Hideyuki Matsumura.
\newblock {\em Commutative ring theory}, volume~8 of {\em Cambridge Studies in
  Advanced Mathematics}.
\newblock Cambridge University Press, Cambridge, 1986.
\newblock Translated from the Japanese by M. Reid.

\bibitem[Mk07]{AG0103120}
Kenji Matsuki.
\newblock Resolution of singularities in characteristic zero with focus on the
  inductive algorithm by {V}illamayor and its simplification by
  {W}{\l}odarczyk.
\newblock {\em preprint, formally a revision of\ \ 
{\tt http://ar{X}iv.org/abs/math.AG/0103120}, but
completely rewritten from scratch}, 2007.

\bibitem[Moh92]{MR1186705}
T.~T. Moh.
\newblock Quasi-canonical uniformization of hypersurface singularities of
  characteristic zero.
\newblock {\em Comm. Algebra}, 20(11):3207--3249, 1992.

\bibitem[Moh96]{MR1395176}
T.~T. Moh.
\newblock On a {N}ewton polygon approach to the uniformization of singularities
  of characteristic {$p$}.
\newblock In {\em Algebraic geometry and singularities (La R\'abida, 1991)},
  volume 134 of {\em Progr. Math.}, pages 49--93. Birkh\"auser, Basel, 1996.

\bibitem[Nag57]{MR0089836}
Masayoshi Nagata.
\newblock Note on a paper of {S}amuel concerning asymptotic properties of
  ideals.
\newblock {\em Mem. Coll. Sci. Univ. Kyoto. Ser. A. Math.}, 30:165--175, 1957.

\bibitem[Nar83a]{MR684627}
Ramanujachari Narasimhan.
\newblock Hyperplanarity of the equimultiple locus.
\newblock {\em Proc. Amer. Math. Soc.}, 87(3):403--408, 1983.

\bibitem[Nar83b]{MR715853}
Ramanujachari Narasimhan.
\newblock Monomial equimultiple curves in positive characteristic.
\newblock {\em Proc. Amer. Math. Soc.}, 89(3):402--406, 1983.

\bibitem[Oda73]{MR0472824}
Tadao Oda.
\newblock Hironaka's additive group scheme.
\newblock In {\em Number theory, algebraic geometry and commutative algebra, in
  honor of Yasuo Akizuki}, pages 181--219. Kinokuniya, Tokyo, 1973.

\bibitem[Oda83]{MR723465}
Tadao Oda.
\newblock Hironaka's additive group scheme. {II}.
\newblock {\em Publ. Res. Inst. Math. Sci.}, 19(3):1163--1179, 1983.

\bibitem[Oda87]{MR894302}
Tadao Oda.
\newblock Infinitely very near singular points.
\newblock In {\em Complex analytic singularities}, volume~8 of {\em Adv. Stud.
  Pure Math.}, pages 363--404. North-Holland, Amsterdam, 1987.

\bibitem[Par99]{MR1714830}
Kapil~H. Paranjape.
\newblock The {B}ogomolov-{P}antev resolution, an expository account.
\newblock In {\em New trends in algebraic geometry (Warwick, 1996)}, volume 264
  of {\em London Math. Soc. Lecture Note Ser.}, pages 347--358. Cambridge Univ.
  Press, Cambridge, 1999.

\bibitem[Pil04]{P_BIRS}
Olivier Piltant.
\newblock Applications of ramification theory to resolution of
  three-dimensional varieties.
\newblock {\em A talk at the ``Workshop on resolution of singularities,
  factorization of birational mappings, and toroidal geometry'' at Banff
  International Research Station, Banff, during 11-16 December}, 2004.

\bibitem[Pom74]{MR0441971}
Marie-Jeanne Pomerol.
\newblock Sur la strate de {S}amuel du sommet d'un c\^one en caract\'eristique
  positive.
\newblock {\em Bull. Sci. Math. (2)}, 98(3):173--182, 1974.

\bibitem[Spi04]{S_BIRS}
Mark Spivakovsky.
\newblock Puiseaux expansions, a local analogue of Nash's space
  of arcs and the local uniformization theorem.
\newblock {\em A talk at the ``Workshop on resolution of singularities,
  factorization of birational mappings, and toroidal geometry'' at Banff
  International Research Station, Banff, during 11-16 December}, 2004.

\bibitem[Tei03]{MR2018565}
Bernard Teissier.
\newblock Valuations, deformations, and toric geometry.
\newblock In {\em Valuation theory and its applications, Vol. II (Saskatoon,
  SK, 1999)}, volume~33 of {\em Fields Inst. Commun.}, pages 361--459. Amer.
  Math. Soc., Providence, RI, 2003.

\bibitem[Vil89]{MR985852}
Orlando E.~Villamayor.
\newblock Constructiveness of {H}ironaka's resolution.
\newblock {\em Ann. Sci. \'Ecole Norm. Sup. (4)}, 22(1):1--32, 1989.

\bibitem[Vil92]{MR1198092}
Orlando E.~Villamayor.
\newblock Patching local uniformizations.
\newblock {\em Ann. Sci. \'Ecole Norm. Sup. (4)}, 25(6):629--677, 1992.

\bibitem[W{\l}o05]{MR2163383}
Jaros{\l}aw W{\l}odarczyk.
\newblock Simple {H}ironaka resolution in characteristic zero.
\newblock {\em J. Amer. Math. Soc.}, 18(4):779--822 (electronic), 2005.

\bibitem[Zar40]{MR0002864}
Oscar Zariski.
\newblock Local uniformization on algebraic varieties.
\newblock {\em Ann. of Math. (2)}, 41:852--896, 1940.
\end{bibpart}

\bigskip

\begin{center}{\it\noindent
References added for Part II
}\end{center}

\bigskip

\begin{bibpart}{HLOQ00}
\bibitem[BV08]{AG08074308}
Ana~M. Bravo and Orlando E.~Villamayor.
\newblock Hypersurface singularities in positive characteristic
and stratification of the singular locus.
\newblock {\em preprint, {\tt http://ar{X}iv.org/abs/0807.4308}}, 2008.
\bibitem[CP07]{CP0703302}
Vincent Cossart and Olivier Piltant.
\newblock Resolution of singularities of threefolds 
in positive characteristic II,\ 
\newblock {\em preprint, {\tt 
http://hal.archives-ouvertes.fr/hal-00139445
}}, 2007.
\bibitem[CP08]{MR2427629}
Vincent Cossart and Olivier Piltant.
\newblock Resolution of singularities of threefolds in positive characteristic.
  {I}. {R}eduction to local uniformization on {A}rtin-{S}chreier and purely
  inseparable coverings.
\newblock {\em J. Algebra}, 320(3):1051--1082, 2008.
\bibitem[EV07]{AG0702836}
Santiago Encinas and Orlando E.~Villamayor.
\newblock Rees algebras and resolution of singularities.
\newblock {\em preprint, {\tt http://ar{X}iv.org/abs/math.AG/0702836}},
 2007.
\bibitem[Kaw07]{MR2361797}
Hiraku Kawanoue.
\newblock Toward resolution of singularities over a field of positive
  characteristic. {I}. {F}oundation; the language of the idealistic filtration.
\newblock {\em Publ. Res. Inst. Math. Sci.}, 43(3):819--909, 2007.
\bibitem[RIMS08]{RIMS2008}
Workshop: 
\newblock ``On the Resolution of Singularities''. 
\quad
\newblock {\em Slides for workshop 
at Research Institute for Mathematical Sciences, Kyoto, 
during 1-5 December}, 2008.\\
{{\tt http://www.kurims.kyoto-u.ac.jp/$\sim$kenkyubu/proj08-mori/index.html}}

\bibitem[Vil06a]{AG0606796}
Orlando E.~Villamayor.
\newblock Rees algebras on smooth schemes: integral closure 
and higher differential operators.
\newblock {\em preprint, {\tt http://ar{X}iv.org/abs/math.AC/0606795}},
 2006.
\bibitem[Vil06b]{AC0606795}
Orlando E.~Villamayor.
\newblock Hypersurface singularities in positive characteristic.
\par
\newblock {\em preprint,\quad
{\tt http://ar{X}iv.org/abs/math.AG/0606796}}, 2006.
\bibitem[W{\l}o07]{W07}
Jaros{\l}aw W{\l}odarczyk.
\newblock Giraud maximal contact in positive characteristic.
\newblock {\em in preparation}, 2007.
\end{bibpart}

\end{document}